\documentclass{amsart}

\usepackage{amssymb,amsmath,amsthm,latexsym,amscd}
\usepackage{epsfig}
\usepackage{psfrag}
\usepackage{graphicx}
\usepackage{tikz}
\usepackage{float}

\newtheorem{Theorem}{\bf Theorem}
\newtheorem{lemma}[Theorem]{\bf Lemma}
\newtheorem{proposition}[Theorem]{\bf Proposition}
\newtheorem{corollary}[Theorem]{\bf Corollary}

\newtheorem{definition}[Theorem]{\bf Definition}

\newtheorem{remark}[Theorem]{\bf Remark}
\newtheorem{theorem}[Theorem]{\bf Theorem}

\def\scfig #1 #2 {\resizebox{#2}{!}{\includegraphics{#1}}}

\newcommand{\be}{\begin{equation}}
\newcommand{\ee}{\end{equation}}

\def\hpic #1 #2 {\mbox{$\begin{array}[c]{l} 
\epsfig{file=#1,height=#2}\end{array}$}}
\def\wpic #1 #2 {\mbox{$\begin{array}[c]{l} 
\epsfig{file=#1,width=#2}\end{array}$}}

\begin{document}
\title{Quantum double inclusions associated to a family of Kac algebra subfactors}
\author{Sandipan De}
\address{Stat-Math Unit\\Indian Statistical Institute, 8th Mile, Mysore Road\\ Bangalore-560059}
\email{sandipan$\_$vs@isibang.ac.in}


\keywords{Subfactors, Hopf algebras, Planar algebras}
\subjclass[2010]{46L37, 16S40, 16T05}
\maketitle
\begin{abstract} 
In \cite{Sde2018} we defined the notion of \textit{quantum double inclusion} associated to a finite-index and finite-depth subfactor and
studied the quantum double inclusion associated to the Kac algebra subfactor $R^H \subset R$ where $H$ is a finite-dimensional Kac algebra
acting outerly on the hyperfinite $II_1$ factor $R$ and $R^H$ denotes the fixed-point subalgebra.  
  In this article we analyse
 quantum double inclusions associated to the family of Kac algebra subfactors given by
 $\{ R^H \subset R \rtimes \underbrace{H \rtimes H^* \rtimes \cdots}_{{\text{$m$ times}}}  : m \geq 1 \}$. For each $m > 2$, we construct  
 a model $\mathcal{N}^m \subset \mathcal{M}$ for the quantum double inclusion of
 $\{ R^H \subset R \rtimes \underbrace{H \rtimes H^* \rtimes \cdots}_{{\text{$m-2$ times}}} \}$ with $\mathcal{N}^m = 
 ((\cdots \rtimes H^{-2} \rtimes H^{-1}) \otimes (H^m \rtimes H^{m+1} \cdots))^{\prime \prime}, \mathcal{M} = (\cdots \rtimes H^{-1} \rtimes H^0
 \rtimes H^1 \rtimes \cdots)^{\prime \prime}$ and where for any integer $i$, $H^i$ denotes $H$ or $H^*$ according as $i$ is odd or even.
 In this article, we give an explicit 
 description of $P^{\mathcal{N}^m \subset \mathcal{M}}$ ($m > 2$), the subfactor planar algebra associated to 
 $\mathcal{N}^m \subset \mathcal{M}$, which turns out to be a planar subalgebra of $^{*(m)}\!P(H^m)$ 
 (the adjoint of the $m$-cabling of the planar algebra of $H^m$).
 We then show that for $m > 2$, depth of $\mathcal{N}^m \subset \mathcal{M}$ is always two. 
 Observing that $\mathcal{N}^m \subset \mathcal{M}$ is reducible for all $m > 2$, we 
  explicitly
  describe the weak Hopf $C^*$-algebra structure on $(\mathcal{N}^m)^{\prime} \cap \mathcal{M}_2$, thus obtaining a family of 
  weak Hopf $C^*$-algebras starting with a single Kac algebra $H$.
\end{abstract}

%

\section*{Introduction}
The motivation for this article primarily stems from the work of the author in \cite{Sde2018}. Given a finite-index
and finite-depth subfactor $N \subset M$ with $N ( = M_0) \subset M ( = M_1) \subset M_2 \subset M_3 \subset \cdots$ 
being the Jones' basic construction tower associated to $N \subset M$, we defined in \cite{Sde2018} the inclusion 
\begin{align*}
 N \vee (M^{\prime} \cap M_{\infty}) \subset M_{\infty}
\end{align*}
to be the \textit{quantum double inclusion} associated to $N \subset M$ where $M_{\infty}$ denotes the $II_1$ factor obtained as the von Neumann
closure $(\cup_{n = 0}^{\infty} M_n)^{\prime \prime}$ in the GNS representation with respect to the trace on $\cup_{n = 0}^{\infty} M_n$
and $N \vee (M^{\prime} \cap M_{\infty})$ denotes the von Neumann algebra generated by
$N$ and $M^{\prime} \cap M_{\infty}$. In \cite{Sde2018} we studied the quantum double inclusion associated to the Kac algebra subfactor 
$R^H \subset R$ where $H$ is a 
finite-dimensional Kac algebra acting outerly on the hyperfinite $II_1$ factor $R$ and $R^H$ denotes the fixed-point subalgebra. The main result of 
\cite{Sde2018} states that the quantum double inclusion of $R^H \subset R$ is isomorphic to $R \subset R \rtimes D(H)^{cop}$ for some outer action of 
$D(H)^{cop}$ on $R$ where $D(H)$ denotes the Drinfeld double of $H$. This result seemed to be quite interesting and motivated us to analyse 
quantum double inclusions associated to a general class of
Kac algebra subfactors 
given by $\{ R^H \subset R \rtimes \underbrace{H \rtimes H^* \rtimes \cdots}_{{\text{$m$ times}}}  : m \geq 1 \}$.

One of the main steps towards understanding the quantum double inclusions associated to the family of subfactors 
$\{ R^H \subset R \rtimes \underbrace{H \rtimes H^* \rtimes \cdots}_{{\text{$m$ times}}}  : m \geq 1 \}$ is to construct their models.
Given any finite-dimensional Kac algebra $H$, let $H^i$, where $i$ is any integer, denote $H$ or $H^*$ according as $i$ is odd or even.
For each positive integer $m > 2$, we construct in $\S 2$ a hyperfinite, finite-index subfactor $\mathcal{N}^m \subset \mathcal{M}$ where 
$\mathcal{N}^m = ((\cdots \rtimes H^{-3} \rtimes H^{-2} \rtimes 
 H^{-1}) \otimes (H^m \rtimes H^{m+1} \rtimes \cdots))^{\prime \prime}, \ \mathcal{M} = (\cdots \rtimes H^{-1} \rtimes H^0 \rtimes H^1 \rtimes 
 \cdots)^{\prime \prime}$ and show that $\mathcal{N}^m \subset \mathcal{M}$ is a model for the quantum double inclusion of $R^H \subset
 R \rtimes \underbrace{H \rtimes H^* \rtimes \cdots}_{{\text{$m-2$ times}}}$.

 The heart of the paper is $\S 3$ where 
 we compute the basic construction tower associated to $\mathcal{N}^m \subset \mathcal{M}$ and also compute the relative commutants. 
 The proofs all rely on explicit pictorial computations in the planar algebra of $H$ or $H^*$.

 In $\S 4$, we explicitly describe the planar algebra 
 associated to the subfactor $\mathcal{N}^m \subset \mathcal{M}$ ($m > 2$) which turns out to be an interesting planar subalgebra of
 $^{*(m)}\!P(H^m)$ (the adjoint of the $m$-cabling of the planar algebra of $H^m$).

It is evident from the main result of \cite{Sde2018} that the quantum double inclusion of $R^H \subset R$ is of depth two. 
It is thus a natural question to ask whether the quantum double inclusions associated to the family of 
subfactors $\{ R^H \subset
 R \rtimes \underbrace{H \rtimes H^* \rtimes \cdots}_{{\text{$m$ times}}} : m \geq 1 \}$ have finite depth. 
 In this article we answer to this question in affirmative by proving that 
for $m > 2$, depth of $\mathcal{N}^m \subset \mathcal{M}$ is always $2$ (Theorem 11). This is the main result of $\S 5$.
One primary ingredient of the proof is Lemma \ref{element} where we identify the commutant of the middle $H$ in $H^* \rtimes H \rtimes H^*$.

In \cite{Sde2018} we constructed a model $\mathcal{N} \subset \mathcal{M}$ for the quantum double inclusion of $R^H \subset R$. 
As an immediate consequence of the main result \cite[Theorem 40]{Sde2018} one obtains that the relative commutant
$\mathcal{N}^{\prime} \cap \mathcal{M}_2$ is isomorphic to $D(H)^{cop*} (= D(H)^{*op})$ as Kac algebras. In $\S 6$
we explicitly describe the structure maps of $\mathcal{N}^{\prime} \cap \mathcal{M}_2$ which will be useful to achieve a
simple and nice description of the weak Hopf $C^*$-algebra structures on $(\mathcal{N}^m)^{\prime} \cap \mathcal{M}_2$ ($m > 2$) in $\S 7$.

It is well-known (see \cite{NikVnr2000}, \cite{Das2004}) that
if $N \subset M$ is a finite-index reducible depth $2$ inclusion of $II_1$ factors 
and if $N (=M_0) \subset M (= M_1) \subset M_2 \subset M_3 \subset \cdots$ is the Jones' basic construction tower associated to $N \subset M$, 
then the relative commutants
$N^{\prime} \cap M_2$ and $M^{\prime} \cap M_3$ admit mutually dual weak Hopf $C^*$-algebra structures.
Now, for each $m > 2$, the subfactor $\mathcal{N}^m \subset \mathcal{M}$ being reducible and of depth $2$, 
$(\mathcal{N}^m)^{\prime} \cap \mathcal{M}_2$ admits a weak Hopf $C^*$-algebra
structure. The final $\S 7$ is devoted to recovering the weak Hopf $C^*$-algebra structure on 
 $(\mathcal{N}^m)^{\prime} \cap \mathcal{M}_2$ for all $m > 2$.
 Here, in Theorem \ref{WHA},  we construct a family $\{K_m : m > 2\}$ of weak Hopf $C^*$-algebras, with underlying vector spaces 
 $A(H)_{m-2}^{op} \otimes D(H)^{*op} \otimes
 A(H)_{m-2}^{op}$ or 
 $A(H^*)_{m-2}^{op} \otimes D(H)^{*op} \otimes
 A(H)_{m-2}^{op}$ according as $m$ is odd or even, 
 such that $K_m \cong (\mathcal{N}^m)^{\prime} \cap \mathcal{M}_2$ as weak Hopf $C^*$-algebras where, for any positive integer $l$ and 
 any finite-dimensional Kac algebra $K$, $A(K)_l$ denotes 
 the finite crossed product algebra $\underbrace{K \rtimes K^* \rtimes \cdots }_{\text{$l$ times}}$.

 \section{Preliminaries}
 The prerequisites for this article can be found in $\S 1$ and $\S 2$ of \cite{Sde2018}. 
For a convenient reading, below we briefly explain the notations and recall some 
 necessary facts that 
 will be frequently used in the sequel.
 
 \subsection{Crossed product by Kac algebras.}
Throughout this article $H (= H(\mu, \eta, \Delta, \varepsilon, S, *))$ will denote a finite-dimensional Kac algebra 
and $\delta$, the positive square root of $dim ~H$. 
We set $H^i = H$ or $H^*$ according as $i$ is odd or even. The unique non-zero idempotent integrals of $H^*$ and $H$ will be denoted by
$\phi$ and $h$ respectively and moreover, for any non-negative integer $i$, the symbols $\phi^i$ and $h^i$ will always denote a copy of 
$\phi$ and $h$ respectively. It is a fact that $\phi(h) = \frac{1}{dim ~H}$. 
The letters $x, y, z, t$ 
will always denote an element of $H$ and for any integer $i$, the symbols $x^i, y^i, z^i, t^i$ will always represent an element of $H$.
The letters $f, g, k$ will always denote an element of $H^*$ and for any integer $i$, the symbols $f^i, g^i$ and $k^i$
will always represent an element of $H^*$.

Given $x \in H, \Delta(x)$ is denoted by $x_1 \otimes x_2$ (a simplified version of the Sweedler coproduct notation). 
We draw the reader's attention to a notational abuse of which we will often be guilty. We denote elements
of a tensor product as decomposable tensors with the understanding that there is an implied omitted
summation (just as in our simplified Sweedler notation). Thus, when we write `suppose $f \otimes x \in H^* \otimes H$',
we mean `suppose $\sum_i f^i \otimes x^i \in H^* \otimes H$' (for some $f^i \in H^*$ and $x^i \in H$, the sum over a finite index set).

We refer to \cite[$\S 1$]{Sde2018} for the notion of action of $H$ on a finite-dimensional complex $*$-algebra say, $A$ and the construction of the
corresponding crossed product algebra, denoted $A \rtimes H$. Though the vector space underlying $A \rtimes H$ is $A \otimes H$, we denote a 
general element of $A \rtimes H$ by $a \rtimes x$ instead of $a \otimes x$.
There is a natural action of $H^*$ on $H$ given by 
 $f .  x = f(x_2) x_1$ for $f \in H^*, x \in H.$ Similarly we have action of $H$ on $H^*$. If $H$ acts on $A$, 
 then $H^*$ also acts on $A \rtimes H$ just by acting on $H$-part and ignoring the $A$-part, meaning that, 
 $f. (a \rtimes x) = a \rtimes f.x = f(x_2) a \rtimes x_1$ for $f$ in $H^*$ and $a \rtimes x \in A \rtimes H$ and consequently, we can construct $A \rtimes H \rtimes
 H^*$. Continuing this way, we may construct $A \rtimes H \rtimes H^* \rtimes \cdots$.
 
 For integers $i \leq j$, we define $H_{[i, j]}$ to be the crossed product algebra $H^i \rtimes H^{i+1} \rtimes \cdots \rtimes H^j$. 
 If $i = j$, we will simply write $H_i$ to denote $H_{[i, i]}$ and if $i > j$, we take $H_{[i, j]}$ to be $\mathbb{C}$.
 A typical element of $H_{[i, j]}$ will be denoted by $x^i / f^i \rtimes f^{i+1} / x^{i+1} \rtimes \cdots$ ($j-i+1$ terms).
 We use the symbol $A(H)_l$, where $l$ is any positive integer, to denote the crossed product algebra
 $\underbrace{H \rtimes H^* \rtimes \cdots}_{\text{$l$ \mbox{terms}}}$.
 
 Following \cite[$\S 1$]{Sde2018}, we denote by $H_{(-\infty, \infty)}$ the algebra which, by definition, is the `union' of all the $H_{[i, j]}$. Note
 that a typical element of $H_{(-\infty, \infty)}$ is a finite sum of terms of the form $\cdots \rtimes x^{-1} \rtimes f^0 \rtimes x^1 \rtimes
 \cdots$ where in any such term all but finitely many of the $f^i$ are $\epsilon$ and all but finitely many of the $x^i$ are $1$. For any 
 integer $m$, $H_{[m, \infty)}$ denotes the subalgebra of $H_{(-\infty, \infty)}$ which consists of all (finite sums of) elements 
 $\cdots \rtimes x^{-1} \rtimes f^0 \rtimes x^1 \rtimes
 \cdots$ of $H_{(-\infty, \infty)}$ where for $i < m$, $f^i = \epsilon$ if $i$ even and $x^i = 1$ if $i$ is odd. Similarly, we define the subalgebra 
 $H_{(-\infty, m]}$ of $H_{(-\infty, \infty)}$.  
 It is worth mentioning that the family $\{H_{(-\infty, -1]} \otimes H_{[m, \infty)} \subset H_{(-\infty, \infty)} : m > 1 \}$ of inclusions
 of infinite iterated crossed product algebras will be used in $\S 2$ to construct models for quantum double inclusions
 associated to the family of Kac algebra subfactors given by 
 $\{ R^H \subset R \rtimes \underbrace{H \rtimes H^* \rtimes \cdots}_{{\text{$m$ times}}} : m \geq 1 \}$.
 The following results will be very useful. We refer to \cite[Theorem 2.1, Corollary 2.3(ii)]{BlaMnt1985} for the proof of Lemma \ref{matrixalg},
\cite[Lemma 4.5.3]{Jjo2008} or \cite[Proposition 3]{DeKdy2015} for the proof of Lemma 
\ref{commutants} and \cite[Lemma 4.2.3]{Jjo2008} for the proof of  Lemma \ref{anti1}.
\begin{lemma}\label{matrixalg}
 $H \rtimes H^* \rtimes H \rtimes \cdots$ ($2k$-terms) is isomorphic to the matrix algebra $M_{n^k}(\mathbb{C})$ where $n = dim ~H$.
\end{lemma}

\begin{lemma}\label{commutants}
For any $p\in {\mathbb Z}$, the subalgebras $H_{(-\infty,p]}$ and $H_{[p+2,\infty)}$ are mutual commutants in $H_{(-\infty, \infty)}$.
\end{lemma} 
Given integers $i \leq j$ and $p \leq q$ such that $j-i = q-p$ and assume that $j$ and $p$ (resp., $i$ and $q$) have the same parity. 
Given $X \in H_{[i, j]}$, let $X^{\prime}$ denote the element obtained by `flipping $X$ about $i$ (equivalently, $j$)' and then applying 
$S^{\otimes(j-i+1)}$ on this flipped element. For instance, if we assume $i$ to be odd and $j$ to be even and if  
$X = x^i \rtimes f^{i+1} \rtimes \cdots \rtimes f^j \in H_{[i, j]}$, then $X^{\prime}$ is given by 
$Sf^j \rtimes Sx^{j-1} \rtimes \cdots \rtimes Sf^{i+1} \rtimes Sx^i$. It is evident that $X^{\prime} \in H_{[p, q]}$. 
\begin{lemma}\label{anti1}
 The map $X \mapsto X^{\prime}$ is a $*$-anti-isomorphism of $H_{[i, j]}$ onto $H_{[p, q]}$ . 
\end{lemma}

 The Fourier transform map $F_H : H \rightarrow H^*$ is defined by 
$F_H(a) = \delta \phi_1(a)\phi_2$ and
satisfies $F_{H^*} F_H = S$. We will usually omit the subscript of $F_H$ and $F_{H^*}$ and write
both as $F$ with the argument making it clear which is meant.

 \subsection{Planar algebras.}
 For the basics of (subfactor) planar algebras, we refer to \cite{Jns1999}, \cite{KdyLndSnd2003}
and \cite{KdySnd2004}. We will use the older notion of planar algebras where $Col$, the set of colours, is given by 
$\{(0, \pm),1,2,\cdots\}$ (note that only $0$ has two variants, namely, $(0, +)$ and $(0,-)$).
This is equivalent to the 
newer notion of planar algebras (see \cite[$\S 2.2$]{DeKdy2016}) where $Col = \{(k, \pm) : k \geq 0 \ \mbox{integer}\}$ and we refer to 
\cite[Proposition 1]{DeKdy2016} for the proof of this equivalence.
We will use the notation $T^{k_0}_{k_1, k_2, \cdots, k_b}$ to denote a tangle $T$ of colour $k_0$ (i.e., the colour of the external
box of $T$ is $k_0$) with $b$ internal
 boxes ($b$ may be zero also) such that the colour of the $i$-th internal box is $k_i$. Given a tangle $ T = T^{k_0}_{k_1, k_2, \cdots, k_b}$ and 
 a planar algebra $P$, $Z^P_T$ will always denote the associated linear map from $P_{k_1} \otimes P_{k_2} \otimes \cdots \otimes P_{k_b}$ to $P_{k_0}$ 
 induced by the tangle $T$. 
 
We will also find it useful to recall the notions of cabling and adjoints for tangles and for planar algebras.
Given any positive integer $m$, and a tangle $T$, say $T = T_{k_1, k_2, \cdots, k_b}^{k_0}$, 
the $m$-cabling of $T$, denoted by $T^{(m)}$, is the tangle obtained from $T$ by replacing each string of $T$ by a parallel cable of $m$-strings.
%
%
%
It is worth noting that the number of internal boxes of $T^{(m)}$ and $T$ are the same and that if $k_i(T^{(m)})$ denotes the colour of the 
$i$-th internal disc of $T^{(m)}$, then
\begin{align*}
        k_i(T^{(m)}) = 
        \begin{cases}
        mk_i, & \mbox{if} \ k_i > 0 \\
        (0, +), & \mbox{if} \ k_i = (0, +) \\
        (0, -), & \mbox{if} \ k_i = (0, -) \ \mbox{and} \  m \  \mbox{is odd}\\
        (0, +), & \mbox{if} \ k_i = (0, -) \ \mbox{and} \  m \ \mbox{is even}.\\
        \end{cases}
        \end{align*}
        Now given any planar algebra $P$, construct a new planar algebra $^{(m)}\!P$, called $m$-cabling of $P$, by setting
        \begin{align*}
        ^{(m)}\!P_k = 
        \begin{cases}
        P_{mk}, & \mbox{if} \ k > 0 \\
        P_{(0, +)}, & \mbox{if} \ k = (0, +) \\
        P_{(0, -)}, & \mbox{if} \ k = (0, -) \ \mbox{and} \  m \  \mbox{is odd}\\
        P_{(0, +)}, & \mbox{if} \ k = (0, -) \ \mbox{and} \  m \ \mbox{is even}\\
        \end{cases}
        \end{align*}
       and defining $Z_T^{^{(m)}\!P} = Z_{T^{(m)}}^P$ for any tangle $T$.  
Similarly,  given a planar algebra $P$, we construct a new planar algebra ${^*\!P}$, called the adjoint of $P$, where for 
      any $k \in Col, {(^*\!P)_k} = P_k$ as vector 
      spaces and given any tangle $T$, the action $Z_T^{^*\!P}$ of $T$ on $^*\!P$ is specified by $Z_{T^*}^P$ where $T^*$ is the tangle obtained by
      reflecting the tangle $T$ across any line in the plane.

\begin{figure}[!h]
\begin{center}
\psfrag{k}{\Huge $k$}
\psfrag{2}{\Huge $2$}
\psfrag{k-2}{\Huge $k-2$}
\resizebox{8.0cm}{!}{\includegraphics{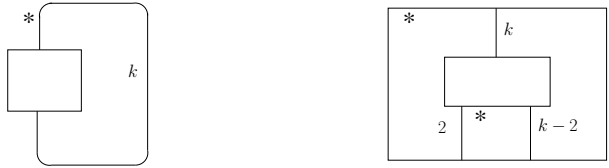}}
\end{center}
\caption{trace tangle : $tr_k^{(0, +)}$(left) and rotation tangle : $R_k^k (k \geq 2)$(right)}
\label{fig:pic67}
\end{figure}

\begin{figure}[!h]
\begin{center}
\psfrag{1}{\Huge $1$}
\psfrag{2}{\Huge $2$}
\psfrag{3}{\Huge $3$}
\resizebox{6.5cm}{!}{\includegraphics{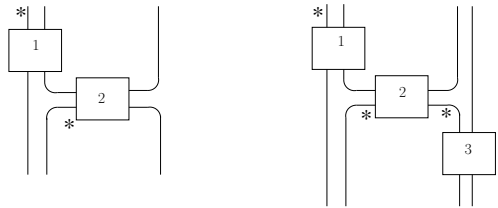}}
\end{center}
\caption{The tangles: $T^3$ (left) and $T^4$ (right)}
\label{fig:pic599}
\end{figure}

\begin{figure}[!h]
\begin{center}
\psfrag{a}{\Huge $1$}
\psfrag{b}{\Huge $2$}
\psfrag{c}{\Huge $2m$}
\psfrag{d}{\Huge $2m+2$}
\psfrag{e}{\Huge $2m+3$}
\psfrag{f}{\Huge $2m+1$}
\psfrag{g}{\Huge $2m+4$}
\psfrag{h}{\Huge $2m+5$}
\psfrag{i}{\Huge $2m+2n+3$}
\psfrag{1}{\Huge $1$}
\psfrag{2}{\Huge $2$}
\psfrag{3}{\Huge $3$}
\resizebox{11.0cm}{!}{\includegraphics{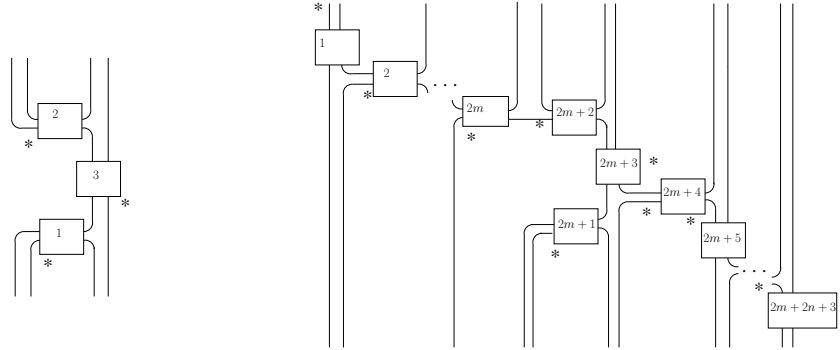}}
\end{center}
\caption{The tangles $A = A(0, 0)$(left) and $A(2m, 2n)$(right) when $m, n \geq 1$}
\label{fig:pic43}
\end{figure}  

Observe that Figures \ref{fig:pic599} and \ref{fig:pic43} show some elements of two families of tangles. In Figure
\ref{fig:pic599} we have the tangles $T^n$ of colour $n$ for $n \geq 2$, with exactly $n-1$ internal $2$-boxes and no internal regions 
illustrated for $n = 3$ and $n = 4$. In Figure \ref{fig:pic43} we have tangles $A(2m,2n)$ defined for $m,n \geq 0$ of colour $2m+2n+4$ with 
exactly $2m+2n+3$ internal $2$-boxes and no internal regions.

If $P$ is a subfactor planar algebra of modulus $d$, then for each $k \geq 1$, we refer to the (faithful, positive, normalised) trace 
 $\tau : P_k \rightarrow \mathbb{C}$ defined for $x \in P_k$ by 
 $\tau(x) = d^{-k} Z_{tr_k^{(0, +)}}(x)$ as the normalised pictorial trace on $P_k$ where $tr_k^{(0, +)}$ denotes the $(0, +)$ tangle with a single 
 internal $k$-box as shown in Figure \ref{fig:pic67}.
  \subsection{Planar algebra of a Kac algebra.}
\begin{figure}[!h]
\begin{center}
\psfrag{zab}{\huge $\zeta a + b$}
\psfrag{eq}{\huge $=$}
\psfrag{a}{\huge $a$}
\psfrag{b}{\huge $b$}
\psfrag{z}{\huge $\zeta$}
\psfrag{+}{\huge $+$}
\psfrag{del}{\huge $\delta$}
\resizebox{10.0cm}{!}{\includegraphics{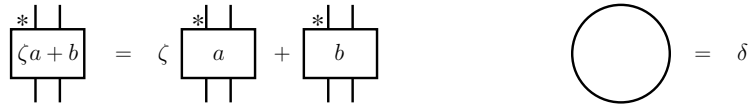}}
\end{center}
\caption{The L(inearity) and M(odulus) relations}
\label{fig:LnrMdl}
\end{figure}

\begin{figure}[!h]
\begin{center}
\psfrag{zab}{\huge $\zeta a + b$}
\psfrag{eq}{\huge $=$}
\psfrag{1h}{\huge $1_H$}
\psfrag{h}{\huge $h$}
\psfrag{z}{\huge $\zeta$}
\psfrag{+}{\huge $+$}
\psfrag{del}{\huge $\delta^{-1}$}
\resizebox{10.0cm}{!}{\includegraphics{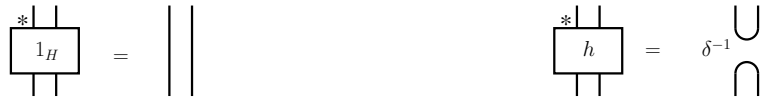}}
\end{center}
\caption{The U(nit) and I(ntegral) relations}
\label{fig:NitNtg}
\end{figure}

\begin{figure}[!h]
\begin{center}
\psfrag{epa}{\huge $\epsilon(a)$}
\psfrag{eq}{\huge $=$}
\psfrag{delinphia}{\huge $\delta \phi(a)$}
\psfrag{h}{\huge $h$}
\psfrag{a}{\huge $a$}
\psfrag{+}{\huge $+$}
\psfrag{del}{\huge $\delta^{-1}$}
\resizebox{10.0cm}{!}{\includegraphics{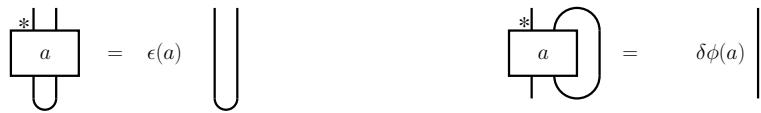}}
\end{center}
\caption{The C(ounit) and T(race) relations}
\label{fig:CntTrc}
\end{figure}

\begin{figure}[!h]
\begin{center}
\psfrag{epa}{\huge $\epsilon(a)$}
\psfrag{eq}{\huge $=$}
\psfrag{a1}{\huge $a_1$}
\psfrag{a2b}{\huge $a_2\,b$}
\psfrag{b}{\huge $b$}
\psfrag{a}{\huge $a$}
\psfrag{sa}{\huge $Sa$}
\psfrag{del}{\huge $\delta$}
\resizebox{10.0cm}{!}{\includegraphics{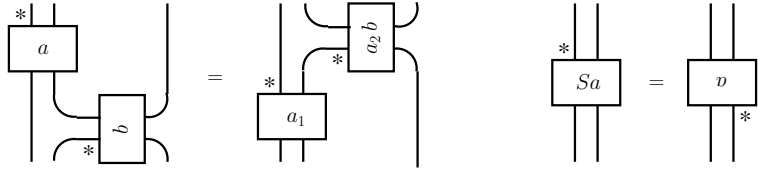}}
\end{center}
\caption{The E(xchange) and A(ntipode) relations}
\label{fig:XchNtp}
\end{figure}

  Suppose that $H$ acts outerly on the hyperfinite $II_1$ factor $M$. Let $P(H, \delta)$ (or, simply, $P(H)$) denote the subfactor planar algebra associated
to $M^H \subset M$ where $M^H$ is the fixed-point subalgebra of $M$. We recall from \cite[Theorem 8]{Sde2018} (see also \cite{KdySnd2006}) 
the construction of $P(H)$. The planar algebra $P(H)$ is defined to be the quotient of the 
universal planar algebra on the label set $L = L_2 = H$ by the set of relation in Figures \ref{fig:LnrMdl} - \ref{fig:XchNtp} 
(where (i) we write the relations
as identities - so the statement $a = b$ is interpreted as $a - b$ belongs to the set of relations; (ii) $\zeta \in k$  and  $a,b \in H;$ and
(iii) the external boxes of all tangles appearing in the relations are left 
undrawn and it is assumed that all external $*$-arcs are the leftmost arcs.

Note that the modulus relation is a pair
of relations - one for each choice of shading the circle. Finally, note that the interchange of $\delta$ and
$\delta^{-1}$ between the (I) and (T) relations here and those of \cite{KdySnd2006} is due to the different
normalisations of $h$ and $\phi$.

 A reformulation of Lemma 16 from \cite{KdySnd2006} will be useful.
Let ${\mathcal T}(k, p) (p \leq k-1)$ denote the set of $k$ tangles (interpreted as $0$ for $k=0$) 
with $p$ internal boxes of colour $2$ and no `internal regions'. If $p = k-1$, we will simply write 
${\mathcal T}(k)$ instead of ${\mathcal T}(k, k-1)$. The result then asserts:

\begin{lemma}\label{iso1}
For each tangle $X \in {\mathcal T}(k, p)$, the map $Z_X^{P(H)}:
(P(H)_2)^{\otimes p} \rightarrow P(H)_k$ is an injective linear map and if $p = k-1$, then $Z_X^{P(H)}:
(P(H)_2)^{\otimes k-1} \rightarrow P(H)_k$ is a linear isomorphism.
\end{lemma}

The following lemma (a reformulation of \cite[Proposition 4.3.1]{Jjo2008}) establishes algebra isomorphisms between $P(H)_k$ and finite iterated 
crossed product algebras.
\begin{lemma}\label{iso}
  For each $k \geq 2$, the map
 from $\underbrace{H \rtimes H^* \rtimes \cdots}_{\text{$k-1$ \mbox{terms}}}$ to $P(H)_k$ given by
 \begin{align*}
  \underbrace{x^1 \rtimes f^2 \rtimes \cdots}_{\text{$k-1$ \mbox{terms}}} \mapsto 
  Z^{P(H)}_{T^k}(\underbrace{x^1 \otimes Ff^2 \otimes \cdots}_{\text{$k-1$ \mbox{terms}}})
 \end{align*}
 is a $*$-algebra isomorphism. 
\end{lemma}
 We will use this identification of $\underbrace{H \rtimes H^* \rtimes \cdots}_{\text{$k-1$ \mbox{terms}}}$ with $P(H)_k$ very frequently
 without mention. 
Finally, for $i \leq j$, $tr_{H_{[i, j]}}$ denotes the faithful, positive, tracial state on $H_{[i, j]}$ given by
\begin{align*}
       \begin{cases}
        h^i \otimes \phi^{i+1} \otimes h^{i+2} \otimes \cdots (j-i+1 \mbox{-terms}), & \mbox{if} \  i \ \mbox{is even}\\
        \phi^i \otimes h^{i+1} \otimes \phi^{i+2} \otimes \cdots (j-i+1 \mbox{-terms}), & \mbox{if} \  i \  \mbox{is odd}.         
        \end{cases}
\end{align*}
Thus, for instance, if we assume $i$ to be odd, $j$ to be even 
and if $X \in H_{[i, j]}$, say,
$X = x^i \rtimes f^{i+1} \rtimes \cdots \rtimes x^{j-1} \rtimes f^j$,
then $tr_{H_{[i, j]}}(X) = \phi^i(x^i) f^{i+1}(h^{i+1}) \cdots \phi^{j-1}(x^{j-1}) f^j(h^j)$. 
\subsection{Drinfeld double construction.}
The Drinfeld double or quantum double construction is a construction that builds a quasitriangular Hopf algebra out of any finite-dimensional Hopf 
algebra. The Drinfeld double of $H$ is denoted by $D(H)$. The definition of $D(H)$ is not uniform in the literature.            
As in \cite{DeKdy2016} what we actually is an isomorphic variant of the version of $D(H)$ in \cite{Mjd2002} which has underlying vector space 
$H^* \otimes H$ and the structure maps are given by the following formulae: 
\begin{eqnarray*}
(f \otimes x)(g \otimes y) &=& g_1(x_1)g_3(Sx_3) (fg_2 \otimes yx_2),\\
\Delta(f \otimes x) &=& (f_2 \otimes x_2) \otimes (f_1 \otimes x_1), {\mbox { and}}\\
S(f \otimes x) &=& f_1(Sx_1)f_3(x_3) (S^{-1}f_2 \otimes Sx_2).
\end{eqnarray*}
One can easily verify that the structure maps of $D(H)^*$ are
given by the following formulae:
\begin{eqnarray*}
  & (f \otimes x) (g \otimes y) = gf \otimes yx,\\
   &\Delta(f \otimes x) = \delta^2 \phi_2(x_2) \phi_4(Sh_2) (\phi_1 f_2 S\phi_3 \otimes x_1) \otimes (f_1 \otimes h_1),\\  
  & S(f \otimes x) = \delta^2 \phi_4(x) \phi_2(h_2) \phi_1SfS\phi_3 \otimes h_1, \\
  & \varepsilon(f \otimes x) = f(1) \epsilon(x).\\
 \end{eqnarray*}
 Consider the linear isomorphism $Id_{H^*} \otimes F_H^{-1} : H^* \otimes H^* \rightarrow D(H)^*$. We can make $H^* \otimes H^*$ into a Kac algebra 
 where the structure maps are obtained by transporting the structure maps on $D(H)^*$ using this linear isomorphism. 
 Thus, by construction, $H^* \otimes H^*$ is isomorphic to 
 $D(H)^*$ as a Kac algebra. The following lemma explicitly describes the structure maps on $H^* \otimes H^*$.
 \begin{lemma}\label{dr}
  The structure maps on $H^* \otimes H^*$ are given by the following formulae:
 \begin{eqnarray*}
  & (g \otimes f) (k \otimes p) = \delta (Sf_2 p)(h) kg \otimes f_1,\\
   &\Delta(g \otimes f) = \delta (\phi_1 g_2 S\phi_3 \otimes f S\phi_2) \otimes (g_1 \otimes \phi_4),\\  
  & S(g \otimes f) = f_1 Sg Sf_3 \otimes Sf_2, \\
  & \varepsilon(g \otimes f) = \delta f(h) g(1).\\
 \end{eqnarray*}
  \end{lemma}
 \begin{proof}
  Easy to verify and is left to the reader.
 \end{proof}

\section{Construction of models for the quantum double inclusion}

In \cite[$\S 3$]{Sde2018} we defined the notion of quantum double inclusion associated to a finite-index and finite-depth subfactor and 
constructed a model for the quantum double inclusion of $R^H \subset R$. In a similar way we construct in this section 
models for the quantum double inclusions of the family of subfactors 
$\{ R^H \subset R \rtimes \underbrace{H \rtimes H^* \rtimes \cdots}_{{\text{$m$ times}}}  : m \geq 1 \}$. 

We begin with recalling from \cite{Sde2018} the notion of quantum double inclusion. 
Given a finite-index and finite-depth subfactor $N \subset M$, let $N ( = M_0) \subset M ( = M_1) \subset M_2 \subset M_3 \subset \cdots$
denote the basic construction tower of $N \subset M$. Let $M_{\infty}$ denote the $II_1$ factor obtained as the von Neumann closure 
$(\cup_{n = 0}^{\infty} M_n)^{\prime \prime}$ in the GNS representation with respect to the trace on $\cup_{n = 0}^{\infty} M_n$. Then the
inclusion 
\begin{align*}
 N \vee (M^{\prime} \cap M_{\infty}) \subset M_{\infty}
\end{align*}
is defined to be the quantum double inclusion associated to $N \subset M$.

It is well-known that for any positive integer $k$, $\mathbb{C} \subset H_k \subset H_{[k-1, k]} \subset H_{[k-2, k]} \subset H_{[k-3, k]} 
\subset \cdots$ is the basic construction tower associated to the initial (connected) inclusion $\mathbb{C} \subset H_k$ so that 
$H_{(-\infty, k]} ( = \cup_{i = 0}^{\infty} H_{[k-i, k]} )$ comes equipped with a tracial state and consequently,
\begin{align*}
 H_{(-\infty, k]}^{\prime \prime} := (\cup_{i=0}^{\infty} H_{[k-i, k]})^{\prime \prime} = (H_{(-\infty, k]})^{\prime \prime}
\end{align*}
turns out to be a hyperfinite $II_1$ factor. It is also well-known (see \cite[Theorem 4.11]{KdyLndSnd2003})
that the basic construction tower associated to $R^H \subset R$ is given by:
\begin{align*}
 R^H \subset R \subset R \rtimes H \subset R \rtimes H \rtimes H^* \subset R \rtimes H \rtimes H^* \rtimes H \subset \cdots.
\end{align*}
The following lemma (\cite[Lemma 17]{Sde2018}) describes models for $R^H \subset R$ as well as for the basic construction tower of 
$R^H \subset R$.
\begin{lemma}\label{z}
 $H_{(-\infty, -1]}^{\prime \prime} \subset H_{(-\infty, 0]}^{\prime \prime}$ is a model for $R^H \subset R$ 
     for some outer action of $H$ on the hyperfinite $II_1$ factor $R$ and $H_{(-\infty, -1]}^{\prime \prime} \subset 
     H_{(-\infty, 0]}^{\prime \prime} \subset H_{(-\infty, 1]}^{\prime \prime} \subset H_{(-\infty, 2]}^{\prime \prime} \subset \cdots$ is a 
     model for the basic construction tower of $R^H \subset R$.
\end{lemma}
As an immediate consequence of Lemma \ref{z} we obtain that $(\cup_{i=-1}^{\infty} H_{(-\infty, i]}^{\prime \prime})^{\prime \prime}$ is a hyperfinite $II_1$ factor. It is not hard to see 
that $(\cup_{i=-1}^{\infty} H_{(-\infty, i]}^{\prime \prime})^{\prime \prime} = (H_{(-\infty, \infty)})^{\prime \prime}$. We set 
\begin{align*}
 H^{\prime \prime}_{(-\infty, \infty)} := (H_{(-\infty, \infty)})^{\prime \prime}. 
\end{align*} 
It follows easily from Lemma \ref{z} and \cite[Proposition 4.3.6]{JnsSnd1997} that: 
\begin{lemma}
 Given any positive integer $m$, $H_{(-\infty, -1]}^{\prime \prime} \subset H_{(-\infty, m]}^{\prime \prime}$ is a model for $R^H \subset 
 R \rtimes \underbrace{H \rtimes H^* \rtimes \cdots}_{\text{$m$ times}}$ 
     and $H_{(-\infty, -1]}^{\prime \prime} \subset 
     H_{(-\infty, m]}^{\prime \prime} \subset H_{(-\infty, 2m+1]}^{\prime \prime} \subset H_{(-\infty, 3m+2]}^{\prime \prime} \subset \cdots$ 
     is a model for the basic construction tower of $R^H \subset R \rtimes \underbrace{H \rtimes H^* \rtimes \cdots}_{\text{$m$ times}}$. 
    \end{lemma}
Thus for any positive integer $m$, a model for the \textit{quantum double inclusion} of 
$R^H \subset R \rtimes \underbrace{H \rtimes H^* \rtimes \cdots}_{{\text{$m$ times}}}$ is given by 
     \begin{align*}
      H_{(-\infty, -1]}^{\prime \prime} \vee ((H_{(-\infty, m]}^{\prime \prime})^{\prime} \cap H^{\prime \prime}_{(-\infty, \infty)}) \subseteq
      H^{\prime \prime}_{(-\infty, \infty)}. 
     \end{align*}
     By an appeal to \cite[Lemma 14(2)]{Sde2018}, one can easily see that 
     \begin{align*}
     (H_{(-\infty, m]}^{\prime \prime})^{\prime} \cap H^{\prime \prime}_{(-\infty, \infty)} = H^{\prime \prime}_{[m+2, \infty)}
     \end{align*}
     and consequently, 
     \begin{align*}
      H_{(-\infty, -1]}^{\prime \prime} \vee 
     ((H_{(-\infty, m]}^{\prime \prime})^{\prime} \cap H^{\prime \prime}_{(-\infty, \infty)})  = 
     H_{(-\infty, -1]}^{\prime \prime} \vee H^{\prime \prime}_{[m+2, \infty)} = 
     (H_{(-\infty, -1]} \otimes H_{[m+2, \infty)})^{\prime \prime}.
     \end{align*}
      \begin{definition}\label{def1}
      For each integer $m > 2$, set $\mathcal{N}^m = (H_{(-\infty, -1]} \otimes H_{[m, \infty)})^{\prime \prime}$ and 
      $\mathcal{M} = H_{(-\infty, \infty)}^{\prime \prime}$.
     \end{definition}
     We have thus shown that:
     \begin{proposition}\label{qdim}
      For each integer $m > 2$, the subfactor $\mathcal{N}^m \subset \mathcal{M}$ is a model for the quantum double inclusion of $R^{H} \subset 
      R \rtimes \underbrace{H \rtimes H^* \rtimes \cdots}_{{\text{$m-2$ times}}}$. 
     \end{proposition}

     \section{Basic construction tower of $\mathcal{N}^m \subset \mathcal{M}, m > 2$ and relative commutants}
The purpose of this section is to construct the basic construction tower associated to $\mathcal{N}^m \subset \mathcal{M} \ (m > 2)$ and also to 
compute the relative commutants. 
\subsection{Some finite-dimensional basic constructions.}
This subsection is devoted to analysing the basic constructions associated to certain unital inclusions of finite-dimensional $C^*$-algebras.
We begin with recalling the following lemma (a reformulation of Lemma 5.3.1 of \cite{JnsSnd1997}) which provides an abstract characterisation of 
the basic construction associated to a unital inclusion of finite-dimensional $C^*$-algebras.
\begin{lemma}\cite[Lemma 5.3.1]{JnsSnd1997} \label{basic}
 Let $A \subseteq B \subseteq C$ be a unital inclusion of finite-dimensional $C^*$-algebras. Let $tr_B$ denote a faithful tracial state on $B$ and
 let $E_A$ denote the $tr_B$-preserving conditional expectation of $B$ onto $A$. Let $f \in C$ be a projection. Then $C$ is isomorphic to the 
 basic construction for  $A \subseteq B$ with $f$ as the Jones projection if the following conditions are satisfied:
 \begin{itemize}
 \item[(i)] $f$ commutes with every element of $A$ and $a \mapsto af$ is an injective map of $A$ into $C$,
 \item[(ii)] $f$ implements the trace-preserving conditional expectation of $B$ onto $A$ i.e., $fbf = E_A(b)f$ for all $b \in B$, and
 \item[(iii)] $BfB = C$.
 \end{itemize}
 \end{lemma}
 In the next lemma, we explicitly compute certain conditional expectation map. 
\begin{lemma}\label{exp}
 Given integers $l, p \geq 1$ and $s \geq 0$, let $\psi_{l, s, p}$ denote the embedding of $H_{[-l, p+s]}$ inside 
   $H_{[-l, -1]} \otimes H_{[p, 3p+s]}$ specified as follows:
   \begin{align*}
 \mbox{Let} \ X = x^{-l}/f^{-l} \rtimes \cdots \rtimes x^{p+s}/f^{p+s} \in H_{[-l, p+s]}, 
\end{align*}
then $\psi_{l, s, p}(X) \in H_{[-l, -1]} \otimes H_{[p, 3p+s]}$ is given by 
\begin{align*}
 (x^{-l}/f^{-l} \rtimes \cdots \rtimes f^{-2} \rtimes x^{-1}_1) \otimes 
 (\underbrace{1 \rtimes \epsilon \rtimes 1 \rtimes \cdots \rtimes \epsilon}_\text{$p-1$ \mbox{terms}}
 \rtimes \ x^{-1}_2 \rtimes f^0 \rtimes \cdots \rtimes x^{p+s}/f^{p+s})
\end{align*}
or
\begin{align*}
 (x^{-l}/f^{-l} \rtimes \cdots \rtimes f^{-2} \rtimes x^{-1}_1) \otimes (\underbrace{\epsilon \rtimes 1 \rtimes \epsilon \rtimes \cdots \rtimes \epsilon}_\text{$p-1$ \mbox{terms}}
 \rtimes \ x^{-1}_2 \rtimes f^0 \rtimes \cdots \rtimes x^{p+s}/f^{p+s})
\end{align*}
according as $p$ is odd or even. Then the trace-preserving conditional expectation $E$ of $H_{[-l, -1]} \otimes H_{[p, 3p+s]}$ onto 
$H_{[-l, p+s]}$ is given by
       \begin{align*}       
       & E((x^{-l}/ f^{-l} \rtimes \cdots \rtimes x^{-1}) \otimes (x^p/f^p \rtimes f^{p+1}/x^{p+1} \rtimes \cdots \rtimes x^{3p+s}/f^{3p+s}))\\
       & = \phi(Sx_2^{-1}x^{2p-1}) tr_{H_{[p, 2p-2]}}(x^p/f^p \rtimes \cdots \rtimes f^{2p-2}) x^{-l}/f^{-l} \rtimes \cdots \rtimes f^{-2} \rtimes
       x^{-1}_1 \rtimes f^{2p} \rtimes \cdots \rtimes x^{3p+s}/f^{3p+s}. 
       \end{align*}
       \end{lemma}
       \begin{proof}
       In \cite[Lemma 21(ii)]{Sde2018} we proved the result for $p = 2$. The proof for the general case will follow in a similar fashion
       and hence, we omit the proof.
       \end{proof}
       Next, we apply Lemma \ref{exp} to explicitly describe certain basic constructions and their associated Jones 
       projections.
       \begin{proposition}\label{basic cons}
 The following are instances of basic constructions with the Jones projections being specified pictorially in appropriate planar algebras.
 \begin{itemize}
 \item[1.]
If $l \geq 1, s \geq 0$  are integers, then given any positive integer $p$, 
$H_{[-l, -1]} \otimes H_{[p, p+s]} \subset H_{[-l, p+s]} \subset H_{[-l, -1]} \otimes H_{[p, 3p+s]}
(\subset H_{[-l, 3p+s]} \cong P(H^l)_{3p+s+l+2})$ is an instance of the basic construction with the Jones projection given by the following figure
\begin{center}
\begin{tikzpicture}
\node at (1.3,0) {\tiny $\delta^{-p}$};
 \draw [black,thick] (1.96,.5) -- (1.96,-.5);
 \node [right] at (2,0) {\tiny $p+l$}; 
 \draw [black,thick] (3,.5) to [out=275, in=265] (3.5,.5); 
  \draw [black,thick] (3,-.5) to [out=85, in=95] (3.5,-.5);
  \node [below] at (3.3,.4) {\tiny $p$};
  \node [above] at (3.3,-.4) {\tiny $p$};
  \draw [black,thick] (4,.5) -- (4, -.5);
  \node [right] at (4.002,0) {\tiny $s+2$};
 \end{tikzpicture}
 \end{center}
where the first inclusion is natural and the second inclusion is given by the map $\psi_{l, s, p}$ as defined in the statement of Lemma 
\ref{exp}. Furthermore, $tr_{H_{[-l, p+s]}}$ is a Markov trace of modulus $\delta^{2p}$ for the inclusion 
$H_{[-l, -1]} \otimes H_{[p, p+s]} \subset H_{[-l, p+s]}$.

 \item[2.]  If $l \geq 1, s \geq 0$ are integers, then given any positive integer $p$, 
      $H_{[-l, p+s]} \subset H_{[-l, -1]} \otimes H_{[p, 3p+s]} \subset H_{[-l, 3p+s]} (\cong P(H^l)_{3p+s+l+2})$ is an instance of the 
      basic construction with the  Jones projection given by
      \begin{center}
      \begin{tikzpicture}
 \node at (.7,0) {\tiny $\delta^{-p}$};
 \draw [black,thick] (1.4,.5) -- (1.4,-.5);
 \node [right] at (1.42,0) {\tiny $l$}; 
 \draw [black,thick] (2,.5) to [out=275, in=265] (2.5,.5); 
  \draw [black,thick] (2,-.5) to [out=85, in=95] (2.5,-.5);
  \node [below] at (2.3,.4) {\tiny $p$};
  \node [above] at (2.3,-.4) {\tiny $p$};
  \draw [black,thick] (3,.5) -- (3, -.5);
  \node [right] at (3.001,0) {\tiny $p+s+2$};
 \end{tikzpicture}
      \end{center}
      where the first inclusion is given by the map $\psi_{l, s, p}$ as described in the statement of Lemma \ref{exp} and the second inclusion is the natural 
      inclusion. Also, $tr_{H_{[-l, -1]} \otimes H_{[p, 3p+s]}}$ is a Markov trace of modulus $\delta^{2p}$ for the inclusion
      $H_{[-l, p+s]} \subset H_{[-l, -1]} \otimes H_{[p, 3p+s]}$.
      \end{itemize}
      \end{proposition}
      \begin{proof}
       In \cite[Proposition 22(2), 22(3)]{Sde2018} we proved the result for $p = 2$. 
      The proof for the general case will follow in a similar fashion. For the sake of completeness, we provide the proof of 
      only one part namely, part 2, of the proposition which is also the harder part.     
      \begin{itemize}
       \item [2.] We only present the proof when
       $l = 1$ and $s=0$, omitting the proof for the general case which is analogous. Let $e$ denote the projection defined in the statement
       of Proposition \ref{basic cons}(2).  We identify as usual $H_{[-1, 6m-3]}$ with $P(H)_{6m}$.
       
        Given $X = x^{-1} \rtimes f^0 \rtimes \cdots \rtimes x^{2m-1} \in H_{[-1, 2m-1]}$,
      its image in $H_{[-1, 6m-3]}$ is given by 
 \begin{align*}
  x^{-1} \rtimes \underbrace{\epsilon \rtimes 1 \cdots \rtimes \epsilon}_{\text{$4m-3$ \ factors}} \rtimes \ x^{-1}_2 \rtimes 
  f^0 \rtimes \cdots \rtimes x^{2m-1}.
 \end{align*}
  The element $eX$ is shown on the left in Figure \ref{fig:D44}. An application of the relation (E) shows that $eX$ equals the element on the 
  right in Figure \ref{fig:D44}. Similarly, by an appeal to the relations (A) and (E), one can easily see that the element $Xe$ equals the element on the 
        right in Figure \ref{fig:D44} so that $e X = X e$. Thus, we conclude that $e$ commutes with $X$. Further, it is evident from the pictorial representation
        of the element $Xe$ as shown on the right in Figure \ref{fig:D44} that the map $X \mapsto Xe$ of $H_{[-1, 2m-1]}$ into $H_{[-1, 6m-3]}$ 
        is injective,  verifying condition (i) of Lemma \ref{basic}.
         \begin{figure}[!h]
\begin{center}
\psfrag{a}{\Huge $\delta^{-(2m-1)}$}
\psfrag{b}{\Huge $x^{-1}_1$}
\psfrag{c}{\Huge $x^{-1}_2$}
\psfrag{d}{\Huge $Ff^0$}
\psfrag{e}{\Huge $=$}
\psfrag{f}{\Huge $x^{2m-1}$}
\psfrag{g}{\Huge $x^{2m-1}$}
\psfrag{h}{\Huge $x^{-1}$}
\psfrag{1}{\Huge $2m-1$}
\psfrag{2}{\Huge $2m-2$}
\psfrag{cc}{\Huge $\cdots$}
\resizebox{10.0cm}{!}{\includegraphics{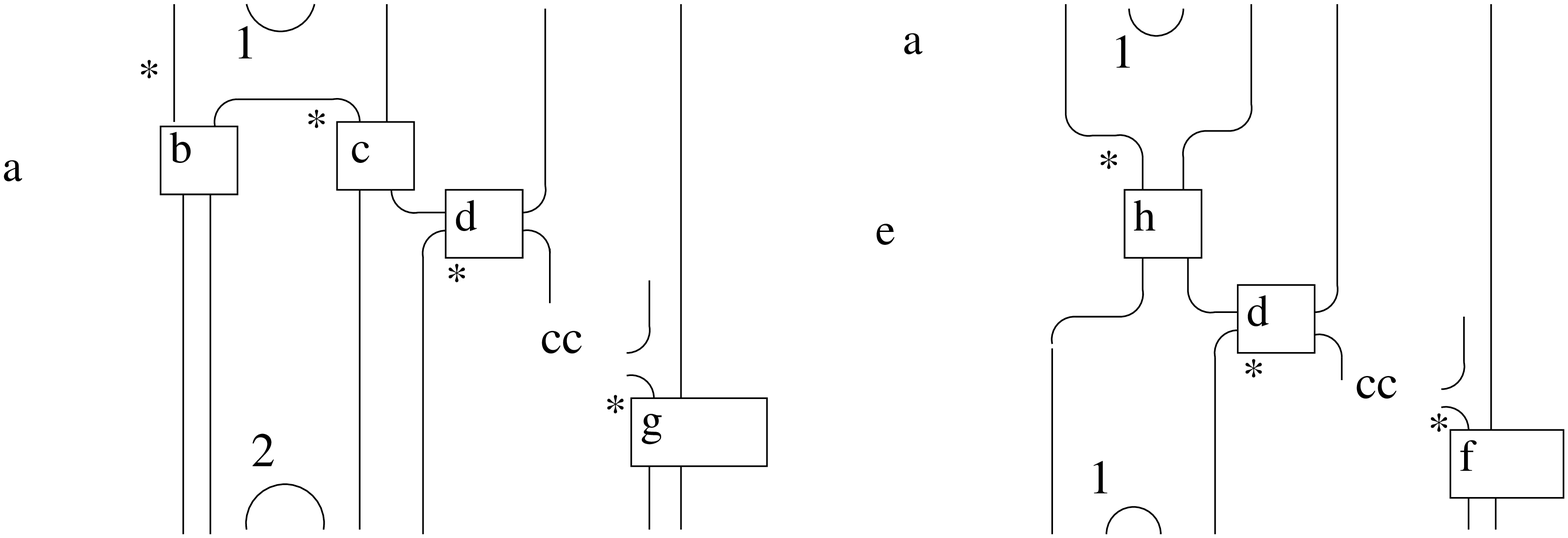}}
\end{center}
\caption{$eX = Xe$}
\label{fig:D44}
\end{figure}

Given $X = x^{-1} \otimes (x^{2m-1} \rtimes f^{2m} \rtimes \cdots \rtimes x^{6m-3})\in H_{-1} \otimes H_{[2m-1, 6m-3]}$,
      the element $eXe$ is shown in Figure \ref{fig:D444}. 
 \begin{figure}[!h]
\begin{center}
\psfrag{a}{\Huge $x^{-1}$}
\psfrag{b}{\Huge $x^{2m-1}$}
\psfrag{c}{\Huge $Ff^{2m}$}
\psfrag{d}{\Huge $x^{2m+1}$}
\psfrag{e}{\Huge $Ff^{4m-4}$}
\psfrag{f}{\Huge $x^{4m-3}$}
\psfrag{g}{\Huge $Ff^{4m-2}$}
\psfrag{h}{\Huge $x^{6m-3}$}
\psfrag{p}{\Huge $\delta^{-(4m-2)}$}
\psfrag{cc}{\Huge $\cdots$}
\psfrag{2}{\Huge $2m-1$}
\resizebox{11.0cm}{!}{\includegraphics{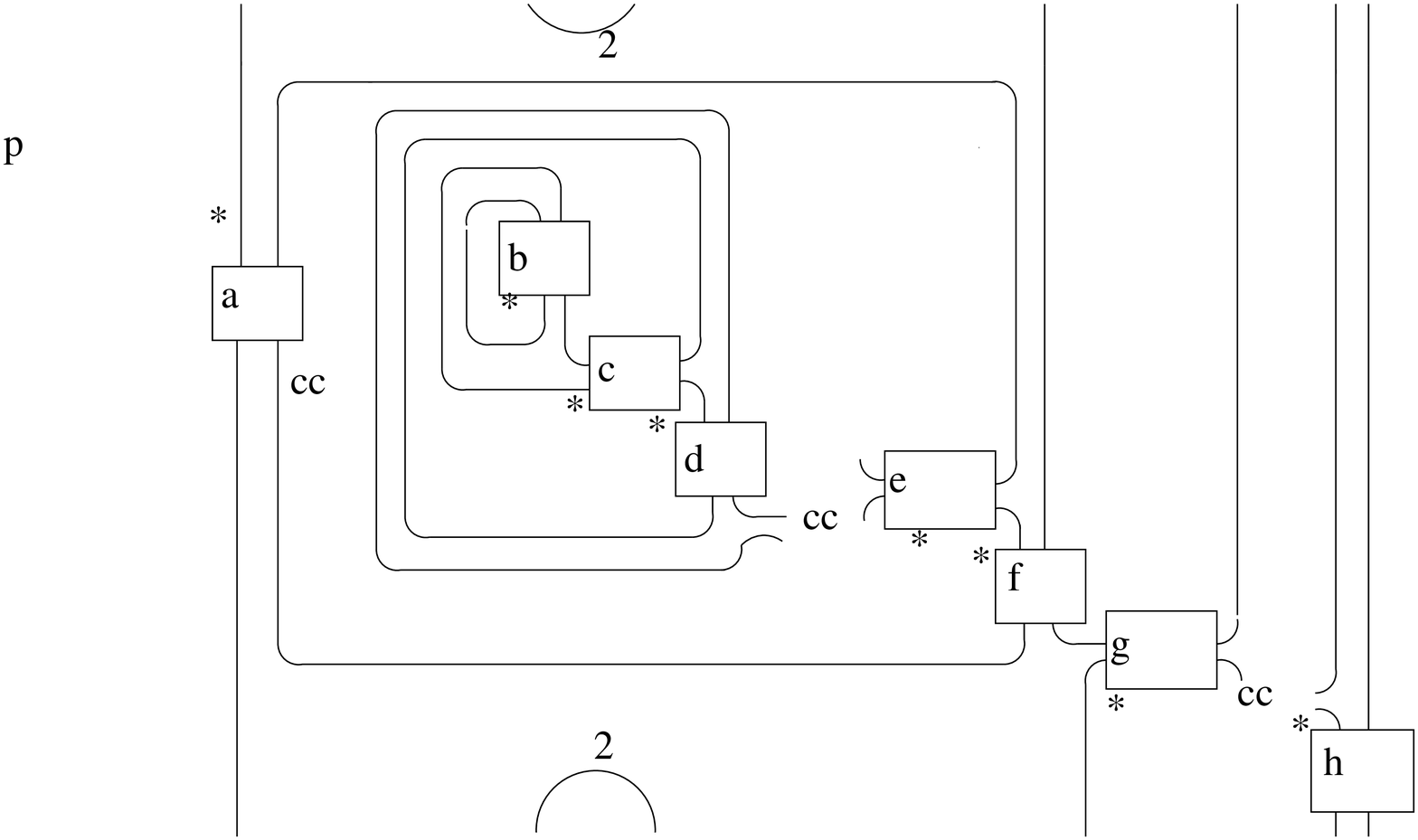}}
\end{center}
\caption{$eXe$}
\label{fig:D444}
\end{figure}        
       Repeated application of the relations (T), (C), and (A) reduces the element in Figure \ref{fig:D444} to that on the left
       in Figure \ref{fig:pic29}
         \begin{figure}[!h]
\begin{center}
\psfrag{a}{\Huge $x^{-1}$}
\psfrag{b}{\Huge $x^{4m-3}$}
\psfrag{c}{\Huge $Ff^{4m-2}$}
\psfrag{d}{\Huge $x^{4m-1}$}
\psfrag{e}{\Huge $x^{6m-3}$}
\psfrag{f}{\Huge $2m-1$}
\psfrag{g}{\Huge $2m-1$}
\psfrag{x}{\Huge $x^{-1}_1$}
\psfrag{i}{\Huge $\delta \ \alpha$}
\psfrag{j}{\Huge $\phi(Sx^{-1}_2 x^{4m-3})$} 
\psfrag{h}{\Huge $\alpha$}
\resizebox{12.0cm}{!}{\includegraphics{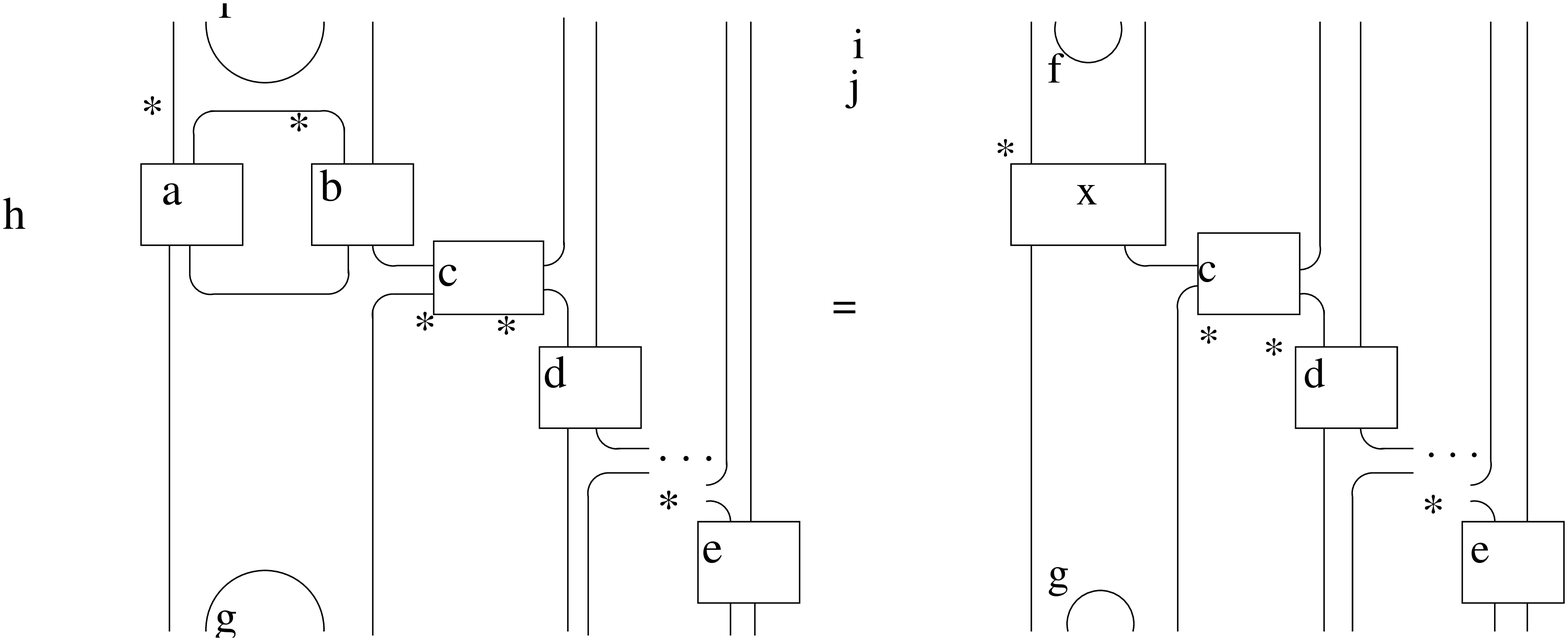}}
\end{center}
\caption{$eXe$}
\label{fig:pic29}
\end{figure}
where $\alpha = \delta^{-2m} tr_{H_{[2m-1, 4m-4]}}(x^{2m-1} \rtimes f^{2m} \rtimes \cdots \rtimes f^{4m-4})$. Again repeated application
of the relations (E) and (A), and finally, an application of the relation (T) reduces the element on the left in Figure \ref{fig:pic29} to that 
on the right in Figure \ref{fig:pic29}.
It follows from Lemma \ref{exp} that if $E$ denotes the trace-preserving conditional expectation of 
$H_{[-1, -1]} \otimes H_{[2m-1, 6m-3]}$ onto $H_{[-1, 2m-1]}$, then
\begin{align*}       
  E(X) = \phi(Sx_2^{-1}x^{4m-3}) tr_{H_{[2m-1, 4m-4]}}(x^{2m-1} \rtimes \cdots \rtimes f^{4m-4}) \ x^{-1}_1 \rtimes f^{4m-2} \rtimes \cdots 
  \rtimes x^{6m-3}. 
        \end{align*}
Now observe that $E(X)e$ equals the element as given by Figure \ref{fig:pic32} 
\begin{figure}[!h]
\begin{center}
\psfrag{a}{\Huge $x^{-1}_1$}
\psfrag{b}{\Huge $x^{-1}_2$}
\psfrag{c}{\Huge $Ff^{4m-2}$}
\psfrag{d}{\Huge $x^{4m-1}$}
\psfrag{e}{\Huge $x^{6m-3}$}
\psfrag{f}{\Huge $2m-2$}
\psfrag{g}{\Huge $2m-1$}
\psfrag{h}{\Huge $\delta \phi(Sx^{-1}_3 x^{4m-3}) \alpha$}
\resizebox{7.0cm}{!}{\includegraphics{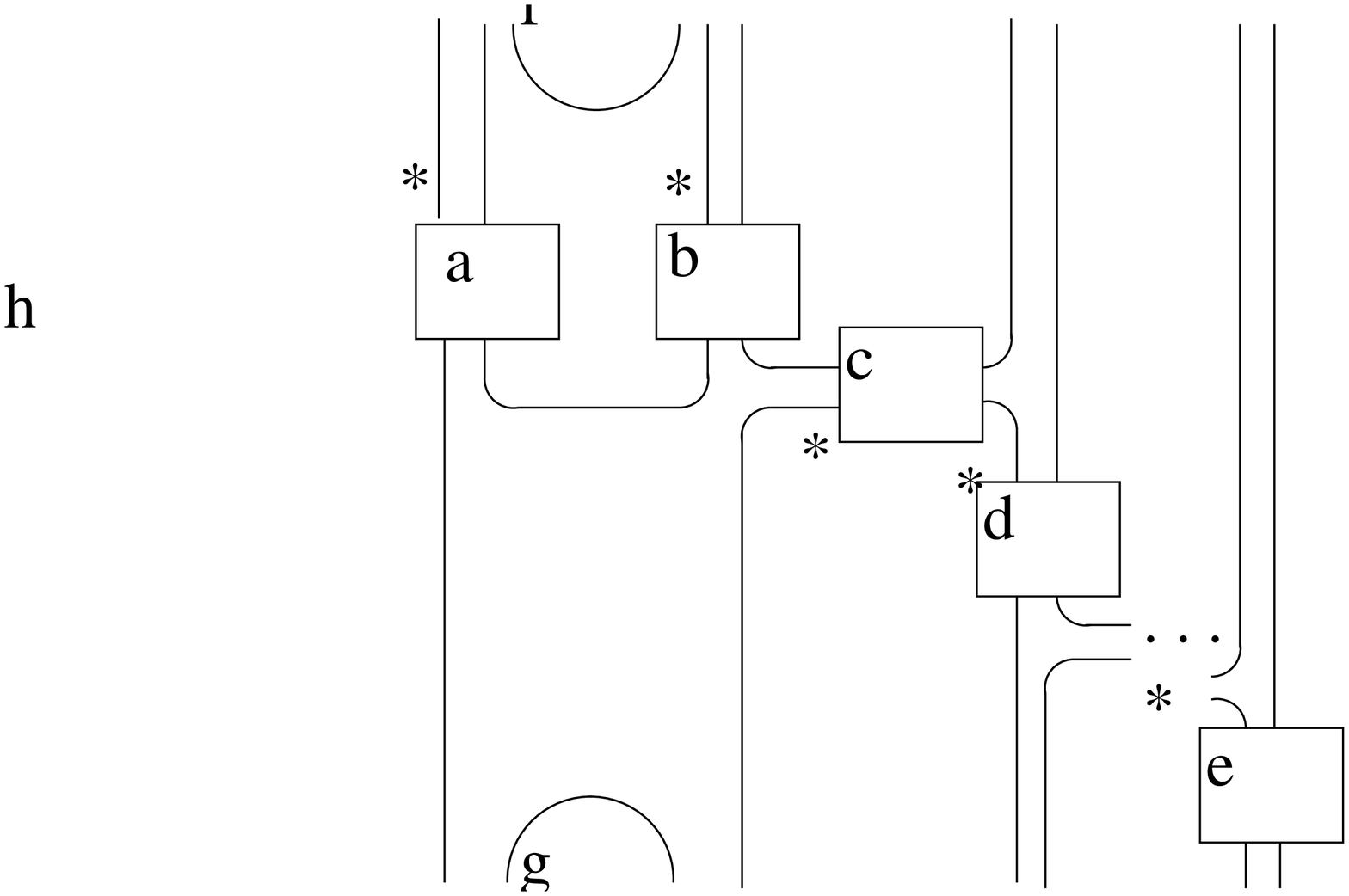}}
\end{center}
\caption{$E(X) e$}
\label{fig:pic32}
\end{figure}
which, after a straightforward computation
using relations (E) and (A), is easily seen to be equal to the element on the right in 
Figure \ref{fig:pic29}. Consequently, $eXe = E(X)e$, verifying condition (ii) of Lemma \ref{basic}.
         In order to verify condition (iii) of Lemma \ref{basic}, we just need to show that 
         $(H_{[-1, -1]} \otimes H_{[2m-1, 6m-3]}) e (H_{[-1, -1]} \otimes H_{[2m-1, 6m-3]}) = H_{[-1, 6m-3]}$.
Consider the elements $X, Y$ in $H_{[-1, -1]} \otimes H_{[2m-1, 6m-3]}$ given by 
        \begin{align*}
        & X = x^{-1} \otimes (x^{2m-1} \rtimes f^{2m} \rtimes \cdots \rtimes x^{6m-3}),\\ 
        &  Y = 1 \otimes (y^{2m-1} \rtimes g^{2m} \rtimes \cdots \rtimes y^{4m-3} \rtimes 
         \underbrace{\epsilon \rtimes 1 \rtimes \cdots \rtimes \epsilon \rtimes 1}_\text{$2m$ \mbox{terms}}).                  
        \end{align*}
        Representing the element $XeY$ pictorially in $P(H)_{6m}$ one can easily see that $XeY$ equals 
        \begin{align*}
  Z_T^{P(H)}(x^{-1}\otimes_{i=m}^{3m-2}(x^{2i-1} \otimes Ff^{2i}) \otimes x^{6m-3} \otimes_{i=m}^{2m-2}(y^{2i-1} \otimes Fg^{2i}) \otimes y^{4m-3})
        \end{align*}        
        where $Z_T$ is the linear 
        isomorphism induced by the tangle $T \in {\mathcal{T}}{(6m)}$ as shown in Figure \ref{fig:pic33}. Thus,
        we see that $(H_{[-1, -1]} \otimes H_{[2m-1, 6m-3]}) e (H_{[-1, -1]} \otimes H_{[2m-1, 6m-3]})$ contains the image of
        $Z_T$. Then by comparing dimensions of spaces we have that
        $ H_{[-1, 6m-3]} = (H_{[-1, -1]} \otimes H_{[2m-1, 6m-3]}) e (H_{[-1, -1]} \otimes H_{[2m-1, 6m-3]})$.         
 \begin{figure}[!h]
\begin{center}
\psfrag{a}{\Huge $1$}
\psfrag{b}{\Huge $2$}
\psfrag{c}{\Huge $3$}
\psfrag{d}{\Huge $2m-1$}
\psfrag{e}{\Huge $2m$}
\psfrag{f}{\Huge $2m+1$}
\psfrag{g}{\Huge $2m+2$}
\psfrag{h}{\Huge $4m$}
\psfrag{i}{\Huge $4m+1$}
\psfrag{j}{\Huge $4m+2$}
\psfrag{k}{\Huge $6m-2$}
\psfrag{l}{\Huge $6m-1$}
\resizebox{8.0cm}{!}{\includegraphics{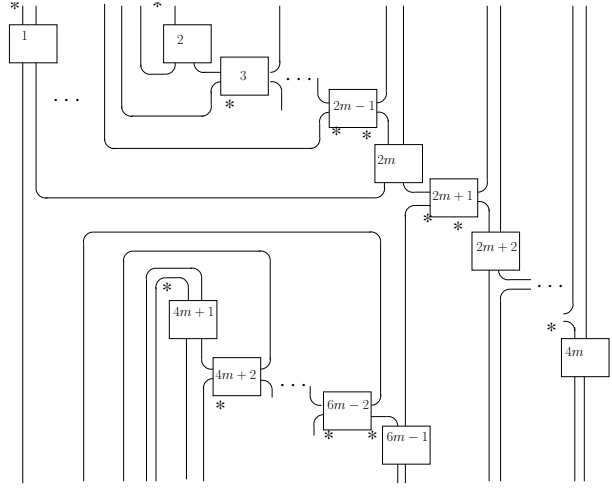}}
\end{center}
\caption{Tangle T}
\label{fig:pic33}
\end{figure}

Finally, a routine computation shows that for any 
$X \in H_{[-1, -1]} \otimes H_{[2m-1, 6m-3]}, \ tr(Xe) = \delta^{-2(2m-1)} tr(X)$,
 so that $tr_{H_{[-1, -1]} \otimes H_{[2m-1, 6m-3]}}$ is a Markov trace of modulus $\delta^{2(2m-1)}$ for the inclusion $H_{[-1, 2m-1]} \subset 
H_{[-1, -1]} \otimes H_{[2m-1, 6m-3]}$, completing the proof.          
      \end{itemize}

      \end{proof}
      
      \subsection{Jones' basic construction tower of $\mathcal{N}^m \subset \mathcal{M}$ and relative commutants}
      Throughout this subsection, $m > 2$ denotes a fixed positive integer.
      The goal of this subsection is to explicitly determine the basic construction tower of $\mathcal{N}^m \subset \mathcal{M}$.

\begin{figure}[!h]
 \begin{center}
 \psfrag{a}{\Large $A_{0, 1}$}
 \psfrag{b}{\Large $\subset$}
 \psfrag{c}{\Large $A_{1, 1}$}
 \psfrag{d}{\Large $\cup$}
 \psfrag{e}{\Large $\cup$} 
 \psfrag{f}{\Large $A_{0, 0}$}
 \psfrag{g}{\Large $\subset$}
 \psfrag{h}{\Large $A_{1, 0}$}
 \resizebox{3.0cm}{!}{\includegraphics{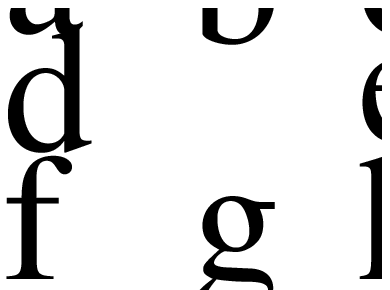}}
 \end{center}
 \caption{Commuting square}
 \label{fig:pic37}
 \end{figure}

\begin{figure}[!h]
\begin{center}
\psfrag{a}{\Huge $\delta^{-2}$}
\psfrag{b}{\Huge $1$}
\psfrag{c}{\Huge $1$}
\psfrag{d}{\Huge $m+2k+1$}
\psfrag{e}{\Huge $1$} 
\psfrag{f}{\Huge $1$}
\psfrag{x}{\Huge $\delta^{-3}$}
\psfrag{y}{\Huge $2$}
\psfrag{z}{\Huge $2$}
\psfrag{w}{\Huge $m+3k+1$}

\resizebox{7cm}{!}{\includegraphics{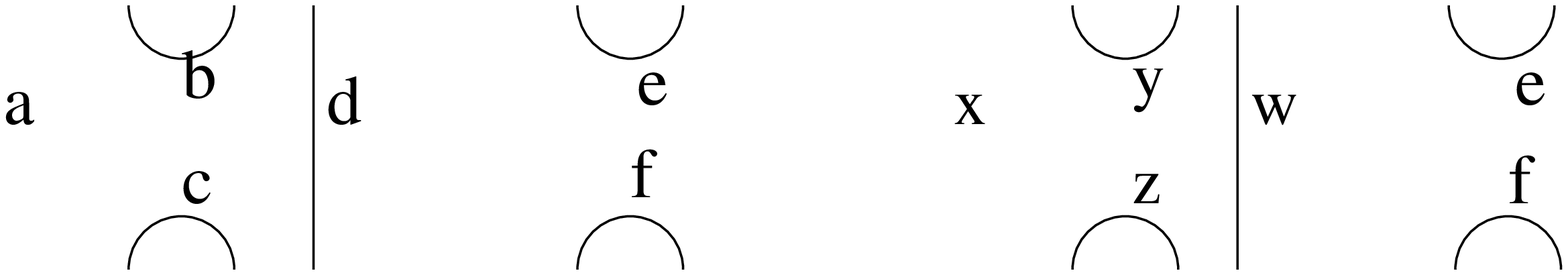}}
\end{center}
\caption{$e^{\prime}_{k+2} :$ with $m$ even (left) and with $m$ odd (right)}
\label{fig:pic100}
\end{figure}

We set $A_{0, 0} = \mathbb{C}, A_{0, 1} = H_{[0, m-1]}, A_{1, 0} = H_{-1} \otimes H_m$ or $H_{[-2, -1]} \otimes H_m$ 
according as $m$ is even or odd and $A_{1, 1} = H_{[-1, m]}$ or $H_{[-2, m]}$ according as $m$ is even  or odd. It follows from 
\cite[Lemma 23]{Sde2018} that
the square in Figure \ref{fig:pic37} is a symmetric commuting square
with respect to $tr_{A_{1, 1}}$ which is a Markov trace for the inclusion $A_{0, 1} \subset A_{1, 1}$. Further, here all the inclusions are 
connected since the lower left corner is $\mathbb{C}$ while the upper right corner is a matrix algebra by Lemma \ref{matrixalg}.
 For $k \geq 2$, we set 
 \begin{align*}
 A_{k, 1} = H_{[-k, m+k-1]} \ \mbox{or} \ H_{[-2k, m+k-1]} \ \mbox{according as} \ m \ \mbox{is even or odd}. 
 \end{align*}
 It is then a consequence of \cite[Proposition 22(1)(i)]{Sde2018} that 
   $A_{0, 1} \subset A_{1, 1} \subset A_{2, 1} \subset A_{3, 1} \subset \cdots$ is the basic construction tower associated to
   the initial inclusion $A_{0, 1} \subset A_{1, 1}$ and for any $k \geq 0$, if $e^{\prime}_{k+2}$ denotes the Jones projection lying in 
   $A_{k+2, 1}$ for the  basic construction of $A_{k, 1} \subset A_{k+1, 1}$, then  $e^{\prime}_{k+2}$ is given by Figure \ref{fig:pic100}.

Further, we define inductively 
   \begin{align*}
    A_{k+2, 0} = < A_{k+1, 0}, e^{\prime}_{k+2} >
   \end{align*}
   for each $k \geq0$.
   It is well-known that $A_{0, 0} \subset A_{1, 0} \overset{e^{\prime}_2} {\subset} A_{2, 0} \overset{e^{\prime}_3} {\subset} 
   A_{3, 0} \overset{e^{\prime}_4} \subset \cdots$ 
   is the basic construction tower of $A_{0, 0} \subset A_{1, 0}$. Proceeding along the same line as in the proof of \cite[Lemma 24]{Sde2018},
   one can show that:
   \begin{lemma}For any $k > 0$,  
   \begin{align*}   
     A_{k, 0} = 
      \begin{cases}
     H_{[-2k, -1]} \otimes H_{[m, m+k-1]}, & \mbox{if} \ m \  \mbox{is odd},\\
     H_{[-k, -1]} \otimes H_{[m, m+k-1]}, & \mbox{if} \ m \  \mbox{is even}.\\
    \end{cases}
   \end{align*}       
   \end{lemma}
   At this point we need to recall from \cite{KdySnd2008} the notion of finite pre-von Neumann algebras. 
By a finite pre-von Neumann algebra, we will mean a pair $(A, \tau)$ consisting of a complex $*$-algebra $A$ that is equipped with a normalised trace $\tau$ such that 
(i) the sesquilinear form defined by $<a, b> = \tau(b^* a)$ defines an inner-product on $A$ and such that (ii) for each $a \in A$, the
 left-multiplication map $\lambda_A(a) : A \longrightarrow A$ is bounded for the trace induced norm of $A$.
 By a compatible pair of finite pre-von Neumann algebras, we will mean
a pair $(A, \tau_A )$ and $(B, \tau_B)$ of finite pre-von Neumann algebras such that $A \subseteq B$ and ${\tau_B|}_{A} = \tau_A$.

If $A$ is a finite pre-von Neumann algebra with trace $\tau_A$, the symbol $L^2(A)$ will always denote
the Hilbert space completion of $A$ for the associated norm. Obviously, the left regular
representation $\lambda_A : A \rightarrow \mathcal{L}(L^2(A))$ is well-defined, i.e., for each $a \in A, \lambda_A(a) : A \rightarrow A$
extends to a bounded operator on $L^2(A)$. The notation $A^{\prime \prime}$ will always denote the von Neumann algebra 
$(\lambda_A(A))^{\prime \prime} \subset \mathcal{L}(L^2(A))$. The following lemma (a reformulation of \cite[Proposition 4.6(1)]{KdySnd2008})
 will be useful. 
\begin{lemma}\label{pre}\cite[Proposition 4.6(1)]{KdySnd2008} 
Let $(A, \tau_A )$ and $(B, \tau_B)$ be a compatible pair of finite pre-von Neumann algebras. The inclusion $A \subseteq B$
extends uniquely to a normal inclusion of $A^{\prime \prime}$ into $B^{\prime \prime}$ with image $(\lambda_B(A))^{\prime \prime}$. 
\end{lemma}

Note that $\cup_{k = 0}^{\infty} A_{k, 0} ( = H_{(-\infty, -1]} \otimes H_{[m, \infty)})$ and $\cup_{k = 0}^{\infty} A_{k, 1}
   ( = H_{(-\infty, \infty)})$ are finite pre-von Neumann algebras and $\cup_{k = 0}^{\infty} A_{k, 0} \subset \cup_{k = 0}^{\infty} A_{k, 1}$
   is a compatible pair so that by Lemma \ref{pre} the inclusion $\cup_{k = 0}^{\infty} A_{k, 0} \subset \cup_{k = 0}^{\infty} A_{k, 1}$ 
   extends uniquely to a normal inclusion $(\cup_{k = 0}^{\infty} A_{k, 0})^{\prime \prime} \subset 
   (\cup_{k = 0}^{\infty} A_{k, 1})^{\prime \prime}$. It follows from Definition \ref{def1} that 
   \begin{align*}
    (\cup_{k = 0}^{\infty} A_{k, 0})^{\prime \prime} = (H_{(-\infty, -1]} \otimes H_{[m, \infty)})^{\prime \prime} =  \mathcal{N}^m \ \mbox{and} \ 
    (\cup_{k = 0}^{\infty} A_{k, 1})^{\prime \prime} = (H_{(-\infty, \infty)})^{\prime \prime} = \mathcal{M}.
   \end{align*}
   Thus we have proved that:
   \begin{lemma}\label{hyper}
    $\mathcal{N}^m$ and $\mathcal{M}$ are hyperfinite $II_1$ factors.
   \end{lemma}
 The following lemma shows that $\mathcal{N}^m \subset \mathcal{M}$ is of finite index  
   equal to $\delta^{2m}$.
   \begin{lemma}\label{index}
    $[\mathcal{M} : \mathcal{N}^m] = \delta^{2m}$.
   \end{lemma}
   
    \begin{proof}
    It is well-known that (see \cite[Corollary 5.7.4]{JnsSnd1997}) $[\mathcal{M} : \mathcal{N}^m]$ equals the square of the norm of the 
    inclusion matrix for 
    $A_{0, 0} \subset A_{0, 1}$ which further equals the modulus of the Markov trace $tr_{A_{0, 1}}$ for the inclusion 
    $A_{0, 0} ( = \mathbb{C}) \subset A_{0, 1} ( = H_{[0, m -1]})$ which, again, by an application of \cite[Proposition 22(1)(ii)]{Sde2018},    
    equals $\delta^{2m}$. 
   \end{proof}
    For each $k \geq 0$ and $n \geq 2$, we now define a finite-dimensional $C^*$-algebra, denoted $A_{k, n}$, as follows.   
    \begin{itemize}
     \item \textit{Case (i): $m$ is odd}
      \begin{align*}
        A_{k, n} = 
  \begin{cases}
        H_{[-2k, -1]} \otimes H_{[m, (n+1)m +k - 1]}, & \mbox{if} \ n  \ \mbox{is even and} \ k > 0,  \\
        H_{[-2k, nm +k - 1]}, & \mbox{if} \ n \ \mbox{is odd and} \ k > 0,   \\
         H_{[-(n-1)m, m-1]}, & \mbox{if} \ k=0. \\
         \end{cases}
         \end{align*}
         \item \textit{Case (ii): $m$ is even}
          \begin{align*}
        A_{k, n} = 
  \begin{cases}
   H_{[-k, -1]} \otimes H_{[m, (n+1)m +k - 1]}, & \mbox{if} \ n \ \mbox{is even and} \ k > 0,  \\
        H_{[-k, nm +k - 1]}, & \mbox{if} \ n \ \mbox{is odd and} \ k > 0,  \\   
         H_{[-(n-1)m, m-1]}, & \mbox{if} \ k=0.   \\
  \end{cases}
  \end{align*}
   \end{itemize}
   We have already seen that for any $k \geq 0$, both the inclusions - the inclusion of $A_{k, 0}$ inside $A_{k, 1}$ and that 
    of $A_{k, n}$ inside $A_{k+1, n}$ for $n = 0, 1$, are natural. We describe below the embedding of $A_{k, n}$ inside $A_{k+1, n}$
   for any $k \geq 0, n \geq 2$ and that of $A_{k, n}$ inside $A_{k, n+1}$ for any $k \geq 0, n \geq 1$.
        \begin{itemize}
         \item For any $n \geq 2, k > 0, A_{k, n}$ sits inside $A_{k+1, n}$ in the natural way.
         \item If $n \geq 1$ is even and $k > 0$, then $A_{k, n}$ sits inside $A_{k, n+1}$ in the natural way.
         \item If $n > 0$ is odd and $k > 0$, the embedding of $A_{k, n}$ inside $A_{k, n+1}$ is given by $\psi_{l, s, p}$ as 
         defined in the statement of Lemma \ref{exp} with $p = m, s = (n-1)m+k-1$, and $l = 2k$ or $k$ according as $m$ is odd or even.
         \item If $n \geq 2$ is odd, then $A_{0, n}$ is identified with the subalgebra $H_{[0, mn-1]}$ of $A_{1, n}$. 
         \item If $n \geq 2$ is even, then $A_{0, n}$ is identified with the subalgebra $H_{[m, (n+1)m - 1]}$ of
         $A_{1, n}$.         
         \item Embedding of $A_{0, n}$ inside $A_{0, n+1}$ is natural for all $n \geq 1$.
        \end{itemize}
         Thus, we have a grid $\{ A_{k, n} : k, n \geq 0 \}$ of finite-dimensional $C^*$-algebras. 
        The following remark contains several useful facts concerning the grid $\{ A_{k, n} : k, n \geq 0 \}$.
        \begin{remark}\label{B}
         \begin{itemize}
          \item[(i)] We have already seen that the square of finite-dimensional $C^*$-algebras as shown 
          in Figure \ref{fig:pic37} is a symmetric commuting square with respect to $tr A_{1,1}$ which is a Markov trace 
          for the inclusion $A_{0,1} \subset A_{1,1}$ and all the inclusions are connected. Further, by Lemmas \ref{hyper} and \ref{index}, 
          $(\cup_{k=0}^{\infty} A_{k, 0})^{\prime \prime} ( = \mathcal{N}^m)$ as well as $(\cup_{k=0}^{\infty} A_{k, 1})^{\prime \prime} 
          ( = \mathcal{M})$
          are hyperfinite $II_1$ factors with $[\mathcal{M} : \mathcal{N}^m] = \delta^{2m}$.
          \item[(ii)]It follows from the embedding prescriptions that the following diagram (see Figure \ref{fig:cd}) commutes for all $k, n \geq 0$. 
         \begin{figure}[!h]
\begin{center}
\psfrag{a}{\Huge $A_{k, n+1}$}
\psfrag{b}{\Huge $\subset$}
\psfrag{c}{\Huge $A_{k+1, n+1}$}
\psfrag{x}{\Huge $A_{k, n}$}
\psfrag{y}{\Huge $\subset$} 
\psfrag{z}{\Huge $A_{k+1, n}$}
\psfrag{p}{\Huge $\cup$}
\psfrag{q}{\Huge $\cup$}
\resizebox{2.5cm}{!}{\includegraphics{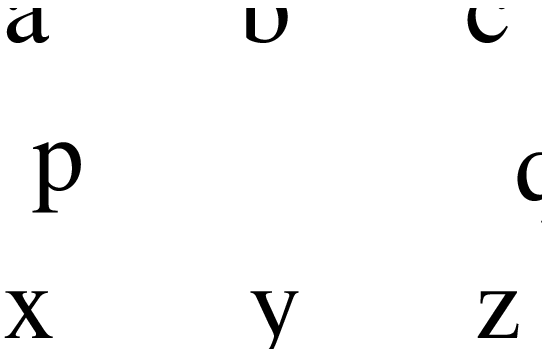}}
\end{center}
\caption{Commutative diagram}
\label{fig:cd}
\end{figure}
     \item[(iii)]It is a direct consequence of \cite[Proposition 22]{Sde2018} that
        for any $k, n \geq 0, A_{k, n} \subset A_{k, n+1} \subset A_{k, n+2}$ is an instance of the basic construction
         and further, $tr_{A_{k, n+1}}$ is a Markov trace of modulus $\delta^{2m}$ for the inclusion 
        $A_{k, n} \subset A_{k, n+1}$. Let $e_{k, n+2}^m$ ($k \geq 0, n \geq 0$) denote the Jones projection lying in $A_{k, n+2}$ applied to
          the basic construction $A_{k, n} \subset A_{k, n+1}$. 
\item[(iv)] For any $k \geq 0, n \geq 2$, the embedding of $A_{k, n}$ inside $A_{k+1, n}$ carries $e_{k, n}^m$ to $e_{k+1, n}^m$. 
         \end{itemize}
        \end{remark}
       
Obviously for any $n \geq 0$, $\cup_{k = 0}^{\infty}A_{k, n}$ is a finite pre-von Neumann algebra.       
 Consider the tower of finite pre-von Neumann algebras 
 \begin{align*}
  \cup_{k = 0}^{\infty}A_{k, 0} \subset \cup_{k = 0}^{\infty}A_{k, 1}
        \subset \cup_{k = 0}^{\infty}A_{k, 2} \subset \cdots.
 \end{align*}
Observe that for any $n \geq 0, \  
        \cup_{k = 0}^{\infty}A_{k, n} \subset \cup_{k = 0}^{\infty}A_{k, n+1}$ is a compatible pair so that by Lemma \ref{pre}, 
        the inclusion $\cup_{k = 0}^{\infty}A_{k, n} \subset \cup_{k = 0}^{\infty}A_{k, n+1}$ extends uniquely to a normal extension 
        $({\cup_{k = 0}^{\infty}A_{k, n}})^{\prime \prime} \subset ({\cup_{k = 0}^{\infty}A_{k, n+1}})^{\prime \prime}$. Note also that  
        $\cup_{k = 0}^{\infty} A_{k, n} = H_{(-\infty, -1]} \otimes H_{[m, \infty)} \ \mbox{or} \ 
    H_{(-\infty, \infty)}$ according as $n$ is even or odd.
    For each $n \geq 0$,  we define $\mathcal{M}_n:= (\cup_{k = 0}^{\infty}A_{k, n})^{\prime \prime}$. Then $\mathcal{M}_n = (H_{(-\infty, -1]} 
    \otimes H_{[m, \infty)})^{\prime \prime} \ \mbox{or} \ H^{\prime \prime}_{(-\infty, \infty)}$ according as $n$ is even or odd.
    In view of the facts concerning the grid $\{ A_{k, n} : k, n \geq 0 \}$ as mentioned in Remark \ref{B}, one can conclude that:
    \begin{proposition}\label{towerbasic}
    $\mathcal{M}_0 (=\mathcal{N}^m) \subset \mathcal{M}_1 (=\mathcal{M}) \subset \mathcal{M}_2 \subset \mathcal{M}_3 \subset \cdots$ is the 
    basic construction tower of $\mathcal{N}^m \subset \mathcal{M}$. 
    \end{proposition}
    \subsection{Computation of the relative commutants.}
\begin{figure}[!h]
\begin{center}
\psfrag{1}{\huge $m$}
\psfrag{2}{\Huge $nm+1$}
\psfrag{3}{\huge $nm$}
\psfrag{a}{\huge $\delta^{-m}$}
\resizebox{8cm}{!}{\includegraphics{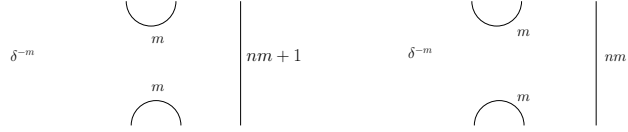}}
\end{center}
\caption{$e_{0, n+2}^m$ (left) and $\tilde{e}_{n+2}^m$ (right), $n \geq 0$}
\label{fig:newpic73}
\end{figure}

     We now proceed to compute the relative commutants. By virtue of Ocneanu's compactness theorem (see \cite[Theorem 5.7.6]{JnsSnd1997}), the relative commutant 
$(\mathcal{N}^m)^{\prime} \cap \mathcal{M}_k$ ($k > 0$) is given 
by
\begin{align*}
 (\mathcal{N}^m)^{\prime} \cap \mathcal{M}_k = A_{0,k} \cap (A_{1,0})^{\prime}, \ k \geq 1.
\end{align*}
 The following proposition describes the spaces $ A_{0, k} \cap (A_{1, 0})^{\prime}, k \geq 1$. 
 The proof of the proposition is similar to that of \cite[Proposition 29]{Sde2018} and we omit its proof.
      \begin{proposition}\label{rel comm}      
      Let $k \geq 1$ be an integer and set
      \begin{align*}
       &\tilde{Q}^m_{2k} = \{X \in H_{[m, (2k+1)m-2]} : 
       X \  \mbox{commutes with} \ \Delta_{k-1}(x) \in \otimes_{i = 1}^{k} H_{2im-1}, \forall x \in H \},\\      
      &\tilde{Q}^m_{2k-1} = \{X \in H_{[0, (2k-1)m-2]} : X \ \mbox{commutes with} \ \Delta_{k-1}(x) \in
       \otimes_{i = 0}^{k-1} H_{2im-1}, \forall x \in H\}.
        \end{align*}
        Then, 
       $A_{0, 2k} \cap (A_{1, 0})^{\prime} \cong \tilde{Q}^m_{2k}$ and $A_{0, 2k-1} \cap (A_{1, 0})^{\prime} \cong
      \tilde{Q}^m_{2k-1}$.
      \end{proposition}
        It follows from Remark \ref{B}(iii) that the Jones projection lying in $\mathcal{N}^{\prime} \cap \mathcal{M}_{n+2} = 
   A_{0, n+2} \cap (A_{1, 0})^{\prime}$ ($n \geq 0$)
   is given by $e_{0, n+2}^m$ (see Figure \ref{fig:newpic73}), which, under the identification of 
   $A_{0, n+2} \cap (A_{1, 0})^{\prime}$ with $\tilde{Q}_{n+2}^m$ as given by Proposition \ref{rel comm}, is easily seen to be identified with the 
   projection $\tilde{e}_{n+2}^m$ in $\tilde{Q}_{n+2}^m$ as shown on the right in Figure \ref{fig:newpic73}.
\begin{remark}
  It is worth knowing the embedding of $\tilde{Q}^m_k$ inside $\tilde{Q}^m_{k+1}$ ($k \geq 1$).
  It follows easily from the embedding formulae of $A_{1, k}$ inside $A_{1, k+1}$ and $A_{0, k}$ (resp., $A_{0, k+1}$) inside 
  $A_{1, k}$ (resp., $A_{1, k+1}$) and Proposition \ref{rel comm} that given 
  $X \in \tilde{Q}^m_k$,
  it sits inside $\tilde{Q}^m_{k+1}$ as 
  \begin{align*}
  \underbrace{\epsilon \rtimes 1 \rtimes \cdots \rtimes 1}_{\text{$m$ factors}} \rtimes X, \ \mbox{if} \ m \ \mbox{is even}
  \end{align*}
  and if $m$ is odd, then the image of $X \in \tilde{Q}^m_k$ inside $\tilde{Q}^m_{k+1}$ is given by
  \begin{align*}
    \underbrace{\epsilon \rtimes 1 \rtimes \cdots \rtimes \epsilon}_{\text{$m$ factors}} \rtimes X  \ \  \mbox{or} \  \  
     \underbrace{1 \rtimes \epsilon \rtimes \cdots \rtimes 1}_{\text{$m$ factors}} \rtimes X
  \end{align*}
  according as $k$ is even or odd. Also, the diagram in Figure \ref{fig:nnpic73} commutes
 where each horizontal arrow indicates the $*$-isomorphism.
  \begin{figure}[!h]
\begin{center}
 \psfrag{x}{\huge $(\mathcal{N}^m)^{\prime} \cap \mathcal{M}_{n+1}$}
\psfrag{y}{\huge $(\mathcal{N}^m)^{\prime} \cap \mathcal{M}_{n}$}
\psfrag{z}{\Huge $\tilde{Q}^m_{n+1}$}
\psfrag{w}{\Huge $\tilde{Q}^m_n$}
\psfrag{b}{\Huge $\longrightarrow$}
\psfrag{c}{\Huge $\cup$}
\resizebox{3.5cm}{!}{\includegraphics{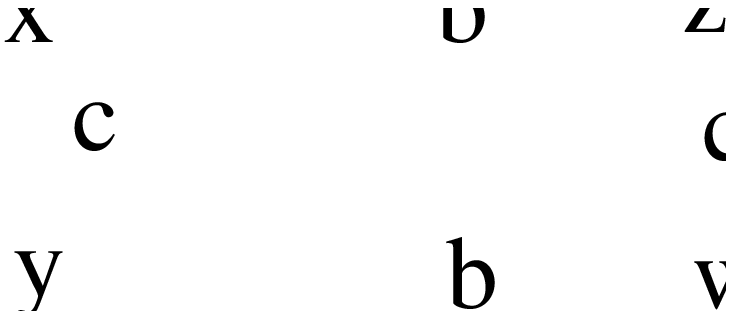}}
\end{center}
\caption{A commutative diagram} 
\label{fig:nnpic73}
\end{figure}
 \end{remark}
 For each integer $n \geq 1$, we define a subspace $Q^m_n$ of $H_{[1, mn-1]}$ or $H_{[0, mn-2]}$ according as $m$ is odd or even as follows:
  \begin{itemize}
     \item \textit{Case (i): $m$ is odd}
     \begin{align*}
       Q^m_{n} := \{ & X \in H_{[1, mn-1]} : X \leftrightarrow \Delta_{k-1}(x) \in \otimes_{i = 1}^k H_{m(2i-1)}, 
 \forall x \in H \text{ where } \\
  &k = \frac{n}{2} \ \mbox{if} \ n \ \mbox{is even or} \ \frac{n+1}{2} \ \mbox{if} \ n \ \mbox{is odd} \}
     \end{align*}
     \item \textit{Case (ii): $m$ is even}
    \begin{align*}                
  Q^m_{n} := \{ & X \in H_{[0,mn-2]} : X \leftrightarrow \Delta_{k-1}(x) \in \otimes_{i = 1}^k H_{m(2i-1)-1}, \forall x \in H \text{ where } \\ 
    &k = \frac{n}{2} \ \mbox{if} \ n \ \mbox{is even or} \ \frac{n+1}{2} \ \mbox{if} \ n \ \mbox{is odd} \}   
    \end{align*}    
     \end{itemize}

  This is an immediate consequence of Lemma \ref{anti1} that for any $n \geq 1$, $\tilde{Q}^m_n$ is $*$-anti-isomorphic to 
  $Q^m_n$ and let $\gamma_n^m: \tilde{Q}_n^m \rightarrow Q_n^m$ denote this anti-isomorphism. We then have the following commutative diagram.
 \begin{figure}[!h]
\begin{center}
\psfrag{a}{\Huge $\tilde{Q}_{n+1}^m$}
\psfrag{b}{\Huge $Q_{n+1}^m$}
\psfrag{c}{\Huge $\tilde{Q}_n^m$}
\psfrag{d}{\Huge $Q_n^m$}
\psfrag{e}{\Huge $\overset{\gamma_{n+1}^m}{\longrightarrow}$}
\psfrag{f}{\Huge $\cup$}
\psfrag{g}{\Huge $\overset{\gamma_n^m}{\longrightarrow}$}
\resizebox{3.0cm}{!}{\includegraphics{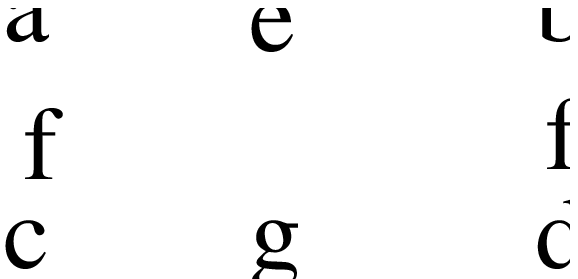}}
\end{center}
\caption{Commutative diagram}
\label{fig:D8}
\end{figure}
  
\begin{figure}[!h]
\begin{center}
\psfrag{1}{\Huge $m$}
\psfrag{3}{\Huge $m(n-2)$}
\psfrag{a}{\Huge $\delta^{-m}$}
\resizebox{4.0cm}{!}{\includegraphics{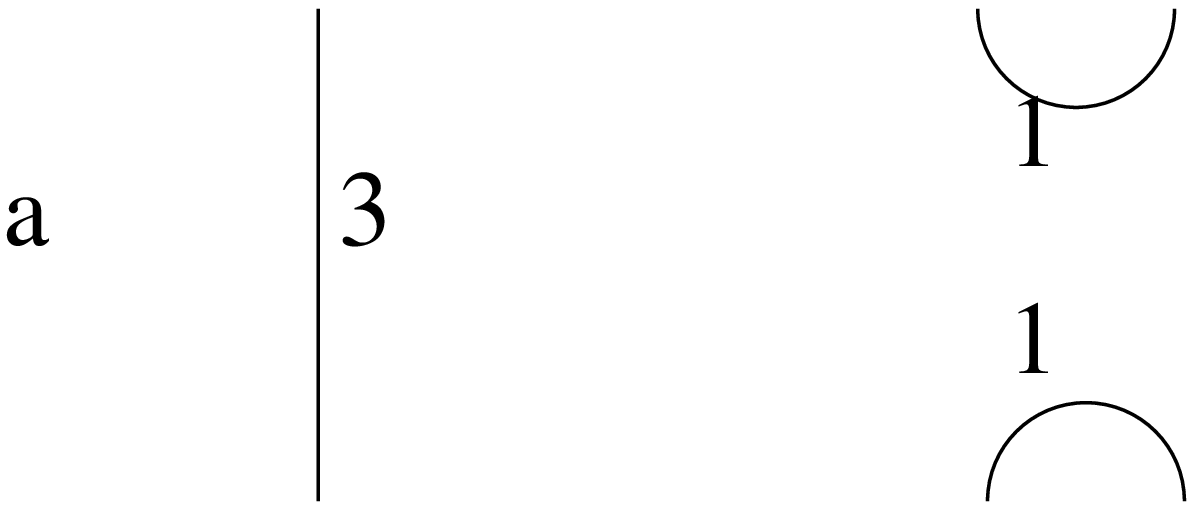}}
\end{center}
\caption{$e_n^m$}
\label{fig:pic77}
\end{figure}

\begin{figure}[!h]
\begin{center}
\psfrag{a}{\Huge $(\mathcal{N}^m)^{\prime} \cap \mathcal{M}_{n+1}$}
\psfrag{b}{\Huge $Q_{n+1}^{m^{op}}$}
\psfrag{c}{\Huge $(\mathcal{N}^m)^{\prime} \cap \mathcal{M}_n$}
\psfrag{d}{\Huge $Q_n^{m^{op}}$}
\psfrag{e}{\Huge $\overset{\Psi_{n+1}^m}{\longrightarrow}$}
\psfrag{f}{\Huge $\cup$}
\psfrag{g}{\Huge $\overset{\Psi_n^m}{\longrightarrow}$}
\resizebox{3.5cm}{!}{\includegraphics{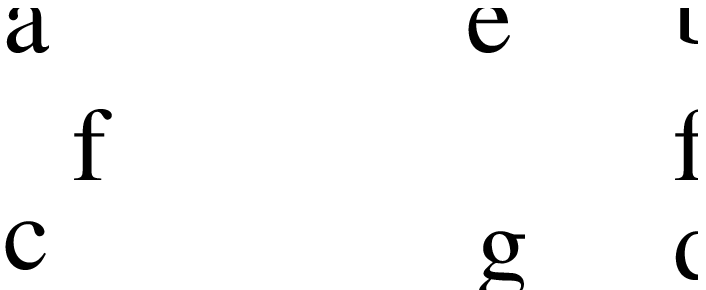}}
\end{center}
\caption{Commutative diagram}
\label{fig:D9}
\end{figure}

Further, if $e_n^m \in Q_n^m (n \geq 2)$ denotes the projection which is the image of $\tilde{e}_n^m \in \tilde{Q}_n^m$ under 
$\gamma_n^m$, 
  it is then not hard to see that $e_n^m$ is given by Figure \ref{fig:pic77}.
  Obviously, the identity map of $Q_n^m$ onto its opposite algebra, denoted $Q_n^{m^{op}}$, is a anti-$*$-isomorphism. 
  For each $n \geq 1$, let $\Psi_n^m: (\mathcal{N}^m)^{\prime} \cap \mathcal{M}_n \rightarrow$ $Q_n^{m^{op}}$ denote the following composite map: 
    \begin{align*}
     (\mathcal{N}^m)^{\prime} \cap \mathcal{M}_n \xrightarrow{*-\mbox{isom}} \tilde{Q}_n^m \xrightarrow{\gamma_n (*- \mbox{anti-isom})} 
     Q_n^m \xrightarrow{\mbox{Identity}} Q_n^{m^{op}}. 
    \end{align*}
    Obviously $\Psi_n^m$ is $*$-isomorphism for each $n \geq 1$ and for $n \geq 2$, it carries $e_{0, n}^m$ to $e_n^m$. The commutative 
    diagrams in Figures \ref{fig:nnpic73} and \ref{fig:D8} together imply commutativity of the diagram
    in Figure \ref{fig:D9}.

It will be useful to identify the spaces $\mathcal{M}^{\prime} \cap \mathcal{M}_n (n \geq 2)$ also. Once again applying Ocneanu's compactness 
theorem we obtain that $\mathcal{M}^{\prime} \cap \mathcal{M}_n (n \geq 2) = A_{0, n} \cap (A_{1, 1})^{\prime} (n \geq 2)$.
Proceeding along the same line of argument as in the proof of Proposition \ref{rel comm}, one can
show that: 
\begin{lemma}
 If $n$ is even, then $A_{0, n} \cap (A_{1, 1})^{\prime}$ can be identified with 
 \begin{align*}
  \{X \in H_{[m, mn-2]}: X \ \mbox{commutes with} \ \Delta_{k-1}(x) \in \otimes_{i=1}^k H_{2mi-1}, k = \frac{n}{2} \},
 \end{align*}
and if $n$ is odd, then $A_{0, n} \cap (A_{1, 1})^{\prime}$ can be identified with 
 \begin{align*}
 \{X \in H_{[0, m(n-1)-2]}: X \ \mbox{commutes with} \ \Delta_{k}(x) \in \otimes_{i=0}^k H_{2mi-1}, k = \frac{n-1}{2} \}. 
 \end{align*}
 \end{lemma}
 As an immediate consequence of this lemma, we obtain that:
 \begin{lemma}\label{second}
The $*$-isomorphism $\Psi_n^m$ of $(\mathcal{N}^m)^{\prime} \cap \mathcal{M}_n$ onto $Q_n^{m^{op}}$ carries 
$\mathcal{M}^{\prime} \cap \mathcal{M}_n (n \geq 2)$ 
onto the subspace of $Q_n^{m^{op}}$ given by
\begin{itemize}
 \item[(i)] $m$ is odd:
 \begin{align*}
  \{X \in H^{op}_{[m+1, mn-1]} : X \ \mbox{commutes with} \ \Delta_k(x) \in \otimes_{i=1}^k H_{(2i-1)m}, \forall x \in H\},
   \end{align*}
   \item[(ii)] $m$ is even:
   \begin{align*}
  \{X \in H^{op}_{[m, mn-2]} : X \ \mbox{commutes with} \ \Delta_k(x) \in \otimes_{i=1}^k H_{(2i-1)m-1}, \forall x \in H\},  
   \end{align*}
   where, in either case, $k = \frac{n}{2}$ or $\frac{n+1}{2}$ according as $n$ is even or odd.
\end{itemize}
\end{lemma}
In the next lemma we consider the question of irreducibility of $\mathcal{N}^m \subseteq M$ for $m \geq 2$.
\begin{lemma}\label{irr}
  $\mathcal{N}^m \subset \mathcal{M}$ is reducible for all $m > 2$.
 \end{lemma}
 \begin{proof}
  Applying Lemma \ref{commutants}, one can easily observe that for $m > 2$, $Q^m_1 = H_{[1, m-2]}$ or 
  $H_{[0, m-3]}$ according as $m$ is odd or even and consequently, $\mathcal{N}^m \subset \mathcal{M}$ is not irreducible.  
 \end{proof}

 \section{Planar algebra of $\mathcal{N}^m \subset \mathcal{M} (m > 2)$}
     Let $m > 2$ be an integer. In this section we explicitly describe the subfactor planar algebra associated to the subfactor 
     $\mathcal{N}^m \subset \mathcal{M}$ 
     which turns out to be a planar subalgebra of $^{(m)*}\!P(H^m)$.
     
\begin{figure}[!h]
\begin{center}
\psfrag{a}{\huge $m-1$}
\psfrag{b}{\Huge $x_1$}
\psfrag{r}{\huge $x_2$}
\psfrag{c}{\huge $2m-2$}
\psfrag{e}{\Huge $Fx_k$}
\psfrag{f}{\huge $\cdots$}
\psfrag{q}{\Huge $Fx_1$}
\psfrag{z}{\huge $x_{k-1}$}
\psfrag{y}{\Huge $Fx_{k-1}$}
\psfrag{v}{\Huge $Fx_2$}
\psfrag{k}{\Huge $x_k$}
\psfrag{n}{\huge $Fx_k$}
\psfrag{X}{\huge $X$}
\psfrag{Y}{\huge $X$}

\resizebox{11cm}{!}{\includegraphics{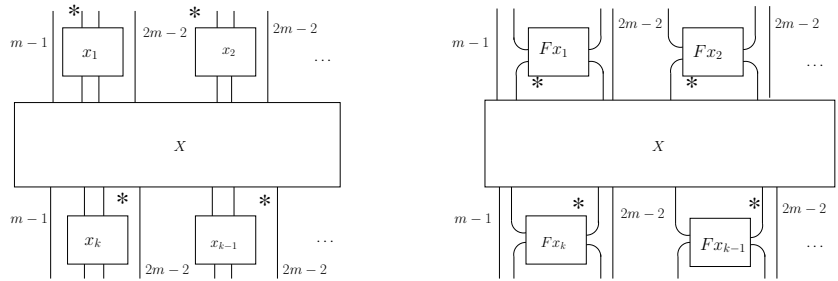}}
\end{center}
\caption{$\alpha^{m, k}_x(X)$, $m$ odd (Left) and $\alpha^{m, k}_x(X)$, $m$ even (Right)}
\label{fig:D10new}
\end{figure}

     For each $n \geq 1$, consider the linear map $\alpha^{m, n} : H \rightarrow End(P(H^m)_{mn})$ defined for $x \in H$ and 
     $X \in P(H^m)_{mn}$ by Figure \ref{fig:D10new} where the notation $\alpha^{m, n}_x$ stands for 
     $\alpha^{m, n}(x)$. 

With the help of the maps $\alpha^{m, n}$ defined above we give an equivalent description of the spaces $Q^m_n$.
    \begin{proposition}\label{descrip}
    For any $k \geq 1$,
    \begin{align*}
   Q_k^m = \{X \in P(H^m)_{mk} : \alpha^{m, k}_h(X) = X \}.
    \end{align*}
    \end{proposition}
    Before we proceed to prove Proposition \ref{descrip}, we pause for a simple Hopf algebraic lemma whose proof is similar to that of 
    \cite[Lemma 34]{Sde2018} and hence, we omit the proof.
    \begin{lemma}\label{Hopf}
    Let $k \geq 1$ be an integer.
    \begin{itemize}
    \item[(a)] If $m$ is odd, then for $X \in H_{[1, m(2k-1) - 1]}$, the following are equivalent :
    \begin{itemize}
    \item[(i)] $X \rtimes 1$ commutes with $\Delta_{k-1}(x) \in \otimes_{i = 1}^{k} H_{m(2i-1)}, \forall x \in H$,
    
    \item[(ii)]$\Delta_{k-1}(h_1)(X \rtimes 1) \Delta_{k-1}(Sh_2) = X \rtimes 1$, where 
    $\Delta_{k-1}(h_1) \otimes \Delta_{k-1}(Sh_2) \in (\otimes_{i = 1}^k H_{m(2i-1)})^{\otimes 2}$.   
    \end{itemize}
    \item[(b)] If $m$ is odd, then for $X \in H_{[1, 2mk - 1]}$, the following are equivalent : 
    \begin{itemize}
    \item[(i)] $X$ commutes with $\Delta_{k-1}(x) \in \otimes_{i = 1}^{k} H_{m(2i-1)}, \forall x \in H$,
    
    \item[(ii)]$\Delta_{k-1}(h_1) X \Delta_{k-1}(Sh_2) = X$, where 
    $\Delta_{k-1}(h_1) \otimes \Delta_{k-1}(Sh_2) \in (\otimes_{i = 1}^k H_{m(2i-1)})^{\otimes 2}$.
    \end{itemize}
     \item[(c)]  If $m$ is even, then for $X \in H_{[0, m(2k-1) - 2]}$, the following are equivalent :    
    \begin{itemize}
    \item[(i)] $X \rtimes 1$ commutes with $\Delta_{k-1}(x) \in \otimes_{i = 1}^{k} H_{m(2i-1)-1}, \forall x \in H$,
    
    \item[(ii)]$\Delta_{k-1}(h_1)(X \rtimes 1) \Delta_{k-1}(Sh_2) = X \rtimes 1$, where 
    $\Delta_{k-1}(h_1) \otimes \Delta_{k-1}(Sh_2) \in (\otimes_{i = 1}^k H_{m(2i-1)-1})^{\otimes 2}.$
    \end{itemize}
     \item[(d)]  If $m$ is even, then for $X \in H_{[0, 2mk - 2]}$, the following are equivalent :    
    \begin{itemize}
    \item[(i)] $X$ commutes with $\Delta_{k-1}(x) \in \otimes_{i = 1}^{k} H_{m(2i-1)-1}, \forall x \in H$,
    
    \item[(ii)]$\Delta_{k-1}(h_1) X \Delta_{k-1}(Sh_2) = X$, where 
    $\Delta_{k-1}(h_1) \otimes \Delta_{k-1}(Sh_2) \in (\otimes_{i = 1}^k H_{m(2i-1)-1})^{\otimes 2}.$    
    
    \end{itemize} 
    \end{itemize}
    \end{lemma}
   We are now ready to prove Proposition \ref{descrip}.
      \begin{proof}[Proof of Proposition \ref{descrip}]
     When $m > 2$ is even, the proof of the proposition is similar to that of \cite[Proposition 33]{Sde2018}. Thus, we prove the 
     proposition only when $m > 2$ is odd, leaving the other case for the reader.
      It is an immediate consequence of Lemma \ref{Hopf}(a) that the space $Q^m_{2k-1}$ can equivalently be described as
      \begin{align*}
       Q^m_{2k-1} = \{ X \in H_{[1, m(2k-1)-1]}: \Delta_{k-1}(h_1) (X \rtimes 1) \Delta_{k-1}(Sh_2) = X \rtimes 1 \}
      \end{align*}
      where $\Delta_{k-1}(h_1) \otimes \Delta_{k-1}(h_2) \in (\otimes_{i=1}^k H_{m(2i-1)})^{\otimes 2}$. 
      Interpreting this equivalent description of $Q^m_{2k-1}$ in the language of the planar algebra of $H$,
      we note that $Q^m_{2k-1}$ consists of
      precisely those elements $X \in P(H)_{m(2k-1)}$ such that the equation of Figure \ref{fig:pic41} holds.
      Now applying the conditional expectation tangle $E^{m(2k-1)}_{m(2k-1)+1}$, we reduce the element on the left in Figure \ref{fig:pic41}
      to that on the
      left in Figure \ref{fig:pic75}. On the other hand an application of the conditional expectation tangle $E^{m(2k-1)}_{m(2k-1)+1}$ 
      to the element on the right 
      in Figure \ref{fig:pic41} and then an appeal to the modulus relation reduces the element on the right in Figure \ref{fig:pic41} 
      to $\delta X$ as shown on the right in Figure \ref{fig:pic75}. Now applying the exchange relation first and then the modulus relation, 
      one can easily see that
      the element on the left in Figure \ref{fig:pic75} indeed equals $\delta \alpha^{m, 2k-1}_h(X)$
       and the desired description of $Q^m_{2k-1}$ follows.

      Similarly, it follows immediately from Lemma \ref{Hopf}(b) that the space $Q^m_{2k}$ can equivalently be described as
      \begin{align*}
       Q^m_{2k} = \{ X \in H_{[1, 2km-1]}: \Delta_{k-1}(h_1) X \Delta_{k-1}(Sh_2) = X \}
      \end{align*}
      where $\Delta_{k-1}(h_1) \otimes \Delta_{k-1}(h_2) \in (\otimes_{i=1}^k H_{m(2i-1)})^{\otimes 2}$.
      Now the desired description of $Q^m_{2k}$ follows at once from the definition of $\alpha^{m, 2k}_h(X)$ and by interpreting this 
      equivalent 
      description of $Q^m_{2k}$ in the language of $P(H)$,
      completing the proof.
      
        \begin{figure}[!h]
\begin{center}
\psfrag{k}{\huge $m-1$}
\psfrag{l}{\Huge $h_1$}
\psfrag{m}{\huge $2m-2$}
\psfrag{n}{\Huge $h_{k-1}$}
\psfrag{1}{\huge $m-1$}
\psfrag{u}{\huge $2m-2$}
\psfrag{p}{\huge $2m-2$}
\psfrag{q}{\Huge $h_{k+2}$}
\psfrag{r}{\huge $2m-2$}
\psfrag{s}{\Huge $h_{2k}$}
\psfrag{t}{\huge $m-1$}
\psfrag{2}{\huge $2m-2$}
\psfrag{Y}{\Huge $X$}
\psfrag{a}{\Huge $h_k$}
\psfrag{b}{\Huge $h_{k+1}$}
\psfrag{1}{\Huge $1$}
\psfrag{e}{\Huge $=$}
\resizebox{8.5cm}{!}{\includegraphics{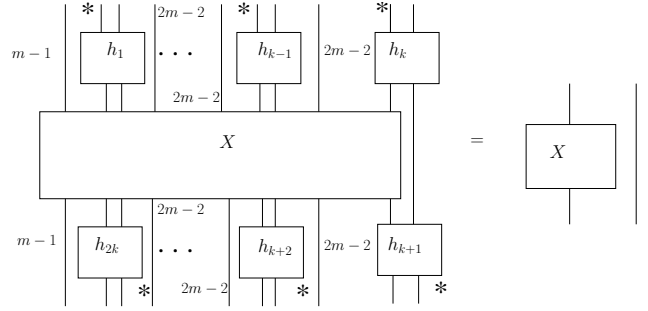}}
\end{center}
\caption{Characterisation of $X$}
\label{fig:pic41}
\end{figure} 
 \begin{figure}[!h]
\begin{center}
\psfrag{k}{\huge $m-1$}
\psfrag{l}{\Huge $h_1$}
\psfrag{m}{\huge $2m-2$}
\psfrag{n}{\Huge $h_{k-1}$}
\psfrag{1}{\huge $m-1$}
\psfrag{u}{\huge $2m-2$}
\psfrag{p}{\huge $2m-2$}
\psfrag{q}{\Huge $h_{k+2}$}
\psfrag{r}{\huge $2m-2$}
\psfrag{s}{\Huge $h_{2k}$}
\psfrag{t}{\huge $m-1$}
\psfrag{2}{\huge $2m-2$}
\psfrag{Y}{\Huge $X$}
\psfrag{a}{\Huge $h_k$}
\psfrag{b}{\Huge $h_{k+1}$}
\psfrag{d}{\Huge $\delta$}
\psfrag{e}{\Huge $=$}

\resizebox{9cm}{!}{\includegraphics{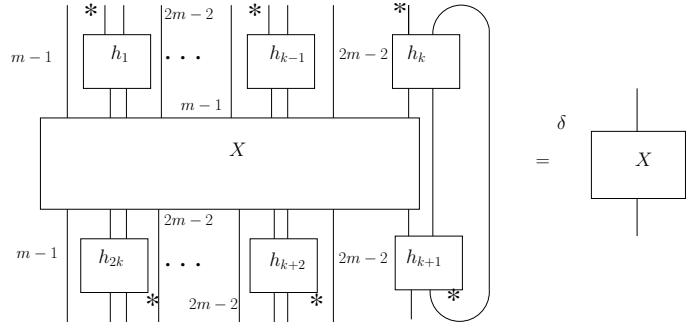}}
\end{center}
\caption{Equivalent characterisation of $X$}
\label{fig:pic75}
\end{figure}     
\end{proof}
  Thus for each $m > 2$, we have a family $\{Q^m_n : n \geq 1 \}$ of vector spaces where for $n \geq 1, Q^m_n$ 
  is a subspace of  $P(H^m)_{mn} = ^{(m)}\!\!\!P(H^m)_{n}$. Setting $Q^m_{0, \pm} = \mathbb{C}$, we note that 
  $Q^m := \{Q_n^m: n \in Col \}$ is a subspace of $^{(m)}\!P(H^m)$. The following proposition, whose proof is similar to that of
  , shows that $Q^m$ is indeed a planar subalgebra of 
  $^{(m)}\!P(H^m)$ and we omit its proof.
   \begin{proposition}\label{planar}    
    For $m > 2$, $Q^m$ is a planar subalgebra of  $^{(m)}\!P(H^m)$. 
    \end{proposition}
     \begin{proof}
     By an appeal to Theorem \cite[Theorem 3.5]{KdySnd2004}, it suffices to prove that $Q^{m}$ is closed under the action of the following  set of tangles
     \begin{align*}
      \{1^{0,+}, 1^{0, -}\} \cup \{R_k^k:k \geq 2\} \cup \{M_{k, k}^k, E^k_{k+1}, I_k^{k+1}: k \in Col \}
     \end{align*}
     where we refer to \cite[Figures 2, 3 and 5]{Sde2018} for the definition of tangles $M_{k, k}^k, E^k_{k+1}, I_k^{k+1}, 
     1^{0,+}$ and $1^{0, -}$.
     
    When $m$ is even, the proof of the proposition is similar to that of \cite[Proposition 35]{Sde2018}. Thus we prove the result only when $m$ 
    is odd.
    
    It is obvious to see that $Q^m$ is closed under the action of the tangles $1^{0, \pm}$ and $M_{k, k}^k, I_k^{k+1} (k \in Col)$.

     To see that $Q^m$ is closed under the action of the rotation tangle $R_k^k \ (k \geq 2)$, we note that for any 
     $X \in Q^m_k \ (k \geq 2)$, we have 
     \begin{align*}
      Z^{^{(m)}\!P(H)}_{R_k^k}(X) =   Z^{^{(m)}\!P(H)}_{R_k^k}(\alpha^{m, k}_h(X)) = \alpha^{m, k}_h(  Z^{^{(m)}\!P(H)}_{R_k^k}(X))  
     \end{align*}
     where the first equality follows from the fact that $\alpha^{m, k}_h(X) = X$ and to see the second equality we need to use the Hopf algebra
     identity $h_1 \otimes h_2 \otimes \cdots \otimes h_l = h_2 \otimes h_3 \otimes 
     \cdots \otimes h_l \otimes h_1$ ($l \geq 2$) which basically follows from $h_1 \otimes h_2 = h_2 \otimes h_1$ (which essentially 
     expresses traciality of $h$).
     
     Verifying that $Q^m$ is closed under the action of $E^k_{k+1} (k \geq 1)$ amounts to verification of the 
     following identity  
     \begin{align*}
      Z_{E^k_{k+1}}^{^{(m)}\!P(H)}(X) = \alpha^{m, k}_h(Z_{E^k_{k+1}}^{^{(m)}\!P(H)}(X))
     \end{align*}
     for any $X \in Q^m_{k+1}$.
     Note that since $\alpha^{m, k+1}_h(X) = X$, we have that 
     \begin{align*}
      Z_{E^k_{k+1}}^{^{(m)}\!P(H)}(X) = Z_{E^k_{k+1}}^{^{(m)}\!P(H)}(\alpha^{m, k+1}_h(X))
     = Z_{{E^k_{k+1}}^{(m)}}^{P(H)}(\alpha^{m, k+1}_h(X)).
     \end{align*}
 When $k$ is odd (resp., even), representing the element 
     $Z_{{E^k_{k+1}}^{(m)}}^{P(H)}(\alpha^{m, k+1}_h(X))$ pictorially in $P(H)$ and then applying relation (E) (resp., (C)) one can easily see 
     that 
     \begin{align*}
    Z_{{E^k_{k+1}}^{(m)}}^{P(H)}(\alpha^{m, k+1}_h(X)) = \alpha^{m, k}_h(Z_{E^k_{k+1}}^{^{(m)}\!P(H)}(X)),  
     \end{align*}
     finishing the proof. 
     \end{proof}
     As an immediate corollary of Proposition \ref{planar} we obtain that
    \begin{corollary} \label{planar1}
     $^*\!Q^m$, the adjoint of $Q^m$, is a planar subalgebra of  $^{*(m)}\!P(H^m)$.
    \end{corollary}
    Finally, similar argument as in the proof of \cite[Proposition 36]{Sde2018} shows that
    $P^{\mathcal{N}^m \subset \mathcal{M}}$, the planar algebra associated to $\mathcal{N}^m \subset \mathcal{M}$,
      is given by the adjoint of the planar algebra $Q^m$. Thus we have:    
    \begin{proposition}\label{pl alg}
     $P^{\mathcal{N}^m \subset \mathcal{M}} = {^*\!Q^m}$, $m > 2$.
    \end{proposition}
     We collect the results of the previous statements into a single main theorem.
    \begin{theorem}\label{maintheorem}
     For any integer $m > 2$,  
     $^*\!Q^m$ is a planar subalgebra of $^{*(m)}\!P(H^m)$ and $^*\!Q^m = P^{\mathcal{N}^m \subset \mathcal{M}}$. 
      If $m$ is odd, $^*\!Q^m_k$ ($k \geq 1$) consists of all $X \in 
      P(H)_{mk}$ such that the element on the left in Figure \ref{fig:charac2} equals $X$        
      and if $m$ is even, $^*\!Q^m_k$ ($k \geq 1$) consists of all $X \in 
      P(H^*)_{mk}$ such that the element on the right in Figure \ref{fig:charac2} equals $X$.
       \begin{figure}[!h]
\begin{center}
\psfrag{d}{\Huge $h_1$}
\psfrag{e}{\Huge $h_2$}
\psfrag{f}{\Huge $h_3$}
\psfrag{g}{\Huge $h_{k-1}$}
\psfrag{h}{\Huge $h_k$}
\psfrag{x}{\Huge $m-1$}
\psfrag{y}{\Huge $m-1$}
\psfrag{a}{\Huge $2m-2$}
\psfrag{z}{\Huge $2m-2$}
\psfrag{c}{\Huge $2m-2$}
\psfrag{7}{\Huge $2m-2$}
\psfrag{s}{\Huge $2m-2$}
\psfrag{m}{\Huge $Fh_{k-1}$}
\psfrag{n}{\Huge $Fh_k$}
\psfrag{o}{\Huge $Fh_1$}
\psfrag{p}{\Huge $Fh_2$}
\psfrag{q}{\Huge $Fh_3$}
\psfrag{r}{\Huge $Fh_{k-1}$}
\psfrag{s}{\Huge $Fh_k$}
\psfrag{m}{\Huge $m-1$}
\psfrag{n}{\Huge $m-1$}
\psfrag{t}{\Huge $2m-2$}
\psfrag{u}{\Huge $2m-2$}
\psfrag{v}{\Huge $2m-2$}
\psfrag{X}{\Huge $X$}
\resizebox{12cm}{!}{\includegraphics{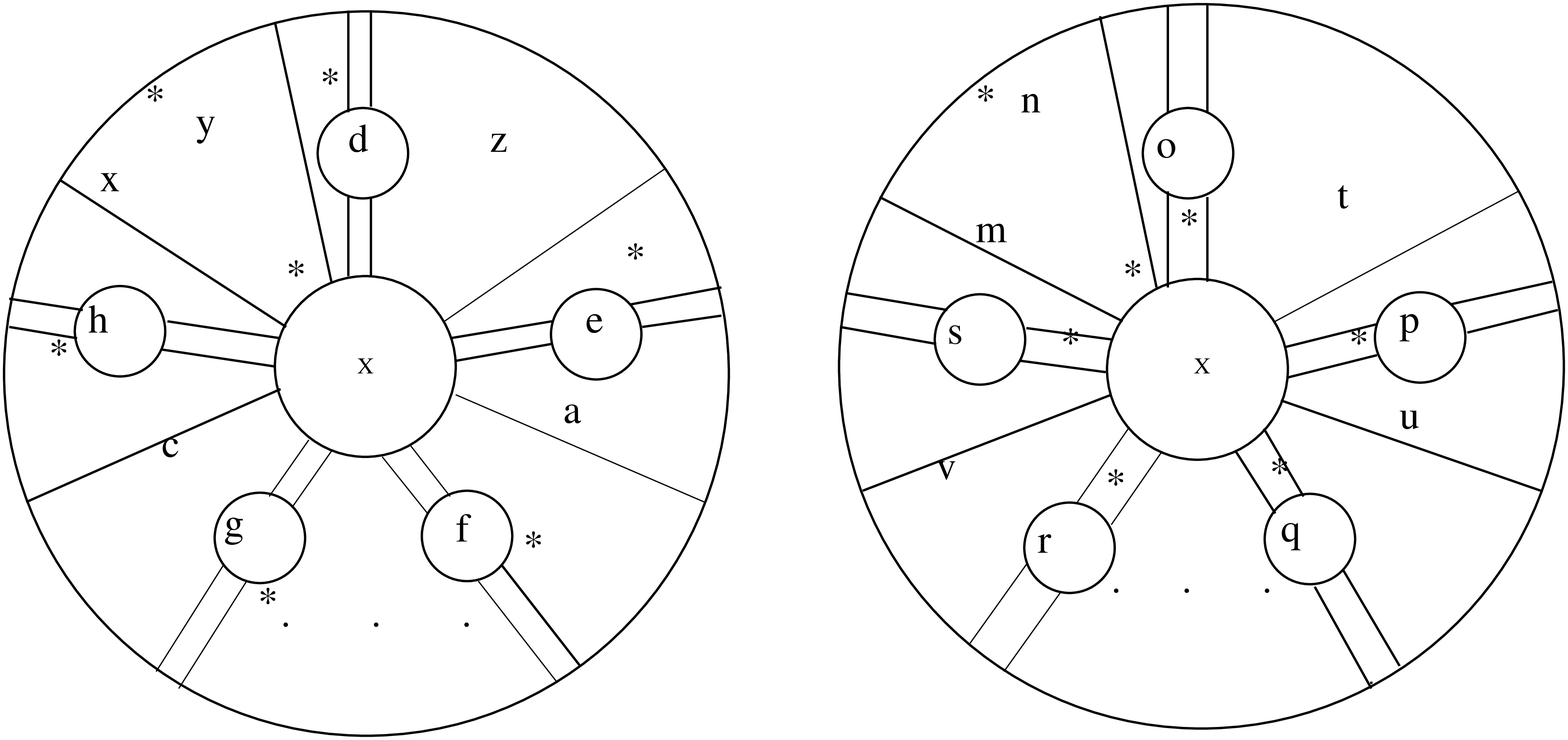}}
\end{center}
\caption{}
\label{fig:charac2}
\end{figure}     

    \end{theorem}
    \begin{proof}
     It follows immediately from Proposition \ref{pl alg} after observing that $\alpha^{m, k}_h(X)$ for $m$ odd 
     (resp., even)
     in Figure \ref{fig:D10new} is equivalent to the element on the left (resp., right)
     in Figure \ref{fig:charac2}.  
    \end{proof}
  \section{Depth of $\mathcal{N}^m \subset \mathcal{M}, m > 2$ }
  In this section we investigate the depth of the subfactors $\mathcal{N}^m \subset \mathcal{M}$ for $m > 2$.
  The main result of this section is contained in the following theorem.
  \begin{theorem}\label{depth}
   For $m > 2$, the subfactor $\mathcal{N}^m \subset \mathcal{M}$ is of depth 2.
  \end{theorem}
  By virtue of the commutative diagram in Figure \ref{fig:D9}, one can easily see that 
   $\mathcal{N}^{m} \subset \mathcal{M}$ has depth $k$, $k \geq 1$ an integer, is equivalent to
$k$ being the smallest positive integer such that 
$Q_{k-1}^{m^{op}} \subset Q_k^{m^{op}} \subset Q_{k+1}^{m^{op}}$ is an instance of the basic construction with the Jones projection $e^m_{k+1}$
which obviously is equivalent to saying that $Q^m_{k-1} \subset Q^m_k \subset Q^m_{k+1}$ is an instance of the basic construction with the same Jones 
projection. Thus, in order to prove Theorem \ref{depth}, it suffices to show that $2$ is the smallest positive integer such that
$Q^m_1 \subset Q^m_2 \subset Q^m_3$ is an instance of the basic construction with the Jones projection $e^m_3$.

We find it necessary to explicitly know the elements of the space $Q^m_2$. The following lemma is the main step to this end.
 \begin{lemma}\label{element}
  Let $\mathcal{S}$ be the space defined by
  \begin{align*}
 \mathcal{S}  = \{X \in H^* \rtimes H \rtimes H^* (\cong P(H^*)_4) :
   X \ \mbox{commutes with} \ \epsilon \rtimes x \rtimes \epsilon, \forall x \in H \}.
  \end{align*}
Then $\mathcal{S}$ precisely consists of elements of the form $Z_A^{P(H^*)}(f_2 \otimes f_1 \otimes g)$ 
  where $A \in {\mathcal{T}}{(4)}$ is the tangle as shown on the left in Figure \ref{fig:pic43} and $f \otimes g \in H^* \otimes H^*$. 
  Consequently, $\mathcal{S}$ has dimension $(dim ~H)^2$.
 \end{lemma}
 \begin{proof}
    Let $X = Z_A^{P(H^*)}(f_2 \otimes f_1 \otimes g) \in P(H^*)_4$ where $f \otimes g \in H^* \otimes H^*$. 
    For any $t \in H,$ let $Y_t$ denote the image of $\epsilon\rtimes t \rtimes
    \epsilon \in H^* \rtimes H \rtimes H^*$ in $P(H^*)_4$ under the algebra isomorphism between $H^* \rtimes H \rtimes H^*$ and $P(H^*)_4$ as
    given by Lemma \ref{iso}, i.e., $Y_t = Z_{T^4}^{P(H^*)}(\epsilon \otimes Ft \otimes \epsilon)$ (see Figure \ref{fig:c9}) where
    $T^4$ is the tangle of colour $4$ as shown on the right in Figure \ref{fig:pic599}. 
    Thus given $t \in H$, we need to show that $X$ commutes with $Y_t$. The element $Y_t X$ (resp., $X Y_t$) is shown on the left (resp., right)
    in Figure \ref{fig:c9}.
     \begin{figure}[!h]
     \begin{center}
 \psfrag{a}{\Huge $f_2$}
 \psfrag{b}{\Huge $f_1$}
 \psfrag{c}{\Huge $g$}
 \psfrag{d}{\Huge $Ft$}
 \resizebox{8.8cm}{!}{\includegraphics{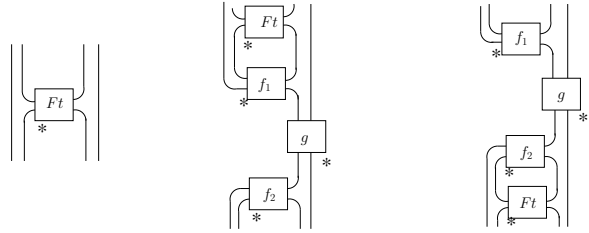}}
 \end{center}
 \caption{Elements: $Y_t$ (left), $Y_t X$ (middle) and $XY_t$ (right)}
 \label{fig:c9}
 \end{figure}
    Set $k = Ft$. An application of the relations (E) and (T) shows
    that $XY_t$ equals the element
    $\delta (Sf_2 k)(h) Z_A^{P(H^*)}(f_3 \otimes f_1 \otimes g)$ whereas $Y_tX$ equals the element 
    $\delta(Sk_1 f_1)(h) Z_A^{P(H^*)}(f_2 \otimes k_2 \otimes g)$.
    Since, by Lemma \ref{iso1}, 
    $Z_A^{P(H^*)}$ is a linear isomorphism of $H^{*\otimes3}$ onto $P(H^*)_4$, in order to see that $XY_t = Y_tX$, it suffices to verify that
    \begin{align*}
    \delta (Sf_2 k)(h) f_3 \otimes f_1 \otimes g = \delta(Sk_1 f_1)(h) f_2 \otimes k_2 \otimes g.
    \end{align*}
    Evaluating the expression on the left-hand side on $a \otimes b \otimes c \in H \otimes H \otimes H$ we obtain 
    $f(b Sh_1 a) k(h_2) g(c)$ whereas evaluating the 
    expression on the right-hand side on the same element we obtain the value $k(Sh_1b) f(h_2a) g(c)$ 
    which, using the Hopf-algebraic formula $Sh_1x \otimes h_2 = Sh_1 \otimes xh_2$ and the fact that $Sh = h$, equals 
    $f(b Sh_1 a) k(h_2) g(c)$.
    Thus we see that
    \begin{align*}
     \{Z_A^{P(H^*)}(f_2 \otimes f_1 \otimes g) : f \otimes g \in H^{*\otimes 2} \} \subseteq \mathcal{S}
    \end{align*}
    and consequently the dimension of the space $S$ is $\geq (dim ~H)^2$. To finish the proof we just need to see that the dimension of 
    $\mathcal{S}$ is $\leq (dim ~H)^2$. First observe that for any $X \in \mathcal{S}$,
    \begin{align*}
    X = X (\epsilon \rtimes h_1 \rtimes \epsilon) (\epsilon \rtimes Sh_2 \rtimes \epsilon) = 
    (\epsilon \rtimes h_1 \rtimes \epsilon) X (\epsilon \rtimes Sh_2 \rtimes \epsilon)
    \end{align*}
    and thus $\mathcal{S}$ is a subspace of 
    \begin{align*}
     W = \{ X \in P(H^*)_4 : X =  Y_{h_1} X Y_{Sh_2}\}.
    \end{align*}
    Thus it suffices to show that $dim ~W \leq (dim ~H)^2$. Let us consider the tangle $S \in {\mathcal{T}}{(4)}$ as given by Figure \ref{fig:pic15}. 
     \begin{figure}[!h]
\begin{center}
\psfrag{a}{\Huge $1$}
\psfrag{b}{\Huge $2$}
\psfrag{c}{\Huge $3$}
\resizebox{2.0cm}{!}{\includegraphics{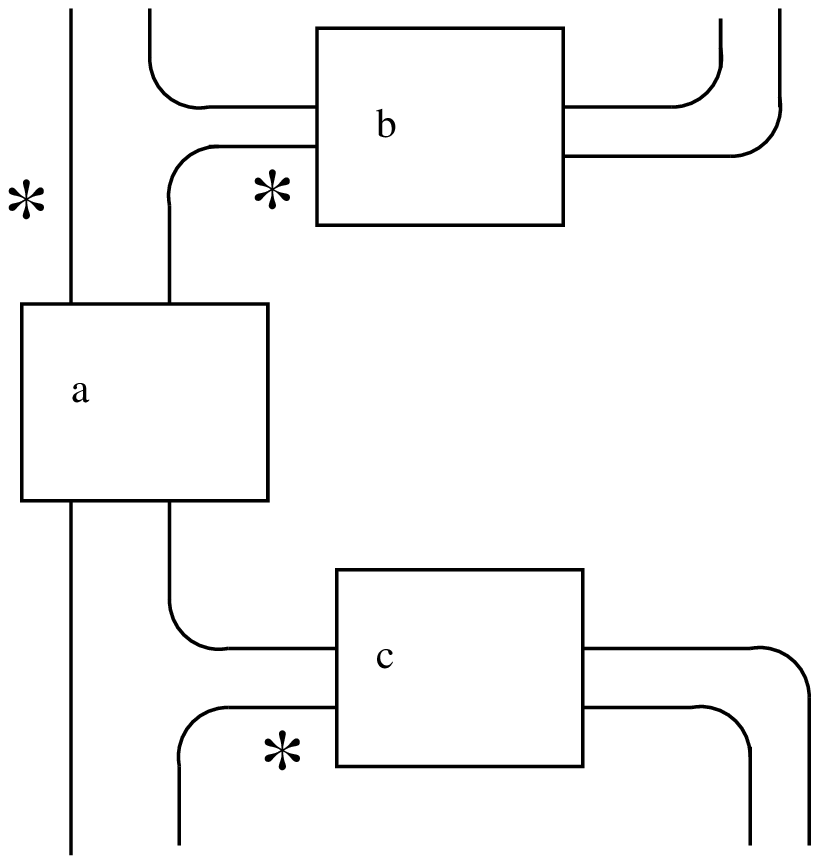}}
\end{center}
\caption{Tangle $S$}
\label{fig:pic15}
\end{figure} 
Since $Z_S$ induces, by Lemma \ref{iso1}, linear isomorphism of $H^{*\otimes 3}$ onto $P(H^*)_4$, we conclude that
 $dim ~W = dim ~\{x \otimes y \otimes z \in H^{\otimes 3} : Z_S(Fx \otimes Fy \otimes Fz) = Y_{h_1} Z_S(Fx \otimes Fy \otimes Fz) Y_{Sh_2} \}.$
One can easily see that $Y_{h_1} Z_S(Fx \otimes Fy \otimes Fz) Y_{Sh_2} = Z_S(Fx \otimes F(h_1 y) \otimes F(z Sh_2))$ and consequently,
$dim ~W = dim ~U$ where
\begin{align*}
 U = \{x \otimes y \otimes z \in H^{\otimes 3} : x \otimes y \otimes z = x \otimes h_1 y \otimes z Sh_2 \}.
\end{align*}
Thus it suffices to see that $dim ~U \leq (dim ~H)^2$. Note that the space $U$, using the Hopf-algebraic formula
$h_1 a \otimes b Sh_2 = h_1 \otimes b a Sh_2$, can alternatively described as 
\begin{align*}
  U = \{x \otimes y \otimes z \in H^{\otimes 3} : x \otimes y \otimes z = x \otimes h_1 \otimes z y Sh_2 \}.
\end{align*}
Finally observe that $U$ is contained in the image of the injective linear map $\theta$ form $H \otimes H$ to $H \otimes H \otimes H$ given by
$x \otimes y \rightarrow x \otimes h_1 \otimes y Sh_2$ for if $x \otimes y \otimes z \in U$, then clearly $x \otimes y \otimes z = 
\theta(x \otimes zy)$ and hence, $dim ~U \leq (dim ~H)^2$, finishing the proof.
\end{proof} 
We now present a technical lemma that will be useful in order to precisely express the elements of $Q^m_2$, $m > 2$.
\begin{lemma}\label{eq}
The following equation holds in $P(H^*)_4$ for $f, g \in H^* :$
\begin{align*}
 Z_A^{P(H^*)}(f_2 \otimes f_1 \otimes g) = Z_{T^4}^{P(H^*)}(f_1Sg_3Sf_3 \otimes f_2g_2 \otimes Sg_1).
\end{align*}
\end{lemma}
\begin{proof}
 Left as an exercise.
\end{proof}
Consequently, the space $\mathcal{S}$ can equivalently be described as: 
\begin{corollary}\label{Eq}
 $\mathcal{S} = \{ f_1Sg_3Sf_3 \rtimes F^{-1}(f_2g_2) \rtimes Sg_1 \in H^* \rtimes H \rtimes H^*: f \otimes g \in H^{*\otimes2} \}.$
\end{corollary}
\begin{proof}
 Follows immediately from Lemmas \ref{element} and \ref{eq}.
\end{proof}
\begin{corollary}\label{Eq1}
 Let $m > 2$. 
 \begin{itemize}
  \item [(i)] If $m$ is odd, then $Q^m_2$ consists of elements of the form 
 \begin{equation}\label{Eq2}
    x^1 \rtimes f^2 \rtimes \cdots \rtimes x^{m-2} \rtimes g_1 Sk_3 Sg_3 \rtimes
  F^{-1}(g_2 k_2) \rtimes Sk_1 \rtimes x^{m+2} \rtimes f^{m+3} \rtimes \cdots \rtimes x^{2m-1} \in A(H)_{2m-1}
  \end{equation}
  with $(x^1 \otimes \cdots \otimes
  x^{m-2} \otimes x^{m+2} \otimes \cdots \otimes x^{2m-1}) \otimes (f^2 \otimes \cdots \otimes f^{m-3} \otimes g \otimes k \otimes f^{m+3} 
  \otimes \cdots \otimes f^{2m-2}) \in H^{\otimes(m-1)} \otimes H^{* \otimes (m-1)}$.
  \item [(ii)] If $m$ is even, $Q^m_2$ consists of elements of the form 
  \begin{equation}\label{Eq3}
   f^1 \rtimes x^2 \rtimes \cdots \rtimes x^{m-2} \rtimes g_1 Sk_3 Sg_3 \rtimes
  F^{-1}(g_2 k_2) \rtimes Sk_1 \rtimes x^{m+2} \rtimes \cdots \rtimes f^{2m-1} \in A(H^*)_{2m-1} 
  \end{equation}
  with $(f^1 \otimes \cdots \otimes f^{m-3} \otimes g \otimes k \otimes f^{m+3} \otimes \cdots \otimes f^{2m-1}) \otimes 
  (x^2 \otimes \cdots x^{m-2} \otimes
  x^{m+2} \otimes \cdots \otimes x^{2m-2}) \in H^{* \otimes m} \otimes H^{\otimes (m-2)}$ . In this case, the elements of $Q^{m}_2$
  can equivalently be expressed as
 \begin{equation}\label{Eq4}
   Z_{A(m-2, m-2)}^{P(H^*)}(f^1 \otimes f^2 \otimes \cdots \otimes f^{m-2} \otimes f^{m-1}_2 \otimes 
  f^{m-1}_1 \otimes f^m \otimes f^{m+1} \otimes \cdots \otimes f^{2m-2}) 
  \end{equation}
 with $\otimes_{i = 1}^{2m-2} f^i \in H^{*\otimes(2m-2)}$ (see Figure \ref{fig:pic43} for the definition of tangles $A(m-2, m-2)$). 
 \end{itemize}  
 \end{corollary}
\begin{proof}
 \begin{itemize}
  \item[(i)] Follows from the definition of $Q^m_2$ as given in $\S 4$ along with an application of Corollary \ref{Eq}.
  \item[(ii)] Follows from the definition of $Q^m_2$ as given in $\S 4$ together with an appeal to Corollary \ref{Eq} and Lemma \ref{eq}. 
\end{itemize}
\end{proof}
\begin{remark}\label{imp}
 Note that $^*\!Q^m_2 = Q^m_2$ as vector spaces, only the multiplication in $^*\!Q^m_2$ is opposite to that of $Q^m_2$. We make 
 no notational distinction between the elements of $^*\!Q^m_2$ and $Q^m_2$. Thus, a general element of $^*\!Q^m_2$ will always be
 expressed in the form as given by \eqref{Eq2} when $m$ is odd, or in the form as given by \eqref{Eq3} or \eqref{Eq4} when $m$ is even. 
\end{remark} 
Note that $\mathcal{N}^m \subset \mathcal{M}$ can not be of depth $1$ for otherwise it must happen that $\mathbb{C} \subset 
Q^m_1 \subset Q^m_2$ is an instance of the basic construction and therefore, $dim ~Q^m_2$ must be equal to
$dim ~(Q^m_1)^2$ which is not possible since $ dim ~Q^m_1 = (dim ~H)^{m-2}$ whereas 
$dim ~Q^m_2$, by an appeal to Corollary \ref{Eq1}, equals $(dim ~H)^{2(m-1)}$. Consequently, depth of $\mathcal{N}^m \subset \mathcal{M}$
is greater than $1$.

In the following two propositions, namely, Proposition \ref{even} and Proposition \ref{odd}, we prove that $Q^m_1 \subset Q^m_2 \subset Q^m_3$
is an instance of the basic construction where Proposition \ref{even} treats the case when $m$ is even while Proposition \ref{odd}
treats the case when $m$ is odd.
We will use Lemma \ref{iso1} frequently without any mention in the proofs of both the propositions.

\begin{proposition}\label{even}
If $m > 1$, then $Q^{2m}_1 \subset Q^{2m}_2 \subset Q^{2m}_3$ is an instance of the basic construction with the Jones projection $e^{2m}_3$.
\end{proposition}
\begin{proof}
 For notational convenience we use the symbol $e$ to denote $e^{2m}_3$.
 Since the conditions (i) and (ii) of Lemma \ref{basic} are automatically satisfied, we only need to verify the condition (iii). 
 Let us consider the elements of $Q^{2m}_2$ given by 
 \begin{align*}
  &X = Z_{A(2m-2, 2m-2)}^{P(H^*)}(f^0 \otimes f^1 \otimes \cdots \otimes f^{2m-3} \otimes f^{2m-2}_2 \otimes f^{2m-2}_1 \otimes 
f^{2m-1} \otimes \cdots \otimes f^{4m-3}), \\
&Y = Z_{A(2m-2, 2m-2)}^{P(H^*)}((\epsilon \otimes F(1))^{\otimes (m-1)} \otimes g^{2m-2}_2 \otimes g^{2m-2}_1 
\otimes g^{2m-1} \otimes g^{2m} \otimes \cdots \otimes g^{4m-3}).
 \end{align*}
A simple computation shows that $XeY$ equals the element 
\begin{align*}
\delta^{-2m} Z_S^{P(H^*)}(& f^0 \otimes \cdots \otimes f^{2m-3} \otimes f^{2m} \otimes \cdots
\otimes f^{4m-3} \otimes g^{2m} \otimes \cdots \otimes g^{4m-3} \otimes g^{2m-1} \otimes f^{2m-1}\\
& \otimes f^{2m-2}_1 \otimes Sg^{2m-2}_1 f^{2m-2}_2 \otimes g^{2m-2}_2)
\end{align*}
where $S \in \mathcal{T}(6m)$ is as shown in Figure \ref{fig:pic42}. Consider the linear map 
$\theta : H^{*\otimes 2} \rightarrow H^{*\otimes 3}$ given by 
\begin{align*}
 f \otimes g \mapsto f_1 \otimes Sg_1 f_2 \otimes g_2.
\end{align*}
 \begin{figure}[!h]
\begin{center}
\psfrag{a}{\huge $1$}
\psfrag{b}{\huge $2$}
\psfrag{c}{\huge $2m-2$}
\psfrag{d}{\huge $6m-3$}
\psfrag{e}{\huge $6m-4$}
\psfrag{f}{\huge $6m-2$}
\psfrag{g}{\huge $2m-1$}
\psfrag{h}{\huge $2m$}
\psfrag{i}{\huge $4m-4$}
\psfrag{j}{\huge $6m-1$}
\psfrag{k}{\huge $6m-5$}
\psfrag{l}{\huge $4m-3$}
\psfrag{m}{\huge $4m$}
\psfrag{n}{\huge $6m-6$}
\resizebox{8.0cm}{!}{\includegraphics{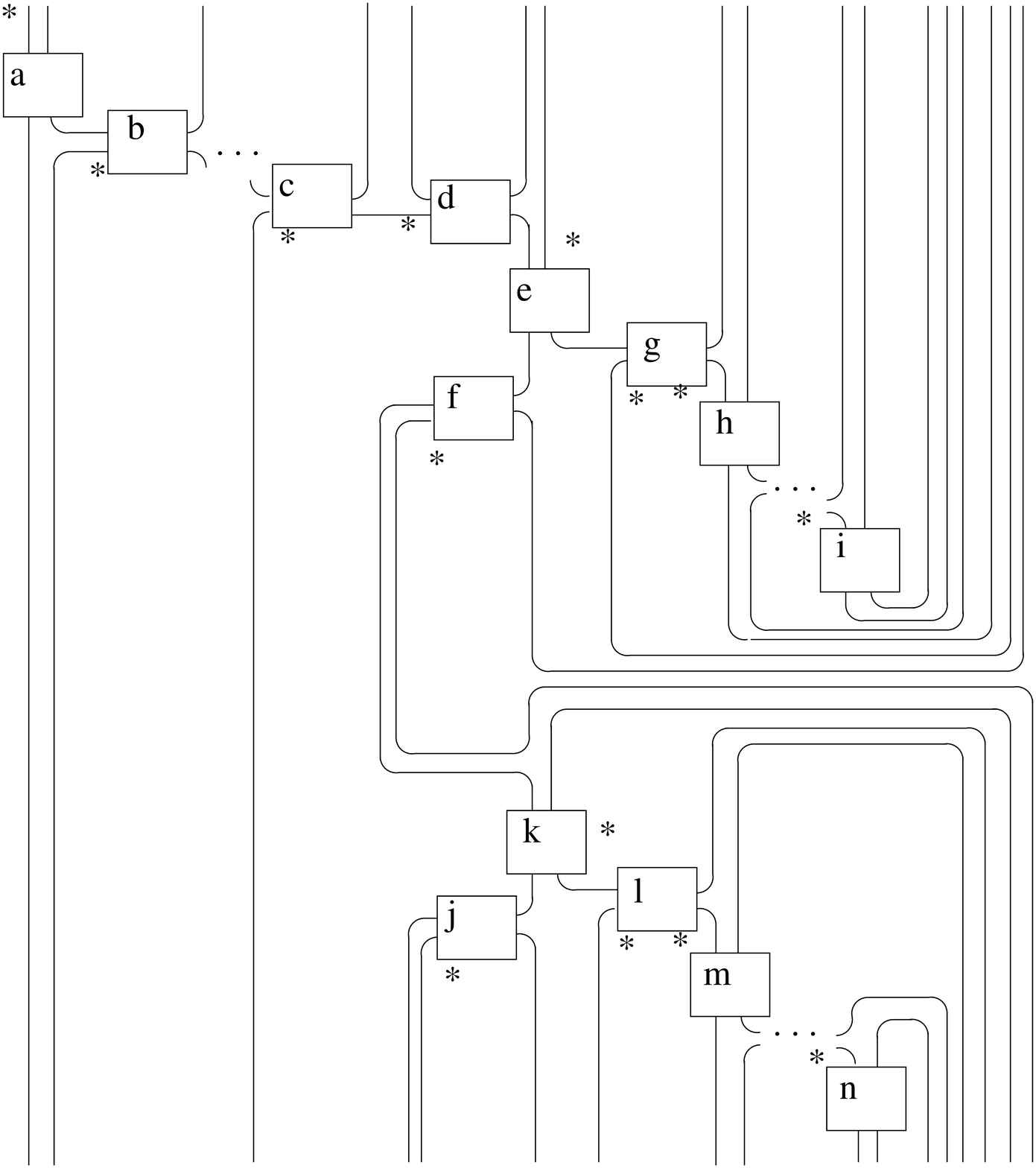}}
\end{center}
\caption{Tangle $S$}
\label{fig:pic42}
\end{figure} 
We assert that $\theta$ is injective. To see this one can easily verify that the map from $H^{*\otimes 3} \rightarrow H^{*\otimes 2}$ given by  
\begin{align*}
 f \otimes g \otimes k \mapsto f(1) k_1 g \otimes k_2
\end{align*}
is a left inverse of $\theta$, proving the assertion. Thus clearly $Q^{2m}_2 e Q^{2m}_2$ contains the image of the injective linear map
\begin{align*}
Z_S \circ ({Id_{H^*}}^{\otimes(6m-4)} \otimes \theta) : H^{*\otimes (6m-2)} \rightarrow P(H^*)_{6m} 
\end{align*}
and consequently, 
$dim ~(Q^{2m}_2 e Q^{2m}_2) \geq (dim ~H)^{6m-2}$. Thus in order to see that $Q^{2m}_2 e Q^{2m}_2 = Q^{2m}_3$ we just need to show that 
$dim ~Q^{2m}_3 \leq (dim ~H)^{6m-2}$. To this end we consider the tangle $P \in \mathcal{T}(6m)$ as shown in Figure \ref{fig:pic44}.
\begin{figure}[!h]
\begin{center}
\psfrag{1}{\huge $1$}
\psfrag{2}{\huge $2$}
\psfrag{3}{\huge $3$}
\psfrag{4}{\huge $4$}
\psfrag{5}{\huge $2m-3$}
\psfrag{6}{\huge $2m-2$}
\psfrag{7}{\huge $6m-3$}
\psfrag{8}{\huge $2m-1$}
\psfrag{9}{\huge $2m$}
\psfrag{10}{\huge $2m+1$}
\psfrag{11}{\huge $2m+2$}
\psfrag{12}{\huge $6m-7$}
\psfrag{13}{\huge $6m-6$}
\psfrag{14}{\huge $6m-5$}
\psfrag{15}{\huge $6m-4$}
\psfrag{16}{\huge $6m-2$}
\psfrag{17}{\huge $6m-1$}
\resizebox{16.0cm}{!}{\includegraphics{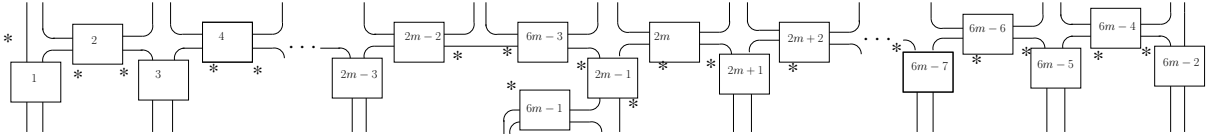}}
\end{center}
\caption{Tangle $P$}
\label{fig:pic44}
\end{figure} 
Since $Z_P^{P(H^*)}$ induces a linear isomorphism of $H^{*\otimes(6m-1)}$ onto $P(H^*)_{6m}$, we observe that, in view of Proposition 
\ref{descrip}, the space $Q^{2m}_3$ is linearly isomorphic to
\begin{align*}
 \{f^1 \otimes \cdots \otimes f^{6m-4} \otimes x \otimes y \otimes z \in H^{*\otimes(6m-4)} \otimes H^{\otimes3} :  
  \alpha_h^{2m, 3}(Z_P(f^1 \otimes \cdots \otimes f^{6m-4}\\  
  \otimes Fx \otimes Fy \otimes Fz)) =
 Z_P(f^1 \otimes \cdots \otimes f^{6m-4} \otimes Fx \otimes Fy \otimes Fz)\}. 
 \end{align*}
 A trivial computation in $P(H^*)$ shows that
 \begin{align*}
   \alpha_h^{2m, 3}(Z_P^{P(H^*)}(\otimes_{i=1}^{6m-4}f^i \otimes Fx \otimes Fy \otimes Fz)) = 
   Z_P^{P(H^*)}(\otimes_{i=1}^{6m-4}f^i \otimes F(h_1x) \otimes F(h_2y) \otimes F(h_3z)).
 \end{align*}
 Now using injectivity of $Z_P$ and invertibility of $F$, we conclude that $Q^{2m}_3$ is linearly isomorphic to the space $W$ defined by
 \begin{align*}
  W = \{x \otimes y \otimes z \otimes f^1 \otimes \cdots \otimes f^{6m-4} \in  H^{\otimes3} \otimes H^{*\otimes(6m-4)} :\\
   h_1x \otimes h_2y \otimes h_3z \otimes f^1 \otimes \cdots \otimes f^{6m-4} = 
 x \otimes y \otimes z \otimes f^1 \otimes \cdots \otimes f^{6m-4}  \}. 
 \end{align*}
Thus it suffices to see that $dim ~W \leq (dim ~H)^{6m-2}$. Note that the space $W$, using the Hopf-algebraic formula $h_1 a \otimes h_2 = h_1 \otimes
h_2 Sa$, can equivalently be described as
\begin{align*}
 W = \{x \otimes y \otimes z \otimes f^1 \otimes \cdots \otimes f^{6m-4} \in H^{\otimes3} \otimes H^{*\otimes(6m-4)} :\\
   h_1 \otimes h_2Sx_2y \otimes h_3Sx_1z \otimes f^1 \otimes \cdots \otimes f^{6m-4} = 
  x \otimes y \otimes z \otimes f^1 \otimes \cdots \otimes f^{6m-4}\}. 
\end{align*}
Further we observe that $W$ is contained in the range of the linear map $\rho: H^{\otimes2} \otimes H^{*\otimes(6m-4)} \rightarrow 
H^{\otimes3} \otimes H^{*\otimes(6m-4)}$ given by
\begin{align*}
 a \otimes b \otimes_{i= 1}^{6m-4}f^i \mapsto h_1 \otimes h_2 a \otimes h_3 b \otimes_{i= 1}^{6m-4}f^i, 
\end{align*}
for, if $X = x \otimes y \otimes z \otimes f^1 \otimes \cdots \otimes f^{6m-4} \in W$, then
$\rho(Sx_2y \otimes Sx_1z \otimes f^1 \otimes \cdots \otimes f^{6m-4}) = 
h_1 \otimes h_2Sx_2y \otimes h_3Sx_1z \otimes f^1 \otimes \cdots \otimes f^{6m-4} = X$
and consequently, $dim ~W \leq$ rank of $\rho \leq (dim ~H)^{6m-2}$,
completing the proof. 
\end{proof}

  \begin{proposition}\label{odd}
  Given $m > 1, Q^{2m-1}_1 \subset Q^{2m-1}_2 \subset Q^{2m-1}_3$ is an instance of the basic construction with the Jones projection $e^{2m-1}_3$.
  \end{proposition}
  \begin{proof}
   Since the conditions (i) and (ii) of Lemma \ref{basic} are 
  automatically satisfied, we just need to verify the condition (iii). 
  Note that $Q^{2m-1}_1 = H_{[1, 2m-3]}$, $Q^{2m-1}_2 = \{ X \in H_{[1, 4m-3]} : X \leftrightarrow H_{2m-1} \}$, $Q^{2m-1}_3 = \{ X \in
  H_{1, 6m-4]} : X \leftrightarrow \Delta(x) \in H_{2m-1} \otimes H_{6m-3}, \forall x \in H \}$. 
  Now an application of Lemma \ref{anti1} shows that the tower $Q^{2m-1}_1 \subset Q^{2m-1}_2 \overset{e^{2m-1}_3}{\subset}Q^{2m-1}_3$ is
   $*$-anti-isomorphic to the tower $A \subset B \overset{e}{\subset} C$ with $e = \tilde{e}^{2m-1}_3$ where
   $A = H_{[-(2m-3), -1]}, B = \{X \in H_{[-(4m-3,) -1]}: X \leftrightarrow H_{-(2m-1)} \}$ and $C = \{ X \in H_{[-(6m-4), -1]} :
  X \leftrightarrow \Delta(x) \in H_{-(6m-3)} \otimes H_{-(2m-1)}, \forall x \in H \}$. Thus it suffices to prove that 
  $BeB = C$, or equivalently, $dim ~BeB = dim ~C$.

 Identify $H_{[-(6m-4), -1]}$ with $P(H^*)_{6m-3}$ and regard $A, B, C$ as subalgebras of $P(H^*)_{6m-3}$. 
In view of Corollary \ref{Eq1}
we see that a general element of $B$ is of the form 
\begin{align*}
Z_U^{P(H^*)}(f^1 \otimes f^2 \otimes \cdots \otimes f^{2m-3} \otimes f^{2m-2}_2 \otimes f^{2m-2}_1 \otimes 
f^{2m-1} \otimes f^{2m} \otimes \cdots \otimes f^{4m-4})
\end{align*}
where $\otimes_{i=1}^{4m-4}f^i \in H^{* \otimes (4m-4)}$ 
and $U$ is the tangle with exactly $4m - 3$ internal $2$-boxes as shown in Figure \ref{fig:pic17}.
\begin{figure}[!h]
\begin{center}
\psfrag{a}{\huge $2m-1$}
\psfrag{b}{\huge $1$}
\psfrag{c}{\huge $2$}
\psfrag{d}{\huge $2m-3$}
\psfrag{e}{\huge $2m-2$}
\psfrag{f}{\huge $2m-1$}
\psfrag{g}{\huge $2m$}
\psfrag{h}{\huge $2m+1$}
\psfrag{i}{\huge $4m-3$}
\psfrag{j}{\huge $2m+2$}
\resizebox{8.0cm}{!}{\includegraphics{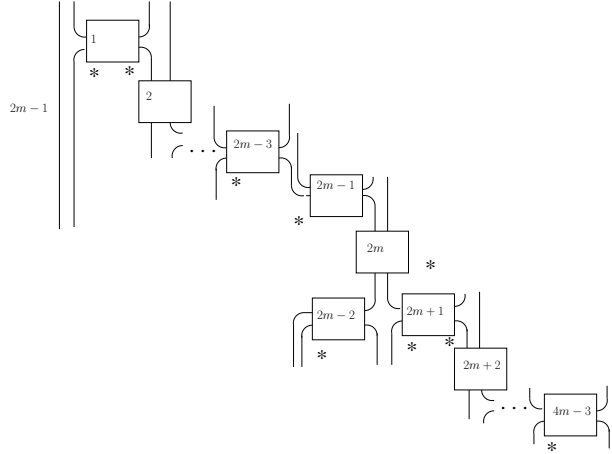}}
\end{center}
\caption{Tangle $U$}
\label{fig:pic17}
\end{figure}\\ 
Now, given 
\begin{align*}
&X = Z_U(f^1 \otimes f^2 \otimes \cdots \otimes f^{2m-3} \otimes f^{2m-2}_2 \otimes f^{2m-2}_1 \otimes 
f^{2m-1} \otimes \cdots \otimes f^{4m-4}), \\
&Y = Z_U(g^1 \otimes g^2 \otimes \cdots \otimes g^{2m-3} \otimes g^{2m-2}_2 \otimes g^{2m-2}_1 
\otimes g^{2m-1} \otimes \underbrace{F(1) \otimes \epsilon \otimes F(1) \otimes \cdots \otimes \epsilon \otimes F(1)}_{\text{$2m-3$ factors}}),
\end{align*}
a little manipulation with the relations (E) and (A) shows that the element $XeY$ equals  
\begin{align*}
 \delta^{-(2m-1)} &Z_Q^{P(H^*)}(\otimes_{i=1}^{2m-3}f^i
\otimes_{i= 2m+1}^{4m-4} f^i \otimes_{i=1}^{2m-3} g^i 
\otimes g^{2m-2}_1Sg^{2m-1} \otimes
g^{2m-2}_2f^{2m}\\ &\otimes g^{2m-2}_3f^{2m-1} \otimes g^{2m-2}_4Sf^{2m-2}_2 \otimes g^{2m-2}_5 \otimes f^{2m-2}_1)
\end{align*}
 where $Q$ 
is the tangle as shown in Figure \ref{fig:pic23}.
\begin{figure}[!h]
\begin{center}
\psfrag{b}{\Huge $1$}
\psfrag{c}{\Huge $2$}
\psfrag{d}{\Huge $2m-3$}
\psfrag{e}{\Huge $6m-6$}
\psfrag{f}{\Huge $6m-4$}
\psfrag{g}{\Huge $6m-7$}
\psfrag{h}{\Huge $6m-8$}
\psfrag{i}{\Huge $4m-7$}
\psfrag{j}{\Huge $2m-2$}
\psfrag{s}{\Huge $4m-6$}
\psfrag{t}{\Huge $4m-5$}
\psfrag{u}{\Huge $6m-10$}
\psfrag{v}{\Huge $6m-5$}
\psfrag{w}{\Huge $6m-3$}
\psfrag{x}{\Huge $6m-9$}
\resizebox{7.5cm}{!}{\includegraphics{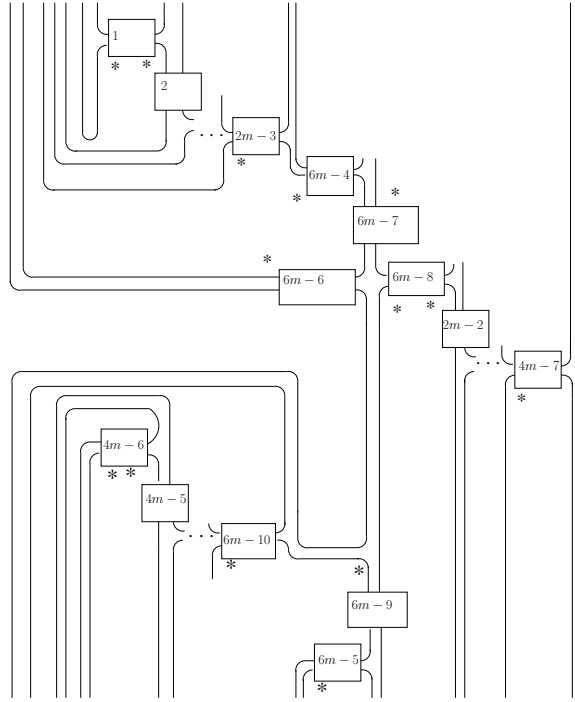}}
\end{center}
\caption{Tangle $Q$}
\label{fig:pic23}
\end{figure}
Observe that $ Q \in {\mathcal{T}}{(6m-3)}$.
Let $\theta : H^{*\otimes 5} \rightarrow H^{* \otimes 6}$ be the linear map defined by 
\begin{align*} 
f \otimes g \otimes k \otimes u \otimes v 
\mapsto f_1 Sg \otimes f_2 k \otimes f_3 u \otimes f_4 Sv_2 \otimes f_5 \otimes v_1.
\end{align*}
We assert that $\theta$ is injective. To see this one
can easily verify that the map $\theta^{\prime} :  H^{*\otimes 6} \rightarrow H^{* \otimes 5}$ given by 
\begin{align*}
f \otimes g \otimes k \otimes p \otimes u \otimes v \mapsto p(1) u_4 \otimes Sf u_3  \otimes Su_2 g \otimes Su_1 k \otimes v
\end{align*}
is a left inverse 
of $\theta$, proving the assertion. Now clearly $BeB$ contains the image of the injective linear map $Z_Q \circ ({Id_{H^*}}^{\otimes (6m-10)} 
\otimes \theta) : H^{*\otimes (6m-5)} \rightarrow P(H^*)_{6m-3}$ and consequently $dim ~BeB \geq (dim ~H)^{6m-5}$. 
To finish the proof we just need to show that $dim ~C \leq (dim ~H)^{6m-5}$, or equivalently, $dim ~Q^{2m-1}_3 \leq (dim ~H)^{6m-5}$. 
\begin{figure}[!h]
\begin{center}
\psfrag{1}{\huge $1$}
\psfrag{2}{\huge $2$}
\psfrag{3}{\huge $3$}
\psfrag{4}{\huge $4$}
\psfrag{m-1}{\huge $m-1$}
\psfrag{m}{\huge $m$}
\psfrag{m+1}{\huge $m+1$}
\psfrag{3m-3}{\huge $3m-3$}
\psfrag{3m-2}{\huge $3m-2$}
\psfrag{3m-1}{\huge $3m-1$}
\psfrag{3m}{\huge $3m$}
\psfrag{4m-5}{\huge $4m-5$}
\psfrag{4m-3}{\huge $4m-3$}
\psfrag{4m-4}{\huge $4m-4$}
\psfrag{6m-6}{\huge $6m-6$}
\psfrag{6m-5}{\huge $6m-5$}
\psfrag{6m-4}{\huge $6m-4$}
\resizebox{20.0cm}{!}{\includegraphics{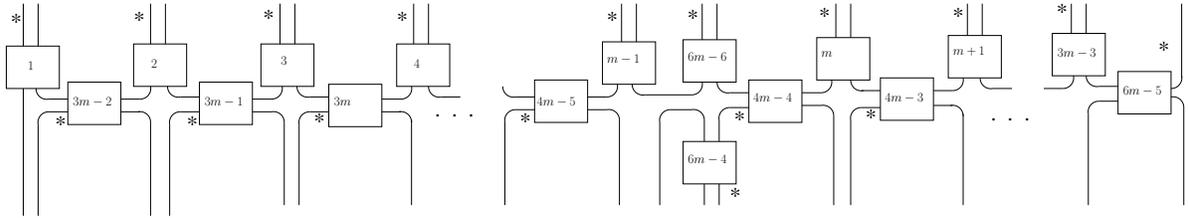}}
\end{center}
\caption{Tangle $R$}
\label{fig:pic63}
\end{figure}
Let us consider the tangle $R \in {\mathcal{T}}{(6m-3)}$ as shown in Figure \ref{fig:pic63}. Since $Z_R^{P(H)}$ is a linear isomorphism of 
$H^{\otimes(6m-4)}$ onto $P(H)_{6m-3}$, we observe, in view of Proposition \ref{descrip}, that the space $Q^{2m-1}_3$ is linearly isomorphic
to the space
\begin{align*}
  \{\otimes_{i = 1}^{6m-4} x^i \in H^{\otimes(6m-4)} : \alpha_h^{2m-1, 3}(Z_R^{P(H)}(\otimes_{i = 1}^{6m-4} x^i)) =
 Z_R^{P(H)}(\otimes_{i = 1}^{6m-4} x^i)\}. 
\end{align*}
A simple computation shows that 
\begin{align*}
 \alpha_h^{2m-1, 3}(Z_R^{P(H)}(\otimes_{i = 1}^{6m-4} x^i)) = Z_R^{P(H)}(\otimes_{i = 1}^{6m-7} x^i \otimes h_1 x^{6m-6} \otimes 
h_2 x^{6m-5} \otimes h_3 x^{6m-4})
\end{align*}
and hence, $Q^{2m-1}_3$ is linearly isomorphic to the space $V$ defined by
\begin{align*}
 V = \{\otimes_{i = 1}^{6m-4} x^i \in H^{\otimes(6m-4)} : h_1x^1 \otimes h_2x^2 \otimes h_3x^3 \otimes_{i=4}^{6m-4} x^i = 
 \otimes_{i = 1}^{6m-4} x^i \}.
\end{align*}
Proceeding in a similar fashion as in the last part of the proof of Proposition \ref{even}, one can easily see that 
$dim ~V \leq (dim ~H)^{6m-5}$, completing the proof.
\end{proof}

We are now ready to conclude Theorem \ref{depth}.
\begin{proof}[Proof of Theorem \ref{depth}]
Follows immediately from Propositions \ref{even} and \ref{odd}.
\end{proof}
 \section{Structure maps on $\mathcal{N}^{\prime} \cap \mathcal{M}_2$}
 The main result of \cite{Sde2018} asserts that the quantum double inclusion of $R^H \subset R$ is isomorphic to 
  $R \subset R \rtimes D(H)^{cop}$ for some outer action of $D(H)^{cop}$ on $R$. We proved this result by constructing a model 
  $\mathcal{N} \subset \mathcal{M}$ (see \cite[Definition 18]{Sde2018} for the definition of 
  $\mathcal{N}$) for the quantum double inclusion of $R^H \subset R$ and then showing that the planar algebras associated to 
  $\mathcal{N} \subset \mathcal{M}$ and $R \subset R \rtimes D(H)^{cop}$ are isomorphic. 
  As an immediate consequence of this result, we obtain that the relative 
 commutant $\mathcal{N}^{\prime} \cap \mathcal{M}_2$ is isomorphic to $D(H)^{cop *} ( = D(H)^{* op})$ as Kac algebras. From the proof of the 
 main result of \cite{Sde2018}, namely \cite[Theorem 40]{Sde2018}, the structure maps on $\mathcal{N}^{\prime} \cap \mathcal{M}_2$ can not 
 directly be derived. In this section 
 we explicitly describe the structure maps of $\mathcal{N}^{\prime} \cap \mathcal{M}_2$ which will be useful in $\S 7$ to achieve a simple 
 and nice description of the weak Hopf $C^*$-algebra structures on $(\mathcal{N}^m)^{\prime} \cap \mathcal{M}_2$ ($m > 2$).

%
 
\begin{figure}[!h]
\begin{center}
\psfrag{1}{\huge $1$}
\psfrag{2}{\huge $2$}
\psfrag{3}{\huge $3$}
\resizebox{10.0cm}{!}{\includegraphics{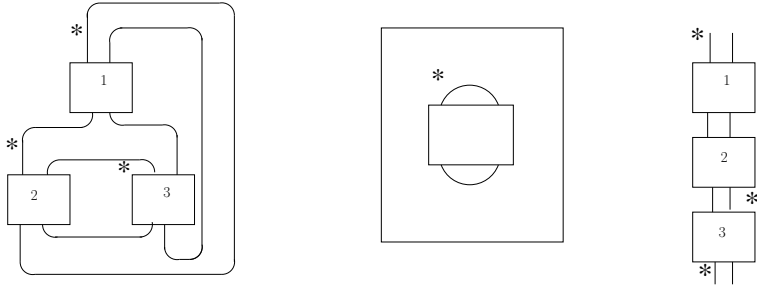}}
\end{center}
\caption{Tangles $E$(left), $F$(middle) and $G$(right)}
\label{fig:pic47}
\end{figure}

Let $N \subset M$ be a finite-index, depth two, irreducible subfactor and
let $N ( = M_0 ) \subset M ( = M_1 ) \subset M_2 \subset 
\cdots$ be the associated tower of basic construction.
Then the relative commutants $N^{\prime} \cap M_2$ and $M^{\prime} \cap M_3$ admit
mutually dual Kac algebra structures. 
Let $P$ denote the subfactor planar algebra associated to $N \subset M$ so that $P_2 = N^{\prime} \cap M_2$.
The next Theorem \ref{Kac} summarises the content of \cite[$\S 2$]{DasKdy2005} where the authors gave pictorial description of
the structure maps on $P_2$. 
 Before we state the theorem, we need to specify certain useful tangles. Let $E, F, G$ denote tangles as shown
in Figure \ref{fig:pic47}.

\begin{theorem}\label{Kac}\cite{DasKdy2005}
The counit $\varepsilon : P_2 \rightarrow \mathbb{C}$ and the antipode $S : P_2 \rightarrow P_2$ are defined for $a \in P_2$ by
\begin{align*}
\varepsilon(a) = [M : N]^{- \frac{1}{2}} Z^P_F(a), \  S(a) = Z^P_{R_2}(a),
\end{align*} 
and the comultiplication $\Delta : P_2 \rightarrow P_2 \otimes P_2$ is the unique linear map such that the equation 
\begin{align*}
 Z_E^P(a \otimes x \otimes y) = [M : N]^{- \frac{1}{2}} tr_2^{0, +}(a_1x) \ tr_2^{0, +}(a_2y)
\end{align*}
holds for all $a, x, y \in P_2$.  
\end{theorem}
%
%

Recall from \cite[Theorem 38]{Sde2018} that the planar algebra associated to $\mathcal{N} \subset \mathcal{M}$, denoted $^*\!Q$, is a planar
subalgebra of $^{*(2)}\!P(H^*)$ and for each integer $k \geq 1$, the space $^*\!Q_k$ is the opposite algebra of 
\begin{align*}
 \{ &X \in H_{[0, 2k-2]}: X \ \mbox{commutes with} \ \Delta_{l-1}(x) \in \otimes_{i = 1}^l H_{4i-3}, \forall x \in H \ \mbox{where} \ 
 l = \frac{k}{2} \ \mbox{or} \ \frac{k+1}{2} \\
 &\mbox{according as} \ k  \ \mbox{is even or odd} \}.  
\end{align*}
Thus, in particular, $^*\!Q_2$ is the opposite algebra of 
 \begin{align*}
  \{X \in H^* \rtimes H \rtimes H^* : X \ \mbox{commutes with the middle} \ H \}.  
 \end{align*}
 That is, $^*\!Q_2$ is the opposite algebra of $\mathcal{S}$ (see Lemma \ref{element} for the definition of $\mathcal{S}$)
 and consequently, by Lemma \ref{element}, 
 $^*\!Q_2$ precisely consists of elements of the form $Z_A^{P(H^*)}(f_2 \otimes f_1 \otimes g)$
where $f \otimes g \in H^{*\otimes2}$. We apply Theorem \ref{Kac} above to derive 
the structure maps for $^*\!Q_2$.
\begin{proposition}\label{structure}
 Let $X = Z_A^{P(H^*)}(f_2 \otimes f_1 \otimes k) \in {^*\!Q_2}$, then 
 \begin{align*}
&\mbox{Comultiplication:} \ \Delta(X) = \delta Z_A((f S \phi_2)_2 \otimes (f S \phi_2)_1 
 \otimes \phi_1 k_2 S \phi_3) \otimes Z_A( (\phi_4)_2 \otimes (\phi_4)_1 \otimes k_1),\\
 & \mbox{Antpode:} \ S(X) = Z_A(Sf_2 \otimes Sf_3 \otimes f_1SkSf_4), \\
 &\mbox{Counit:} \ \varepsilon(X) = \delta f(h) k(1), \\
 &\mbox{Involution:} \ X^* = Z_A^{P(H^*)}(Sf_1^* \otimes Sf_2^* \otimes g^*).
 \end{align*}
\end{proposition}
\begin{proof}
 Applying Theorem \ref{Kac}, the formula for $S(X)$ follows directly by using the relations (E) and (A) whereas to verify the formula for $\varepsilon(X)$
 one needs to use the relations (T) and (M). To verify the involution formula we just need to observe that $X^* = Z_{A^*}^{P(H^*)}(f_2^* \otimes f_1^*
 \otimes g^*) = Z_A^{P(H^*)}(Sf_1^* \otimes Sf_2^* \otimes g^*)$.
 
 We now verify the formula for $\Delta(X)$. It follows from Theorem \ref{Kac} that given any $W \in {^*\!Q_2}$, then $\Delta_{^*\!Q_2}(W)
  = W_1 \otimes W_2$ is that element of $(^*\!Q_2)^{\otimes 2}$ such the equation
 \begin{align}\label{cc}
  Z_E^{^*\!Q}(W \otimes Y \otimes Z) = \delta^{-2}(= [\mathcal{M} : \mathcal{N}]^{- \frac{1}{2}}) tr_2^{(0, +)}(W_1 \ Y) \ tr_2^{(0, +)}(W_2 \ Z)
 \end{align}
 holds for all $Y, Z \in {^*\!Q_2}$.
  Let $Y, Z$ be arbitrary elements in $^*\!Q_2$, say, $Y = Z_A(g_2 \otimes g_1 \otimes p), Z = Z_A(u_2 \otimes u_1 \otimes v)$. Then 
 the element $Z^{^*\!Q}_E(X \otimes Y \otimes Z) = Z_{E^{*(2)}}^{P(H^*)}(X \otimes Y \otimes Z)$ is as shown in Figure \ref{fig:pic48}.
 A very lengthy but completely routine computation in $P(H^*)$ along with repeated application of the well-known Hopf-algebraic formulae such as
$h_1 a \otimes h_2 = h_1 \otimes h_2 Sa$ shows that
\begin{align*}
 &Z^{^*\!Q}_E(X \otimes Y \otimes Z) = \delta^5 f(Sh^1_1 Sh^2_1) g(h^2_2 h^1_2) k(h^3_1 h^4_1) u(h^4_4 h^2_3 h^1_3 Sh^4_2) p(h^4_3) v(h^3_2)\\ 
 &= \delta^{-2} \delta tr_2((Z_A((f S \phi_2)_2 \otimes (f S \phi_2)_1 
 \otimes \phi_1 k_2 S \phi_3) Y) tr_2((Z_A( (\phi_4)_2 \otimes (\phi_4)_1 \otimes k_1)Z),
\end{align*}
verifying the formula for $\Delta(X)$.
 \begin{figure}[!h]
\begin{center}
\psfrag{1}{\huge $g_2$}
\psfrag{2}{\huge $g_1$}
\psfrag{3}{\huge $p$}
\psfrag{4}{\huge $u_2$}
\psfrag{5}{\huge $u_1$}
\psfrag{6}{\huge $v$}
\psfrag{7}{\huge $f_2$}
\psfrag{8}{\huge $f_1$}
\psfrag{9}{\huge $k$}
\resizebox{6.0cm}{!}{\includegraphics{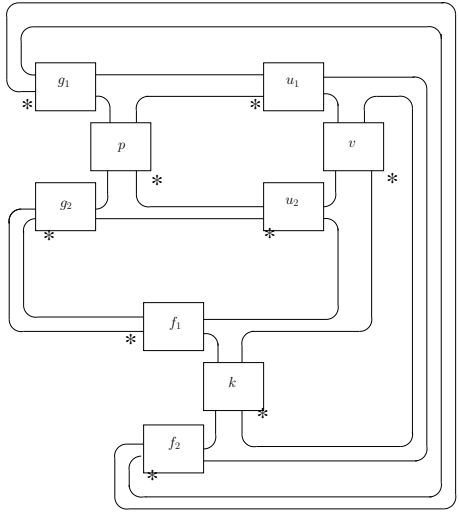}}
\end{center}
\caption{ $Z^{^*\!Q}_E(X \otimes Y \otimes Z)$}
\label{fig:pic48}
\end{figure}
\end{proof}
Using Lemma \ref{eq}, it follows immediately from Proposition \ref{structure} that the structure maps of $^*\!Q_2$ can also be expressed as:  
\begin{lemma}\label{sstructure}
 Given $X = f_1Sk_3Sf_3 \rtimes F^{-1}(f_2k_2) \rtimes Sk_1 \in {^*\!Q_2}$ where $f \otimes k \in {H^*}^{\otimes 2}$, then
 \begin{align*}
 \mbox{Comultiplication:} \ \Delta(X) = \ &\delta \ ((fS\phi_2)_1 \ S(\phi_1 k_2 S\phi_3)_3 \ S(fS\phi_2)_3 \rtimes F^{-1}((fS\phi_2)_2 \ 
 (\phi_1 k_2 S\phi_3)_2) \rtimes S(\phi_1 k_2 S\phi_3)_1)\\
 &\otimes \ ((\phi_4)_1 \ S(k_1)_3 \ S(\phi_4)_3 \rtimes F^{-1}((\phi_4)_2 (k_1)_2) \rtimes S(k_1)_1),\\
 \mbox{Antipode:} \ S(X) = \ & k_1 \rtimes F^{-1}S(f_2k_2) \rtimes S(f_1Sk_3Sf_3),\\
  \mbox{Counit:} \ \epsilon(X) = \ & \delta f(h) k(1),\\
  \mbox{Involution:} \ X^* = \ & (Sf^*)_1 S((g^*)_3) (Sf^*)_3 \rtimes F^{-1}((Sf^*)_2 (g^*)_2) \rtimes S((g^*)_1).
 \end{align*} 
\end{lemma}

\begin{remark}\label{Dr}
In $\S 1.4$, we considered a version of $D(H)^*$ whose underlying vector space is $H^* \otimes H^*$ and the structure maps are given by Lemma 
\ref{dr}. 
Consider the linear isomorphism
\begin{align*}
 \nu: D(H)^{* op} = H^* \otimes H^* \mapsto {^*\!Q_2} 
\end{align*}
that takes $g \otimes f \mapsto Z_A^{P(H^*)}(f_2 \otimes f_1 \otimes g)$. It follows immediately from Proposition \ref{structure} and the  
structure maps on $D(H)^*$ as given by Lemma \ref{dr}
that $\nu$ is an isomorphism of Kac algebras.

\end{remark}

 \section{Weak Hopf $C^*$-algebra structure on $(\mathcal{N}^m)^{\prime} \cap \mathcal{M}_2, m > 2$}
 It is well-known (see \cite{Das2004}, \cite{NikVnr2000}) that if $N \subset M$ is a depth two, reducible, finite-index inclusion of $II_1$-factors
 and if $N (= M_0) \subset M (= M_1) \subset M_2 \subset M_3 \subset \cdots$ is the Jones' basic construction tower associated to $N \subset M$,
  then the relative commutants
  $N^{\prime} \cap M_2$ and $M^{\prime} \cap M_3$ admit mutually dual weak Hopf $C^*$-algebra structures.
  The following Theorem \ref{weak} (reformulation of Proposition 4.7 of \cite{Das2004}) explicitly describes the weak Hopf $C^*$-algebra 
  structures on 
  $N^{\prime} \cap M_2$.
  Before we state the theorem, we need to fix some notations. Let $P$ be the planar algebra associated to $N \subset M$ so that 
  $P_2 = N^{\prime} \cap M_2$, $P_{1, 2} =  M^{\prime} \cap 
  M_2$. Set $[M:N] = d^2$. Further, 
  let $z_R$ denote the unique element of $P_{1, 2}$ for which $tr_L(x) = tr_2(z_R x)$
  for all $x \in P_{1, 2}$ where $tr_L$ is the trace of left regular representation of $P_{1, 2}$ and $tr_2$ denotes the normalised pictorial 
  trace on $P_2$. One can easily see that $z_R$ is a well-defined, central, positive, invertible element of $P_{1, 2}$. 
  By $\omega_R$ we denote the unique positive
  square root of $z_R$ and let $\omega_L$ be $Z_{R_2^2}(\omega_R)$. We will use $\omega$ to denote $\omega_L \omega_{R}^{-1}$. The following theorem 
  describes the structure maps of $P_2$.
  \begin{theorem}\label{weak} \cite{Das2004}
   The comultiplication $\Delta : P_2 \rightarrow P_2 \otimes P_2$ is the unique linear map such that the equation
   \begin{align*}
    Z_E^P(\omega_R a \omega_L \otimes x \otimes y) = d^{-1} \ Z_{tr_2^{0, +}}^P(\omega_R a_1 \omega_L x)\  .\ Z_{tr_2^{0, +}}^P(\omega_R a_2 \omega_L y) 
   \end{align*}
   holds in $P$ for all $a, x ,y \in P_2$.
   The counit $\varepsilon:P_2 \rightarrow \mathbb{C}$ and the antipode $S: P_2 \rightarrow P_2$ are defined by
   \begin{align*}
    \varepsilon(a) = d^{-1} Z_F^P(w_R a w_L),\ \mbox{and}\ \ S(a) = Z_G^P(w \otimes a \otimes w^{-1}),
   \end{align*}
   for all $a$ in $P_2$ where $E, F, G$ are the tangles as shown in the Figure \ref{fig:pic47}.
 \begin{figure}[!h]\label{comult}
 \begin{center}
 \psfrag{1}{\Huge $\omega_R a \omega_L$}
 \psfrag{2}{\Huge $x$}
 \psfrag{3}{\Huge $y$}
 \psfrag{4}{\Huge $\omega_R a_1 \omega_L$}
 \psfrag{5}{\Huge $x$}
 \psfrag{6}{\Huge $\omega_R a_2 \omega_L$}
 \psfrag{7}{\Huge $y$}
 \psfrag{8}{\Huge $d^{-1}$}
 \psfrag{9}{\Huge $d^{-1}$}
 \psfrag{10}{\Huge $d^{-1}$}
 \psfrag{11}{\Huge $=$}
 \resizebox{9.0cm}{!}{\includegraphics{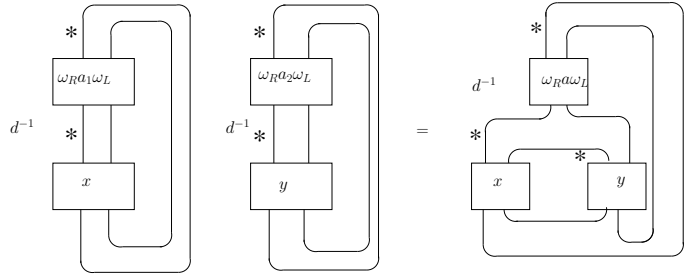}}
 \end{center}
 \caption{Definition of $\Delta$}
 \label{fig:pic52}
 \end{figure}
 \begin{figure}[!h]\label{antico}
\begin{center}
\psfrag{a}{\Huge $\omega$}
\psfrag{b}{\Huge $a$}
\psfrag{c}{\Huge $\omega^{-1}$}
\psfrag{s}{\Huge $S(a) = $}
\psfrag{e}{\Huge $\epsilon(a) = $}
\psfrag{d}{\Huge $d^{-1}$}
\psfrag{w}{\Huge $\omega_R a \omega_L$}
\psfrag{,}{\Huge $,$}
\resizebox{10.0cm}{!}{\includegraphics{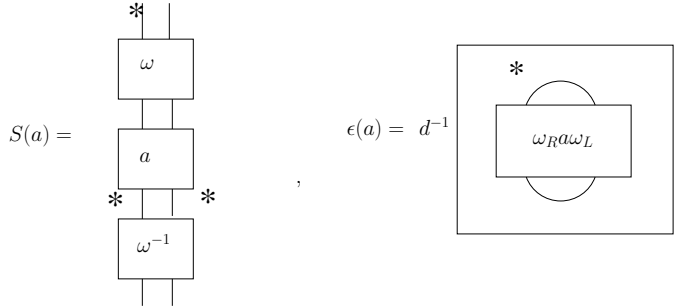}}
\end{center}
\caption{Definitions of $S$ and $\epsilon$}
\label{fig:pic101}
\end{figure}

  \end{theorem}
  We use Theorem \ref{weak} to recover the weak Hopf $C^*$-algebra structure on $^*\!Q^m_2$ for all $m > 2$. Let us
 use the symbols $\omega_R^m, \omega_L^m$ and $\omega^m$ to denote the elements $\omega_R, \omega_L$ and $\omega$ respectively 
 of $^*\!Q^m_2$.

  Note that in order to find the structure maps of $^*\!Q^m_2$ using Theorem \ref{weak}, we must know the elements 
  $\omega_R^m, \omega_L^m$ and $\omega^m$. It follows from Lemma \ref{second} that 
  \begin{align*}
 ^*\!Q^m_{1, 2} = 
 \begin{cases}
  \{X \in H^{op}_{[m+1, 2m-1]}: X \ \mbox{commutes with} \ H_m \}, \ \mbox{if} \ m \ \mbox{is odd}\\
  \{X \in H^{op}_{[m, 2m-2]}: X \ \mbox{commutes with} \ H_{m-1} \}, \ \mbox{if} \ m \ \mbox{is even}.
 \end{cases}
  \end{align*}
  Then by an appeal to Lemma \ref{commutants} it follows immediately that $^*\!Q^m_{1, 2}$ is identified with the subalgebra 
  \begin{align*}
   H^{op}_{[m+2, 2m-1]} \ \mbox{or} \ H^{op}_{[m+1, 2m-2]} \ \mbox{of} \ ^*\!Q^m_2     
  \end{align*}
  according as $m$ is odd or even.
  Thus, for any $m > 2$, $^*\!Q^m_{1, 2} \cong A(H)_{m-2}^{op}$ 
  as algebras. 
  One can easily see that 
  if $A$ is any multi-matrix algebra over the complex field, then for any $a \in A, tr_L^A(a) = tr_R^A(a)$ and consequently,
  $tr_L^A(a) = tr_L^{A^{op}}(a)$ where $tr_L^A(a)$ (resp., $tr_R^A(a)$) denotes the trace of the linear endomorphism of $A$ given by left (resp., right)
  multiplication by $a$. Thus, given any $X$ in $^*\!Q^m_{1, 2} = A(H)_{m-2}^{op}$, we have $tr_L^{^*\!Q^m_{1, 2}}(X) = tr_L^{A(H)_{m-2}^{op}}(X)
  = tr_L^{A(H)_{m-2}}(X)$. The following lemma computes $tr_L^{A(H)_k}(X)$ for any $X \in A(H)_k$ where $k \geq 1$ is an integer.
  \begin{lemma}\label{4}
   $tr_L^{A(H)_k}(X) = (dim ~H)^k$ (normalised pictorial trace of $X$).
  \end{lemma}
  \begin{proof}
   If $k$ is even, then $A(H)_k$ is a matrix algebra by Lemma \ref{matrixalg} and
   the result follows immediately.
   Now suppose that $H$ is a finite-dimensional Hopf algebra acting on a finite-dimensional algebra $A$. A simple exercise in linear algebra
   shows that for any $a \in A$, $tr^{A \rtimes H}_L(a \rtimes 1) = dim ~H \ tr^A_L(a)$. Thus if $k$ is odd, then given $X \in A(H)_k$, 
   we note that $tr_L^{A(H)_k}(X) = \frac{1}{dim ~H} tr_L^{A(H)_{k+1}}(X \rtimes 1) = \frac{1}{dim ~H} dim ~H^{k+1}$ (normalised pictorial
   trace of $X \rtimes 1) = (dim ~H)^k$ (normalised pictorial trace of $X$) where the second equality follows since $k+1$ is even, completing the proof.
  \end{proof}
 
  As an immediate corollary we have
  \begin{corollary}\label{omega}
  $\omega^m_L = \omega^m_R = (dim ~H)^{\frac{m-2}{2}} 1_{^*\!Q^m_2}$ and hence, $\omega^m = 1_{^*\!Q^m_2}$ for $m > 2$.  
  \end{corollary}
\begin{proof}
 It follows from Lemma \ref{4} and the discussion preceding Lemma \ref{4} that for any $X$ in $^*\!Q^m_{1, 2}$, $tr_L^{^*\!Q^m_{1, 2}}(X) =
 tr_L^{A_{m-2}}(X) = (dim ~H)^{m-2}$ normalised pictorial trace of $X$. Hence, $\omega_R^m = (dim ~H)^{\frac{m-2}{2}} 1_{^*\!Q^m_{1, 2}}$ and consequently, 
 $\omega_L^m = (dim ~H)^{\frac{m-2}{2}} 1_{^*\!Q^m_2}$, $\omega^m = 1_{^*\!Q^m_2}$. 
 \end{proof}
 We now proceed towards recovering the structure maps of $^*\!Q^m_2$. The entire procedure solely relies on Hopf-algebraic as well as pictorial computations in 
 $P(H)$ or $P(H^*)$. At this point, it is worth recalling from Corollary \ref{Eq1} and Remark \ref{imp} that a general element of $^*\!Q^m_2$ is of the form
 \begin{align}\label{1}
  x^1 \rtimes f^2 \rtimes \cdots \rtimes x^{m-2} \rtimes k_1 Sp_3 Sk_3 \rtimes
  F^{-1}(k_2 p_2) \rtimes Sp_1 \rtimes x^{m+2} \rtimes f^{m+3} \rtimes \cdots \rtimes x^{2m-1} 
 \end{align}
or 
\begin{align}\label{2}
 f^1 \rtimes x^2 \rtimes \cdots \rtimes x^{m-2} \rtimes k_1 Sp_3 Sk_3 \rtimes
  F^{-1}(k_2 p_2) \rtimes Sp_1 \rtimes x^{m+2} \rtimes \cdots \rtimes f^{2m-1} 
\end{align}
according as $m$ is odd or even. Moreover, if $m$ is even, the elements of $^*\!Q^m_2$ can also be expressed as 
\begin{align}\label{3}
 Z_{A(m-2, m-2)}^{P(H^*)}(f^1 \otimes f^2 \otimes \cdots \otimes f^{m-2} \otimes f^{m-1}_2 \otimes 
  f^{m-1}_1 \otimes f^m \otimes f^{m+1} \otimes \cdots \otimes f^{2m-2}).
\end{align}

First we find the formula for antipode in $^*\!Q^m_2$. It follows from Theorem \ref{weak} and Lemma \ref{omega} that  
for any $X \in {^*\!Q^m_2}$, 
\begin{align*}
 S(X) &= Z^{^*\!Q^m}_{R^2_2} (X)\\
 &= Z_{(R_2^2)^{* m}}^{P(H^*)}(X)\\
  &= Z_{(R_2^2)^m}^{P(H^*)}(X) \ (\mbox{since} \ (R_2^2)^* = R_2^2). 
\end{align*}
 Let us consider a general element, say $X$, of $^*\!Q^m_2$ as given by \eqref{1} or 
\eqref{2} according as $m$ is odd or even. Assume that $m$ is even. Then note that $X$ is identified
 with 
 \begin{align*}
 Z_{T^{2m}}^{P(H^*)}(f^1 \otimes Fx^2 \otimes \cdots \otimes Fx^{m-2} \otimes k_1 Sp_3 Sk_3 \otimes k_2 p_2 \otimes Sp_1 
 \otimes Fx^{m+2} \otimes f^{m+3} \otimes
  \cdots \otimes f^{2m-1}). 
  \end{align*}
  Consequently, 
  \begin{align*}
   S(X) = Z_{(R_2^2)^m}^{P(H^*)}(Z_{T^{2m}}^{P(H^*)}(&f^1 \otimes Fx^2 \otimes \cdots \otimes Fx^{m-2}
   \otimes k_1 Sp_3 Sk_3 \otimes k_2 p_2 \otimes Sp_1 \\
 &\otimes Fx^{m+2} \otimes f^{m+3} \otimes \cdots \otimes f^{2m-1}))  
  \end{align*}
  which, by repeated application of 
  the relation (A) in $P(H^*)$ is easily seen to be equal to
  \begin{align*}
  Z_{T^{2m}}^{P(H^*)}(Sf^{2m-1} \otimes \cdots \otimes SFx^{m+2} \otimes p_1 \otimes S(k_2 p_2) \otimes S(k_1 Sp_3 Sk_3) \otimes 
  SFx^{m-2} \otimes 
  \cdots \otimes SFx^2 \otimes Sf^1)
 \end{align*}
 which, by virtue of the fact that $FS = SF$, is identified with 
 \begin{align*}
  Sf^{2m-1} \rtimes \cdots \rtimes Sx^{m+2} \rtimes p_1 \rtimes F^{-1}S(k_2 p_2) \rtimes S(k_1 Sp_3 Sk_3) \rtimes 
  Sx^{m-2} \rtimes 
  \cdots \rtimes Sx^2 \rtimes Sf^1.
 \end{align*}
 Thus, we obtain the formula for $S(X)$ when $m$ is even. Proceeding exactly the same way, one can show that the formula for $S(X)$, when $m$ is odd, 
 is given by 
 \begin{align*}
   Sx^{2m-1} \rtimes \cdots \rtimes Sx^{m+2} \rtimes p_1 \rtimes F^{-1}S(k_2 p_2) \rtimes S(k_1 Sp_3 Sk_3) \rtimes  
  Sx^{m-2} \rtimes 
  \cdots \rtimes Sf^2 \rtimes Sx^1.
 \end{align*}
 Thus we have proved that:
 \begin{lemma}\label{wanti}
  Let $X \in {^*\!Q^m_2}$ be as given by \eqref{1} or \eqref{2} according as $m$ is odd or even. Then $S(X)$ is given by 
  \begin{align*}
 Sf^{2m-1} \rtimes \cdots \rtimes Sx^{m+2} \rtimes p_1 \rtimes F^{-1}S(k_2 p_2) \rtimes S(k_1 Sp_3 Sk_3) \rtimes 
  Sx^{m-2} \rtimes 
  \cdots \rtimes Sx^2 \rtimes Sf^1  
  \end{align*}
  or
  \begin{align*}
  Sx^{2m-1} \rtimes \cdots \rtimes Sx^{m+2} \rtimes p_1 \rtimes F^{-1}S(k_2 p_2) \rtimes S(k_1 Sp_3 Sk_3) \rtimes  
  Sx^{m-2} \rtimes 
  \cdots \rtimes Sf^2 \rtimes Sx^1
 \end{align*}
 according as $m$ is odd or even.
 \end{lemma}

 Among all the structure maps the hardest is to 
 recover the comultiplication formula and we do it in steps.
 By an appeal to Corollary \ref{omega} and Lemma \ref{index}, it follows immediately from Theorem \ref{weak} that given $X \in ^*\!Q^m_2, 
 \Delta_{^*\!Q^m_2}(X) = 
 X_1 \otimes X_2$ is that element of $(^*\!Q^m_2)^{\otimes 2}$ such that the equation  
 \begin{align}\label{c}
  Z_E^{^*\!Q^m}(X \otimes Y \otimes Z) = \delta^{-m} \ \delta^{2m-4} \ Z_{tr_2^{(0, +)}}^{^*\!Q^m}( X_1 Y)\  .\ 
  Z_{tr_2^{(0, +)}}^{^*\!Q^m}(X_2 Z), 
 \end{align}
 holds for all $Y, Z \in ^*\!Q^m_2$.

 We begin with finding the comultiplication formula for a special class of elements of
 $^*\!Q^m_2$. 
 Recall from the discussion preceding Proposition \ref{structure} in $\S 6$ that the space $^*\!Q_2$, where $^*\!Q$ is the planar algebra associated
 to $\mathcal{N} \subset \mathcal{M}$, is same as $\mathcal{S}$ as vector spaces but as an algebra it is the opposite algebra of $\mathcal{S}$.
 Given $f \rtimes x \rtimes g \in {^*\!Q_2}$, we define $X^m_{f \rtimes x \rtimes g}$ to be the element of $^*\!Q^m_2$ given by
 $\underbrace{1 \rtimes \epsilon \rtimes \cdots \rtimes 1}_{\text{$m-2$ \mbox{terms}}} \rtimes 
 f \rtimes x \rtimes g \rtimes  \underbrace{1 \rtimes \epsilon \rtimes \cdots \rtimes 1}_{\text{$m-2$ \mbox{terms}}}$ or  
  $\underbrace{\epsilon \rtimes 1 \rtimes \cdots \rtimes 1}_{\text{$m-2$ \mbox{terms}}} \rtimes f \rtimes x \rtimes g \rtimes
  \underbrace{1 \rtimes \epsilon \rtimes \cdots \rtimes \epsilon}_{\text{$m-2$ \mbox{terms}}}$ according as $m$ is odd or even. 
  Let $\Delta_{^*\!Q_2}(f \rtimes x \rtimes g) = (f \rtimes x \rtimes
 g)_1 \otimes (f \rtimes x \rtimes g)_2$. The following lemma computes $\Delta_{^*\!Q^m_2}(X^m_{f \rtimes x \rtimes g})$.
 \begin{lemma}\label{antipode}
  $\Delta_{^*\!Q^m_2}(X^m_{f \rtimes x \rtimes g}) = 1_1 X^m_{(f \rtimes x \rtimes g)_1} \otimes 1_2 X^m_{(f \rtimes x \rtimes g)_2} \in 
  (^*\!Q^m_2)^{\otimes 2}$ where 
  $1_1 \otimes 1_2 = \Delta_{^*\!Q^m_2}(1)$. 
 \end{lemma}
 \begin{proof}
  To avoid notational clumsiness and to elucidate the computational procedure, instead of treating the general case, we explicitly work out the 
  particular case when $m = 4$. The general case when $m > 2$ is even follows in a similar fashion and the case when $m > 2$ is odd 
  follows almost in a similar way with slight modifications.  
  
  Let $f \rtimes x \rtimes g \in {^*\!Q_2}$. Let us consider the element $X^4_{f \rtimes x \rtimes g} \in {^*\!Q^4_2}$. 
  Let $Y, Z$ be arbitrary elements of $^*\!Q^4_2$, say, 
 $Y = Z_{A(2, 2)}^{P(H^*)}(k^1 \otimes k^2 \otimes u_2 \otimes u_1 \otimes v \otimes k^3 \otimes k^4)$ and 
 $Z = Z_{A(2, 2)}^{P(H^*)}(p^1 \otimes p^2 \otimes \tilde{u}_2 \otimes \tilde{u}_1 \otimes \tilde{v} \otimes p^3 \otimes p^4)$. 
It follows from \eqref{c} that in order to verify the comultiplication formula 
for $X^4_{f \rtimes x \rtimes g}$, we just need to verify that 
\begin{align*}
Z^{^*\!Q^4}_{tr_2^{(0, +)}}(1_1 X^4_{(f \rtimes x \rtimes g)_1}Y) \  
Z^{^*\!Q^4}_{tr_2^{(0, +)}}(1_2 X^4_{(f \rtimes x \rtimes g)_2}Z) = Z^{^*\!Q^4}_{E}(X^4_{f \rtimes x \rtimes g} \otimes Y \otimes Z). 
\end{align*}
Another appeal to \eqref{c} shows that 
\begin{align*}
 & Z^{^*\!Q^4}_{tr_2^{(0, +)}}(1_1 X^4_{(f \rtimes x \rtimes g)_1}Y) \  
Z^{^*\!Q^4}_{tr_2^{(0, +)}}(1_2 X^4_{(f \rtimes x \rtimes g)_2}Z)  \\
&= Z^{^*\!Q^4}_{E}(1 \otimes X^4_{(f \rtimes x \rtimes g)_1}Y \otimes X^4_{(f \rtimes x \rtimes g)_2}Z)\\
&= Z^{P(H^*)}_{E^{* (4)}}(1 \otimes Y X^4_{(f \rtimes x \rtimes g)_1} \otimes Z X^4_{(f \rtimes x \rtimes g)_2}).
\end{align*}
where the last equality follows from the definition of adjoint and cabling of a planar algebra.
Now a pleasant but lengthy computation in $P(H^*)$ involving sphericality of $P(H^*)$ and repeated application of the 
relations (E), (T), (C) and (A), shows that 
\begin{align*}
&Z^{P(H^*)}_{E^{* (4)}}(1 \otimes Y X^4_{(f \rtimes x \rtimes g)_1} \otimes Z X^4_{(f \rtimes x \rtimes g)_2})\\
 & = \delta^4 \ p^4(h^1)\ p^3(1)\ k^1(h^2)\ k^2(1)\ (Sp^1_1 k^4)(h^3)\ (v S (k^3 p^1_2 p^2)_2(h^4) \\
& Z_{E^{* (2)}}^{P(H^*)} (1 \otimes Z_A^{P(H^*)}(u_2 \otimes u_1
 \otimes (k^3p^1_2 p^2)_1) \ (f \rtimes x \rtimes g)_1 \otimes Z_A^{P(H^*)}(\tilde{u}_2 \otimes\tilde{u}_1 \otimes
 \tilde{v}) \ (f \rtimes x \rtimes g)_2)\\
 & = \delta^4 \ p^4(h^1)\ p^3(1)\ k^1(h^2)\ k^2(1)\ (Sp^1_1 k^4)(h^3)\ (v S (k^3 p^1_2 p^2)_2(h^4) \\
 & Z_E^{^*\!Q} (1 \otimes (f \rtimes x \rtimes g)_1 \ Z_A^{P(H^*)}(u_2 \otimes u_1
 \otimes (k^3p^1_2 p^2)_1) \otimes (f \rtimes x \rtimes g)_2 \ Z_A^{P(H^*)}(\tilde{u}_2 \otimes\tilde{u}_1 \otimes
 \tilde{v})).
\end{align*}
Now repeated application of Equation \eqref{cc} shows that
\begin{align*}
&\delta^4 \ p^4(h^1)\ p^3(1)\ k^1(h^2)\ k^2(1)\ (Sp^1_1 k^4)(h^3)\ (v S (k^3 p^1_2 p^2)_2(h^4) \\
 & Z_E^{^*\!Q} (1 \otimes (f \rtimes x \rtimes g)_1 \ Z_A^{P(H^*)}(u_2 \otimes u_1
 \otimes (k^3p^1_2 p^2)_1) \otimes (f \rtimes x \rtimes g)_2 \ Z_A^{P(H^*)}(\tilde{u}_2 \otimes\tilde{u}_1 \otimes
 \tilde{v}))\\
 &= \delta^{-2} \ \delta^4 \ p^4(h^1)\ p^3(1)\ k^1(h^2)\ k^2(1)\ (Sp^1_1 k^4)(h^3)\ 
 (v S (k^3 p^1_2 p^2)_2(h^4) \\ 
 &Z_{tr_2^{(0, +)}}^{^*\!Q} ((f \rtimes x \rtimes g)_1 \ Z_A^{P(H^*)}(u_2 \otimes u_1
 \otimes (k^3p^1_2 p^2)_1)) \ Z_{tr_2^{(0, +)}}^{^*Q\!} ((f \rtimes x \rtimes g)_2 \ Z_A^{P(H^*)}(\tilde{u}_2 \otimes\tilde{u}_1 \otimes
 \tilde{v})) \\
 &= \delta^4 \ p^4(h^1)\ p^3(1)\ k^1(h^2)\ k^2(1)\ (Sp^1_1 k^4)(h^3)\ (v S (k^3 p^1_2 p^2)_2(h^4) \\
 & Z_E^{^*\!Q}((f \rtimes x \rtimes g) \otimes Z_A^{P(H^*)}(u_2 \otimes u_1
 \otimes (k^3p^1_2 p^2)_1)) \otimes Z_A^{P(H^*)}(\tilde{u}_2 \otimes\tilde{u}_1 \otimes
 \tilde{v})). 
\end{align*}
Finally, a routine computation in $P(H^*)$ shows that   
\begin{align*}
&Z^{^*\!Q^4}_{E}(X^4_{f \rtimes x \rtimes g} \otimes Y \otimes Z)\\
&= \delta^4 \ p^4(h^1)\ p^3(1)\ k^1(h^2)\ k^2(1)\ (Sp^1_1 k^4)(h^3)\ (v S (k^3 p^1_2 p^2)_2(h^4) \\
 & Z_E^{^*\!Q}((f \rtimes x \rtimes g) \otimes Z_A^{P(H^*)}(u_2 \otimes u_1
 \otimes (k^3p^1_2 p^2)_1)) \otimes Z_A^{P(H^*)}(\tilde{u}_2 \otimes\tilde{u}_1 \otimes
 \tilde{v})). 
\end{align*}
Hence, the formula for $\Delta_{{^*\!Q^4_2}}(X^4_{f \rtimes x \rtimes g})$ is verified.  
 \end{proof}
 We now proceed towards establishing the comultiplication formula for a general element of $^*\!Q^m_2$. Let us take a general element $X$ 
 of $^*\!Q^m_2$ as given by \eqref{1} or \eqref{2} according as $m$ is odd or even. 
 The multiplication in $^*\!Q^m_2$ shows that $X$ can be expressed as 
 \begin{align*}
 X = X^m_3 \ X^m_{k_1 Sp_3 Sk_3 \rtimes F^{-1}(k_2 p_2) \rtimes Sk_1} \ X^m_1
 \end{align*}
 with 
  $X^m_1 \otimes X^m_{k_1 Sp_3 Sk_3 \rtimes F^{-1}(k_2 p_2) \rtimes Sk_1} \otimes X^m_3$ ($\in (^*\!Q^m_2)^{\otimes 3}$) being given by
 \begin{align*}
  (f^1 \rtimes \cdots \rtimes x^{m-2} \rtimes
 \underbrace{\epsilon \rtimes 1 \rtimes \cdots \rtimes \epsilon}_{\text{$m+1$ \mbox{terms}}}) \otimes 
 X^m_{k_1 Sp_3 Sk_3 \rtimes F^{-1}(k_2 p_2) \rtimes Sk_1}
 \otimes 
 (\underbrace{\epsilon \rtimes \cdots \rtimes \epsilon}_{\text{$m+1$ \mbox{terms}}} \rtimes 
 x^{m+2} \rtimes \cdots \rtimes f^{2m-1})
 \end{align*}
 or 
 \begin{align*}
 (x^1 \rtimes \cdots \rtimes x^{m-2} \rtimes \underbrace{\epsilon \rtimes 1 \rtimes \cdots \rtimes 1}_{\text{$m+1$ \mbox{terms}}}) \otimes 
 X^m_{k_1 Sp_3 Sk_3 \rtimes F^{-1}(k_2 p_2) \rtimes Sk_1} \otimes 
 (\underbrace{1 \rtimes \cdots \rtimes \epsilon}_{\text{$m+1$ \mbox{terms}}} \rtimes x^{m+2} \rtimes \cdots \rtimes x^{2m-1}) 
 \end{align*}
 according as $m$ is even or odd. It is then not hard to show using Equation \eqref{c} and Lemma \ref{antipode} that:  
 \begin{proposition}\label{wcomul}
 $\Delta_{^*\!Q^m_2}(X) = \Delta_{^*\!Q^m_2}(1)\ (X^m_{(k_1 Sp_3 Sk_3 \rtimes F^{-1}(k_2 p_2) \rtimes Sk_1)_1} \ X^m_1 \ \otimes \ 
 X^m_3 \ X^m_{(k_1 Sp_3 Sk_3 \rtimes F^{-1}(k_2 p_2) \rtimes Sk_1)_2})$. 
 \end{proposition}
 Observe that the comultiplication formula involves $\Delta(1)$. Certain useful facts regarding $\Delta(1)$ are contained in the following lemma.
 \begin{lemma}\label{del}\cite[Proposition 4.12, Corollary 4.13]{Das2004}
  If $P = P^{N \subset M}$ denotes the subfactor planar algebra associated to the finite-index, reducible, depth two inclusion $N \subset M$
  of $II_1$-factors, then $\Delta_{P_2}(1) = f^1 \otimes Sf^2$ where $f^1 \otimes f^2$ is the unique symmetric separability element of 
  $P_{1, 2}$ and
  $\Delta_{P_2}$ denotes the comultiplication in the weak Hopf $C^*$-algebra $P_2$.
 \end{lemma}
 Let $U \otimes V$ denote the unique symmetric separability element of $^*\!Q^m_{1, 2}$ and let $U \otimes V ( \in {^*\!Q^m_{1, 2}}^{\otimes 2})$
 be given by
 \begin{align*}
 (\underbrace{\epsilon \rtimes 1 \rtimes \epsilon \rtimes \cdots \rtimes \epsilon}_{\text{$m+1$ \mbox{terms}}} \rtimes \ y^1 \rtimes g^2 \rtimes
 \cdots \rtimes g^{m-2}) \otimes (\underbrace{\epsilon \rtimes 1 \rtimes \epsilon \rtimes \cdots \rtimes \epsilon}_{\text{$m+1$ \mbox{terms}}} 
  \ \rtimes \tilde{y}^1 \rtimes \tilde{g}^2 \rtimes
 \cdots \rtimes \tilde{g}^{m-2})
  \end{align*}
or  
 \begin{align*}
 (\underbrace{1 \rtimes \epsilon \rtimes 1 \rtimes \cdots \rtimes \epsilon}_{\text{$m+1$ \mbox{terms}}} \rtimes \ y^1 \rtimes g^2 \rtimes
 \cdots \rtimes y^{m-2}) \otimes (\underbrace{1 \rtimes \epsilon \rtimes 1 \rtimes \cdots \rtimes \epsilon}_{\text{$m+1$ \mbox{terms}}} 
  \ \rtimes \tilde{y}^1 \rtimes \tilde{g}^2 \rtimes
 \cdots \rtimes \tilde{y}^{m-2})
 \end{align*}
 according as $m$ is even or odd.
 It then follows from Lemma \ref{del} and Lemma \ref{wanti} that $\Delta_{^*\!Q^m_2}(1) ( = U \otimes SV)$ equals
 \begin{align*}
 (\underbrace{\epsilon \rtimes 1 \rtimes \epsilon \rtimes \cdots \rtimes \epsilon}_{\text{$m+1$ \mbox{terms}}} \rtimes \ y^1 
 \rtimes g^2 \rtimes
 \cdots \rtimes g^{m-2}) \otimes (S\tilde{g}^{m-2} \rtimes \cdots \rtimes S\tilde{g}^2 \rtimes S\tilde{y}^1 \rtimes 
 \underbrace{\epsilon \rtimes 1 \rtimes \epsilon \rtimes \cdots \rtimes \epsilon}_{\text{$m+1$ \mbox{terms}}})   
 \end{align*} or
 \begin{align*}
 (\underbrace{1 \rtimes \epsilon \rtimes 1 \rtimes \cdots \rtimes \epsilon}_{\text{$m+1$ \mbox{terms}}} \rtimes \ y^1 
 \rtimes g^2 \rtimes
 \cdots \rtimes y^{m-2}) \otimes (S\tilde{y}^{m-2} \rtimes \cdots \rtimes S\tilde{g}^2 \rtimes S\tilde{y}^1 \rtimes 
 \underbrace{\epsilon \rtimes 1 \rtimes \epsilon \rtimes \cdots \rtimes 1}_{\text{$m+1$ \mbox{terms}}}),   
 \end{align*}
 according as $m$ is even or odd.
 Hence, it follows from Proposition \ref{wcomul} and Lemma \ref{del} that 
 \begin{align*}
 \Delta_{^*\!Q^m_2}(X) = U \ X^m_{(k_1 Sp_3 Sk_3 \rtimes F^{-1}(k_2 p_2) \rtimes Sp_1)_1} \ X^m_1 \ \ \otimes \ 
  S(V) \ X^m_3 \ X^m_{(k_1 Sp_3 Sk_3 \rtimes F^{-1}(k_2 p_2) \rtimes Sp_1)_2}.
 \end{align*}
 Using the comultiplication formula in $^*\!Q^2_2$ as given by Lemma \ref{sstructure} , we see that
 \begin{align*}
 \Delta_{^*\!Q^m_2}(X) = \ \delta \ &U \ X^m_{((k S\phi_2)_1 \ S(\phi_1p_2S\phi_3)_3 \ S(k S\phi_2)_3) \ \rtimes \ F^{-1}((k S\phi_2)_2 \ 
 (\phi_1p_2S\phi_3)_2) \ \rtimes \ S(\phi_1p_2S\phi_3)_1)} \ X^m_1 \ \otimes\\ 
  &S(V) \ X^m_3 \ X^m_{((\phi_4)_1 \ S(p_1)_3 S(\phi_4)_3 \ \rtimes \ F^{-1}((\phi_4)_2 (p_1)_2) \ \rtimes \ S(p_1)_1)}.
 \end{align*}
 Assume now that $m$ is even.
 A tedious computation using the multiplication rule in $^*\!Q^m_2$ shows that the formula for $\Delta_{^*\!Q^m_2}(X)$ is given by: 
 \begin{align}\label{eqq}
 \Delta_{^*\!Q^m_2}(X) &= \delta \ (\phi_4 Sp_2 S\phi_6)(S\tilde{y}^1_1) \ (f^1 \rtimes \cdots \rtimes x^{m-2} \rtimes 
 (k S\phi_2)_1 \ S(\phi_1p_3S\phi_3)_3 
  \ S(k S\phi_2)_3 \ \rtimes \\   
 & F^{-1}((k S\phi_2)_2 \ (\phi_1p_3S\phi_3)_2) \ \rtimes \ S(\phi_1p_3S\phi_3)_1 
  \rtimes y^1 
 \rtimes g^2 \rtimes \cdots \rtimes g^{m-2}) \ \otimes \nonumber \\
 & \ (S\tilde{g}^{m-2} \rtimes \cdots \rtimes S\tilde{g}^2 \rtimes S\tilde{y}^1_2  
   \rtimes  
 (\phi_5)_1 \ S(p_1)_3 S(\phi_5)_3 \ \rtimes \ F^{-1}((\phi_5)_2 (p_1)_2) \ \rtimes \ S(p_1)_1 \ \rtimes \nonumber \\
 & x^{m+2} \rtimes f^{m+3} \rtimes \cdots \rtimes f^{2m-1}) \nonumber.
 \end{align}
Similarly, when $m$ is odd, one can show that the formula for $\Delta_{^*\!Q^m_2}(X)$ is given by:
\begin{align}\label{eqqq}
 \Delta_{^*\!Q^m_2}(X) &= \delta \ (\phi_4 Sp_2 S\phi_6)(S\tilde{y}^1_1) \ (x^1 \rtimes \cdots \rtimes x^{m-2} \rtimes
 (k S\phi_2)_1 \ S(\phi_1p_3S\phi_3)_3 
  \ S(k S\phi_2)_3 \ \rtimes \\   
 & F^{-1}((k S\phi_2)_2 \ (\phi_1p_3S\phi_3)_2) \ \rtimes \ S(\phi_1p_3S\phi_3)_1 
  \rtimes y^1 \rtimes \cdots \rtimes y^{m-2}) \ \otimes \nonumber \\
  & \ (S \tilde{y}^{m-2} \rtimes \cdots \rtimes S\tilde{y}^1_2 
   \rtimes  
 (\phi_5)_1 \ S(p_1)_3 S(\phi_5)_3 \ \rtimes \ F^{-1}((\phi_5)_2 (p_1)_2) \ \rtimes \ S(p_1)_1 \ \rtimes \nonumber \\
 & x^{m+2} \rtimes f^{m+3} \rtimes \cdots \rtimes x^{2m-1}) \nonumber.
 \end{align}
 For each integer $m > 2$, let $K_m$ denote the vector space $A(H)_{m-2}^{op} \otimes D(H)^{*op} \otimes A(H)_{m-2}^{op}$ or  
$A(H^*)_{m-2}^{op} \otimes D(H)^{*op} \otimes A(H)_{m-2}^{op}$ according as $m$ is odd or even.
Consider the linear isomorphism $\psi$ of $K_m$ onto $^*\!Q^m_2$ given by
  \begin{align*}
   a \otimes (g \otimes f) \otimes b \mapsto a \rtimes f_1Sg_3Sf_3 \rtimes F^{-1}(f_2 g_2) \rtimes Sg_1 \rtimes b.
  \end{align*}
  We make $K_m$ into a 
  weak Hopf $C^*$-algebra by transporting the structure maps on $^*\!Q^m_2$ to $K_m$ using this linear isomorphism. Thus, by 
  construction, $K_m$ is isomorphic to $^*\!Q^m_2$ as weak Hopf $C^*$-algebras.
  The next theorem, which is the main result of this section, explicitly describes the structure maps of $K_m$.
  \begin{theorem}\label{WHA}
  For each $m > 2$, $K_m$ is a weak Hopf $C^*$- algebra with the structure maps given by the following formulae.
\begin{align*}
  &\mbox{Multiplication:} \ \ (a \otimes (g \otimes f) \otimes b) ( \tilde{a} \otimes (\tilde{g} \otimes \tilde{f}) \otimes \tilde{b}) = 
  (\tilde{f_1} S \tilde{g_2} S \tilde{f_3}.a) \tilde{a} \otimes (g_2 \otimes f)(\tilde{g_1} \otimes \tilde{f_2}) \otimes 
  b \rho_{g_1}(\tilde{b}),\\
  &\mbox{Comultiplication:} \ \ \Delta(a \otimes (g \otimes f) \otimes b) = \delta (a \otimes (\phi_1 g_3 S\phi_3 \otimes fS\phi_2) \otimes u)
  \otimes ((\phi_4 Sg_2 S\phi_6). v^{\prime} \otimes (g_1 \otimes \phi_5) \otimes b), \\
   & \mbox{Counit:} \ \varepsilon(a \otimes (g \otimes f) \otimes b) = 
       \delta^{m-2} f(h) Z_{F^{(m-1)}}^{P(H^m)}(a \rtimes Sg \rtimes b),\\
  &\mbox{Antipode:} \ \ S(a \otimes (g \otimes f) \otimes b) = b^{\prime} \otimes S_{D(H)^{*op}}(g \otimes f) \otimes a^{\prime},
   \end{align*} 
   with $a \otimes (g \otimes f) \otimes b$, $\tilde{a} \otimes (\tilde{g} \otimes \tilde{f}) \otimes \tilde{b}$ being
   elements of $K_m$, where
   \begin{itemize}
   \item for any $k \in H^*$ and $\underbrace{(\cdots \rtimes f \rtimes x)}_{\text{$m-2$ factors}}$ in $A(H)^{op}_{m-2} \ \mbox{or} \ A(H^*)^{op}_{m-2}$ according as $m$ is odd or 
   even, $k. (\cdots \rtimes f \rtimes x) := k(x_2)(\cdots \rtimes f \rtimes x_1)$,
   \item $\rho$ denotes the algebra action of $H^*$ on $A(H)^{op}_{m-2}$ defined for $k \in H^*$ and $x \rtimes p \rtimes \cdots \in 
   A(H)^{op}_{m-2}$
   by  $\rho_k(x \rtimes p \rtimes \cdots) = 
   k(Sx_1) x_2 \rtimes p \rtimes \cdots$, 
   \item $u \otimes v$ denotes the unique symmetric separability element of $A(H)_{m-2}$, 
   \item for any positive integer $k$, $F^{(k)}$ denotes the $k$-cabling of the tangle $F$ (see Figure \ref{fig:pic47}) and 
   \item for any $X$ in $H_{[i, j]}$, the symbol $X^{\prime}$ denotes the element as defined in $\S 1.1$ preceding Lemma \ref{anti1}. 
   \end{itemize}      
 \end{theorem}
 \begin{proof}
 We first verify the formula for antipode. Assume without loss of generality that $m$ is even. Let $X \in K_m$ be given by
 \begin{align*}
  (f^1 \rtimes x^2 \rtimes \cdots \rtimes x^{m-2}) \otimes (p \otimes k) \otimes (x^{m+2} \rtimes \cdots \rtimes f^{2m-1}) 
 \end{align*}
 so that 
 \begin{align*}
  \psi(X) = f^1 \rtimes x^2 \rtimes \cdots \rtimes x^{m-2} \rtimes k_1 Sp_3 Sk_3 \rtimes
  F^{-1}(k_2 p_2) \rtimes Sp_1 \rtimes x^{m+2} \rtimes \cdots \rtimes f^{2m-1}.
 \end{align*}
 By Lemma \ref{wanti},
 \begin{align*}
  S(\psi(X)) &= Sf^{2m-1} \rtimes \cdots \rtimes Sx^{m+2} \rtimes p_1 \rtimes F^{-1}S(k_2 p_2) \rtimes S(k_1 Sp_3 Sk_3) \rtimes  
  Sx^{m-2} \rtimes \cdots \rtimes Sx^2 \rtimes Sf^1\\
  &= (x^{m+2} \rtimes \cdots \rtimes f^{2m-1})^{\prime} \rtimes p_1 \rtimes F^{-1}S(k_2 p_2) \rtimes S(k_1 Sp_3 Sk_3) \rtimes
  (f^1 \rtimes x^2 \rtimes \cdots \rtimes x^{m-2})^{\prime}.
 \end{align*}
Finally, it follows from the formula for antipode in $^*\!Q_2$ as given by Lemma \ref{sstructure} and Remark \ref{Dr} that 
\begin{align*}
  \psi^{-1}(S(\psi(X))) = (x^{m+2} \rtimes \cdots \rtimes f^{2m-1})^{\prime} \otimes S_{D(H)^{* op}}(p \otimes k) \otimes
  (f^1 \rtimes x^2 \rtimes \cdots \rtimes x^{m-2})^{\prime}.
\end{align*}
Thus the formula for antipode is verified. In a similar way, using the comultiplication formula in $^*\!Q^m_2$ as given by 
\eqref{eqq} or \eqref{eqqq} according as $m$ is even or odd, the comultiplication formula in $K_m$ can easily be verified. 
The verification of the multiplication
 and counit formula in $K_m$ involves tedious computation and we leave these verifications for the reader.   
 \end{proof}



%
%

 \begin{remark}
  It follows immediately from the formula for antipode as given in Theorem \ref{WHA} that each $K_m$ ($m > 2$) has involutive anipode i.e., 
  square of the antipode is the identity and consequently, each $K_m$ is a weak Kac algebra.
 \end{remark}
 \begin{remark}
  The sole importance of Theorem \ref{WHA} lies in the fact that it constructs a family of weak Kac algebras out of a given finite-dimensional 
  Kac algebra.
 \end{remark}

 \section*{acknowledgement}
 The author sincerely thanks Prof. Vijay Kodiyalam for a careful reading of the manuscript and suggestions for improvement and also for  
 his support, constructive comments and several fruitful discussions during the whole project. He also thanks Prof. David Evans for his
 support in making it possible to attend the inspiring programme ``Operator Algebras: Subfactors and their Applications'' held at INI, 2017 
 while this paper was being written up. The author was supported by 
 National Board of Higher Mathematics, India.

\end{document}